%% file: Phasefield.tex
\documentclass[review,hidelinks,onefignum,onetabnum]{siamart220329}
\usepackage[top=2in,bottom=1in, left=1in, right=1in]{geometry}
\usepackage{subfigure,multirow}
\usepackage{amsmath,amssymb}
\usepackage{graphicx}
\usepackage{float}
\usepackage{booktabs}
\usepackage{multirow}
\usepackage[justification=centering]{caption}

\input{ex_shared}

\ifpdf
\hypersetup{
  pdftitle={A thermodynamically consistent phase-field model for mass transport with interfacial reaction and deformation},
  pdfauthor={Zhaoyang Wang, Huaxiong Huang, Ping Lin and Shixin Xu}
}
\fi

\begin{document}

\maketitle

\begin{abstract}
In this paper, a thermodynamically consistent phase-field model is proposed to describe the mass transport and reaction processes of multiple species in a fluid. A key feature of this model is that reactions between different species occur only at the interface, and may induce deformation of the interface. For the governing equations derived based on the energy variational method, we propose a structure-preserving numerical scheme that satisfies the mass conservation and energy dissipation laws at the discrete level. Furthermore, we carry out a rigorous error analysis of the time-discrete scheme for a simplified case. A series of numerical experiments are conducted to validate the effectiveness of the model as well as the accuracy and stability of the scheme. In particular, we simulate microvessels with straight and bifurcated structures to illustrate the risk of microaneurysm formation.
\end{abstract}



\section{Introduction}
\label{section1}
Transmembrane transport of substances and reactions occurring at the membrane (interface) are common processes in both industrial and biological systems. In chemical engineering, membrane separation technologies \cite{strathmann1981membrane} and transmembrane transport \cite{bernardo2020recent} in reaction processes are essential components of many industrial operations. The transmembrane transport of species and the reactions occurring at the vessel wall \cite{claesson2021permeability} in the human body are closely associated with health. A typical example is the formation of retinal microaneurysms \cite{hirai2007retinopathy}. Under hyperglycemic conditions, structural proteins adhered to the vascular wall react with glucose, leading to the formation of advanced glycation end-products (AGEs). The accumulation of AGEs induces vascular wall stiffening and loss of elasticity, eventually contributing to the formation of retinal microaneurysms \cite{wautier2015advanced}. 

These processes, characterized by the dynamic interaction of mass transport, interfacial reactions, and the potentially resulting mechanical deformations in fluid environments, present significant challenges for mathematical modeling and numerical simulation. Conventional reaction-diffusion models usually consider reactions to occur uniformly within a domain, and rarely pay attention to the reaction characteristics and topological changes at the interface or membrane. However, changes in the shape of the solid interface caused by the reaction may affect the dynamic behavior of other coupled fields \cite{yang2016mathematical}. 
For such fluid-structure interaction (FSI) problems, various methods have been developed, including the level set method \cite{osher2004level, wang2003level}, the immerse boundary method \cite{peskin2002immersed}, the front-tracking method \cite{unverdi1992front, tryggvason2001front}, and the phase-field (diffuse-interface) method \cite{du2004phase, guo2015thermodynamically}. 

In particular, the phase-field method describes the membrane as a thin diffusive layer, which serves as an interface separating two incompressible fluids. The interface can be represented by the zero level-set of an order parameter. This parameter takes the values of $1$ and $-1$, corresponding to the bulk fluids inside and outside the interface, respectively, thus eliminating the need to explicitly track the interface position. Due to its ability to easily deal with the evolution of the interface over time, the phase-field method has been successfully applied to various interface deformation problems induced by fluid dynamics, such as vesicle motion \cite{shen2022energy}, tumor growth \cite{oden2010general}, and moving contact line dynamics \cite{zhu2019thermodynamically}. Recently, the phase-field approach has been employed to model transmembrane transport processes. Qin et al. \cite{qin2022phase} proposed a thermodynamically consistent phase-field model to describe mass transfer across the permeable moving interface, where the diffusion flux across the membrane is affected by the conductivity and concentration difference on each side. A multi-component solute diffusion model was developed in \cite{kou2023thermodynamically}, which involves the crossing influences between different solutes. To the best of our knowledge, there is no work that considers reactions occurring on the membrane that cause its deformation during mass transport, despite this being quite common in the field of biology.

In this paper, we first derive a thermodynamically consistent phase-field model, which incorporates mass transport, interfacial reactions and the associated interface deformation under flow conditions into a unified framework. Unlike previous models (e.g., surfactant models \cite{laradji1992effect, yang2017linear}) of interface and substances interactions, our model explicitly incorporates reactions confined to the interface that drive morphological changes. By defining energy and dissipation functionals based on the physical problem considered, and combining them with the kinematic conservation laws, we use the energy variational method to derive this model that satisfies the mass conservation and energy dissipation laws. Furthermore, we carefully design an efficient structure-preserving numerical scheme for the proposed nonlinear system, ensuring that it maintains mass conservation and energy dissipation laws at the discrete level. Given the scarcity of numerical analysis work on such problems, for the simplified case of the model, which consists of a coupled system of the Cahn-Hilliard-Navier-Stokes equations and the reaction-diffusion equation, we establish rigorous error estimates for the time-discrete scheme in $L^{\infty}(0,T;H^1(\Omega))\cap L^{2}(0,T;H^2(\Omega))$ for the velocity, $L^{2}(0,T;H^1(\Omega))$ for the pressure, $L^{\infty}(0,T;H^1(\Omega))$ for the phase function and $L^{\infty}(0,T;L^2(\Omega))\cap L^{2}(0,T;H^1(\Omega))$ for the concentration. 

To demonstrate the practical application of the proposed model, we simulate the formation of microaneurysms in a straight vessel and further illustrate the risk of microaneurysm formation in a bifurcated vessel. We would like to point out that this work provides a theoretical and numerical framework for studying FSI problems affected by interfacial reactions and deformations, which can be easily applied to study other fluid-related interfacial dynamics problems, such as atherosclerotic plaque growth \cite{wang2023fast} and high-temperature corrosion of metals \cite{jin2022corrosion}.

\subsection*{Outline}
The rest of this manuscript is organized as follows. In Section \ref{section2}, we introduce the energy functionals and derive the phase-field model with reactions at the interface. The mass conservation and energy dissipation laws of the model are given. In Section \ref{section3}, we construct a time-discrete scheme and the corresponding fully discrete finite element scheme, and prove their mass conservation and energy stability. In Section \ref{section4}, we carry out error estimates for the time-discrete scheme. Numerical results are presented in Section \ref{section5} to validate the proposed model and numerical method. The simulation of a microaneurysm is shown to demonstrate the practical application of the model. Some conclusions and remarks are given in final Section \ref{section6}.

\subsection*{Notation}
For domain $\Omega$ and $1\leq p\leq\infty$, we use the standard notation for the Banach space $L^{p}(\Omega)$ and the Sobolev space $W^{k,p}(\Omega)$ or $H^p(\Omega)$. The symbol $(\cdot,\cdot)$ indicates the standard scalar product in $L^2$. Throughout this paper, the letter $C$ denotes a generic positive constant, with or without subscript, its value may change from one line of an estimate to the next. We will write the dependence of the constant on parameters explicitly if it is essential.

\section{Phase-field model with reactions at the interface}
\label{section2}
In this section, we use the energy variational method to derive a thermodynamically
consistent phase-field model with reactions at the phase interface in incompressible fluids.

\subsection{Model derivation}
Let the problem domain be $\Omega$ with boundary $\partial \Omega$. We use 
$\phi\in [-1,1]$ as the phase-field label function, where $\{\phi=1\}$ represents one phase, $\{\phi=-1\}$ represents the other phase, and $\{\phi=0\}$ denotes the interface or membrane existing between the two phases. 

We assume that there are three species $Y_1$, $Y_2$ and $Y_3$ in $\Omega$, where $Y_1$ can be distributed throughout the entire domain and undergo transmembrane transport, while $Y_2$ and $Y_3$ are confined to the interface. A reaction $a_1Y_1+a_2Y_2 \rightleftharpoons a_3Y_3$ only occurs at the interface, where $a_1$, $a_2$ and $a_3$ are stoichiometric coefficients and satisfy the conservation relation $a_1+a_2=a_3$. Let $c_i \ (i=1,2,3)$ denote the concentration of each species, and define $\mathcal{R}\geq 0$ as the net reaction rate function.

We start with the model derivation based on the incompressibility condition of the fluids and the mass conservation relations:
\begin{subequations}\label{the2.1}
	\begin{align}
&\rho\frac{\textbf{D}\textbf{\textit{u}}}{\textbf{D}t}=\nabla \cdot \sigma_{\nu}+\nabla \cdot \sigma_{\phi}, \quad \nabla \cdot \textbf{\textit{u}}=0, \label{the2.1a}\\
&\frac{\textbf{D}\phi}{\textbf{D}t}=-\nabla \cdot \textbf{\textit{j}}_{\phi}, \label{the2.1b}\\
&\frac{\textbf{D}c_1}{\textbf{D}t}=-\nabla \cdot \textbf{\textit{j}}_{1}-a_1\mathcal{R}, \label{the2.1c}\\
&\frac{\textbf{D}c_2}{\textbf{D}t}=-\nabla \cdot \textbf{\textit{j}}_{2}-a_2\mathcal{R}, \label{the2.1d}\\
&\frac{\textbf{D}c_3}{\textbf{D}t}=-\nabla \cdot \textbf{\textit{j}}_{3}+a_3\mathcal{R}. \label{the2.1e}
\end{align}
\end{subequations}
Here, $\rho$ is the density of the fluid, $\textbf{\textit{u}}$ is the fluid velocity, $\frac{\textbf{D}}{\textbf{D} t}=\frac{\partial}{\partial t}+\textbf{\textit{u}}\cdot\nabla$ is the material derivative, and $\textbf{\textit{j}}$ is the mass flux that will be determined later. $\sigma_{\nu}$ and $\sigma_{\phi}$ are stresses induced by the viscosity of the fluid and by phase interface, respectively. For variables $\phi$, $c_1$, $c_2$, $c_3$, $\textbf{\textit{j}}_{\phi}$, $\textbf{\textit{j}}_{1}$, $\textbf{\textit{j}}_{2}$, the homogeneous Neumann boundary conditions are used. For fluid velocity $\textbf{\textit{u}}$, we use the no-slip boundary condition.

We define kinetic energy $I_1$, phase mixing energy $I_2$, entropy energy $I_3$ of $Y_1$, and mixing energies $I_4$ and $I_5$ of $Y_2$ and $Y_3$ in the domain $\Omega$ as 
\begin{equation}\label{the2.2}
	\begin{split}
&I_1=\frac{1}{2}\int_\Omega \rho \textbf{\textit{u}}^2 d\textbf{x}, \ I_2=\int_{\Omega} \lambda(c_3)\left(F(\phi)+\frac{\epsilon^2}{2}|\nabla\phi|^2\right)d\textbf{x}, \ I_3= \int_\Omega \mathcal{G}c_1\left(\text{ln}\left(\frac{c_1}{c_0}\right)-1\right) d\textbf{x}, \\
&I_4=\int_\Omega \mathcal{G}c_2\left(\text{ln}\left(\frac{c_2}{c_0}\right)-1\right)d\textbf{x}+M\int_\Omega \mathcal{G}c_2\phi^2 d\textbf{x}-\frac{N}{4}\int_{\Omega}\mathcal{G} c_2\left(\phi^2-1\right)^2d\textbf{x},\\
&I_5=\int_\Omega \mathcal{G}c_3\left(\text{ln}\left(\frac{c_3}{c_0}\right)-1\right)d\textbf{x}+M\int_\Omega \mathcal{G}c_3\phi^2 d\textbf{x}-\frac{N}{4}\int_{\Omega}\mathcal{G} c_3\left(\phi^2-1\right)^2d\textbf{x},
	\end{split}
\end{equation}
where $\lambda>0$ is the mixing energy density function related to the concentration of generated $Y_3$, and its variation implies a change in interfacial tension and may cause interface deformation. $F(\phi)=\frac{1}{4}\left(\phi^2-1\right)^2$ is the double well potential, $c_0$ is the reference concentration, $\epsilon$ is the thickness of the interface. $\mathcal{G}$ is the energy scale in thermodynamics, and the unit is $\text{J/mol}$. 

We would like to note that in order to ensure $Y_2$ and $Y_3$ only exist at the interface, local coupling terms $c_i\phi^2$ and $-c_i\left(\phi^2-1\right)^2$ $(i=2,3)$ are added to $I_4$ and $I_5$, inspired by the idea of constructing a surfactant model \cite{engblom2013diffuse}. The role of terms $M\mathcal{G}c_2\phi^2$ and $M\mathcal{G}c_3\phi^2$ is to penalize the free $Y_2$ and $Y_3$ in the bulk phase with the penalty parameter $M\geq 0$. $-\frac{N}{4}\mathcal{G} c_2\left(\phi^2-1\right)^2$ and $-\frac{N}{4}\mathcal{G} c_3\left(\phi^2-1\right)^2$ enhance the adhesion of $Y_2$ and $Y_3$ to the interface with the parameter $N\geq 0$. 

For the total energy $E_{total}=\sum_{i=0}^{5} I_i$, we can define the chemical potentials:
\begin{equation}\label{the2.3}
\begin{split}
&\mu_\phi=\frac{\delta E_{total}}{\delta \phi}=-\epsilon^2\nabla \cdot \left(\lambda\nabla \phi\right)+\lambda F'(\phi)+2M\mathcal{G}c_2\phi-N\mathcal{G}c_2(\phi^3-\phi)+2M\mathcal{G}c_3\phi-N\mathcal{G}c_3(\phi^3-\phi),
\end{split}
\end{equation}
\begin{equation}\label{the2.4}
	\begin{split}
\mu_1&=\frac{\delta E_{total}}{\delta c_1}=\mathcal{G}\text{ln}\left(\frac{c_1}{c_0}\right),
	\end{split}
\end{equation}
\begin{equation}\label{the2.5}
\begin{split}
\mu_2&=\frac{\delta E_{total}}{\delta c_2}=\mathcal{G}\text{ln}\left(\frac{c_2}{c_0}\right)+M\mathcal{G}\phi^2-\frac{N}{4}\mathcal{G}(\phi^2-1)^2,
\end{split}
\end{equation}
\begin{equation}\label{the2.6}
	\begin{split}
	\mu_3&=\frac{\delta E_{total}}{\delta c_3}=\frac{\partial \lambda}{\partial c_3}\left(F(\phi)+\frac{\epsilon^2}{2}\left|\nabla \phi\right|^2\right)+\mathcal{G}\text{ln}\left(\frac{c_3}{c_0}\right)+M\mathcal{G}\phi^2-\frac{N}{4}\mathcal{G}(\phi^2-1)^2.
	\end{split}
\end{equation}

The dissipation functional of the system consists of the dissipation due to fluid friction, mixing of two phases in bulk, mixing and reaction of species
\begin{equation}\label{the2.7}
\begin{split}
\Delta&=\int_\Omega 2\nu |D_\nu|^2d\textbf{x}+\int_\Omega \mathcal{M} |\nabla\mu_\phi|^2d\textbf{x}+\int_\Omega \frac{D_1c_1}{\mathcal{G}} |\nabla\mu_1|^2d\textbf{x}+\int_\Omega \frac{D_2c_2}{\mathcal{G}} |\nabla\mu_2|^2d\textbf{x}\\
&+\int_\Omega \frac{D_3c_3}{\mathcal{G}} |\nabla\mu_3|^2d\textbf{x}+\int_\Omega \frac{P(\phi)}{\mathcal{G}}\left(a_1\mu_1+a_2\mu_2-a_3\mu_3\right)^2d\textbf{x},
\end{split}
\end{equation}
where $\nu$ is the fluid density, $D_\nu(\textbf{\textit{u}})=\frac{1}{2}(\nabla \textbf{\textit{u}}+\nabla \textbf{\textit{u}}^{\dag})$ is the strain rate, $\mathcal{M}$ is the mobility rate, $D_i \ (i=1,2,3)$ are diffusion coefficients. In particular, considering the restricted diffusion of $Y_1$ across the membrane due to the potential permeability of the interface, $D_1$ is defined as \cite{qin2022phase}:
\begin{equation}\label{the2.8}
	\begin{split}
\frac{1}{D_1}= \frac{1+\phi}{2D_1^{+}}+\frac{1-\phi}{2D_1^{-}}+\frac{(1-\phi^2)^2}{q\epsilon },
	\end{split}
\end{equation}
where $D_1^+$ and $D_1^-$ are the diffusion coefficients of the two phases inside and outside the interface, respectively. $q=\frac{dQ(c_1)}{dc_1}$, where $Q$ is a function that describes the membrane permeability \cite{xu2018osmosis}.

We would like to point out that the setting of the reaction dissipation (the last term in (\ref{the2.7})) is motivated by the linear phenomenological constitutive laws for chemical reactions (see. e.g., \cite{de2013non}). This form of dissipation has also been applied in the tumor growth model \cite{hawkins2012numerical, wang2025stability}. 

According to the law of energy dissipation, the rate of change of total energy
equals the dissipation $\frac{d}{dt}E_{total}=-\Delta\leq 0$. By taking the time derivative of each term in the energy $E_{total}$ sequentially, we obtain
\begin{equation}\label{the2.9}
	\begin{split}
\frac{dI_1}{dt}=-\int_{\Omega} \sigma_{\nu}:\nabla\textbf{\textit{u}}d\textbf{x}-\int_\Omega \sigma_{\phi}: \nabla\textbf{\textit{u}}   d\textbf{x}-\int_{\Omega}p\nabla\cdot \textbf{\textit{u}}d\textbf{x},
	\end{split}
\end{equation}

\begin{equation}\label{the2.10}
\begin{split}
\frac{dI_2}{dt}&=\int_{\Omega}  \lambda F'(\phi)\frac{\partial \phi}{\partial t}d\textbf{x}-\int_\Omega\epsilon^2\nabla\cdot (\lambda\nabla\phi)\frac{\partial \phi}{\partial t} d\textbf{x}+\int_\Omega\frac{\partial \lambda}{\partial c_3}\frac{\partial c_3}{\partial t}\left(F(\phi)+\frac{\epsilon^2}{2}|\nabla\phi|^2\right)d\textbf{x}\\
&=\int_{\Omega}  \mu_\phi \left(\frac{\partial \phi}{\partial t}+\textbf{\textit{u}}\cdot \nabla \phi\right)d\textbf{x}-\int_\Omega \mu_\phi \textbf{\textit{u}}\cdot \nabla\phi d\textbf{x}+\int_\Omega\frac{\partial \lambda}{\partial c_3}\frac{\partial c_3}{\partial t}\left(F(\phi)+\frac{\epsilon^2}{2}|\nabla\phi|^2\right)d\textbf{x}\\
&-\int_\Omega 2M\mathcal{G}c_2\phi\frac{\partial \phi}{\partial t}d\textbf{x}+\int_\Omega N\mathcal{G}c_2(\phi^3-\phi)\frac{\partial \phi}{\partial t}d\textbf{x}-\int_\Omega 2M\mathcal{G}c_3\phi\frac{\partial \phi}{\partial t}d\textbf{x}+\int_\Omega N\mathcal{G}c_3(\phi^3-\phi)\frac{\partial \phi}{\partial t}d\textbf{x}\\
&=\int_\Omega \nabla\mu_\phi \cdot \textbf{\textit{j}}_\phi d\textbf{x}+\int_\Omega \epsilon^2\left(\nabla\cdot(\lambda\nabla\phi\otimes\nabla\phi)\right)\cdot \textbf{\textit{u}}d\textbf{x}-\int_\Omega\lambda\left(\nabla\left(\frac{\epsilon^2}{2}\left|\nabla\phi\right|^2+F(\phi)\right)\right)\cdot \textbf{\textit{u}}d\textbf{x}\\
&+\int_\Omega\frac{\partial \lambda}{\partial c_3}\frac{\partial c_3}{\partial t}\left(F(\phi)+\frac{\epsilon^2}{2}|\nabla\phi|^2\right)d\textbf{x}\\
&+\int_\Omega \mathcal{G}c_2\left(\frac{N}{4}\nabla\left(\left(\phi^2-1\right)^2\right)-M\nabla\phi^2\right)\cdot \textbf{\textit{u}}d\textbf{x}+\int_\Omega\mathcal{G} c_3\left(\frac{N}{4}\nabla\left(\left(\phi^2-1\right)^2\right)-M\nabla\phi^2\right)\cdot \textbf{\textit{u}}d\textbf{x}\\
&-\int_\Omega 2M\mathcal{G}c_2\phi\frac{\partial \phi}{\partial t}d\textbf{x}+\int_\Omega N\mathcal{G}c_2(\phi^3-\phi)\frac{\partial \phi}{\partial t}d\textbf{x}-\int_\Omega 2M\mathcal{G}c_3\phi\frac{\partial \phi}{\partial t}d\textbf{x}+\int_\Omega N\mathcal{G}c_3(\phi^3-\phi)\frac{\partial \phi}{\partial t}d\textbf{x},
\end{split}
\end{equation}
where we use the fact that $\nabla\cdot(\lambda\nabla\phi\otimes\nabla\phi)=\nabla\lambda\left|\nabla\phi\right|^2+\lambda\Delta \phi\nabla\phi+\frac{\lambda}{2}\nabla\left|\nabla\phi\right|^2$.

\begin{equation}\label{the2.11}
\begin{split}
\frac{dI_3}{dt}&=\int_\Omega\mu_1\frac{\partial c_1}{\partial t}d\textbf{x}=\int_{\Omega} \nabla \mu_1\cdot\textbf{\textit{j}}_{1} d\textbf{x}-a_1\int_{\Omega}\mu_1\mathcal{R}  d\textbf{x}
\end{split}
\end{equation}

\begin{equation}\label{the2.12}
\begin{split}
\frac{dI_4}{dt}&=\int_{\Omega} \mathcal{G}\text{ln}\left(\frac{c_2}{c_0}\right)\frac{\partial c_2}{\partial t}d\textbf{x}+M\int_{\Omega}\frac{\partial c_2}{\partial t}\mathcal{G}\phi^2d\textbf{x}+2M\int_{\Omega} \mathcal{G}c_2\phi\frac{\partial \phi}{\partial t}d\textbf{x}\\
&-\frac{N}{4}\int_\Omega \frac{\partial c_2}{\partial t}\mathcal{G}(\phi^2-1)^2d\textbf{x}-\int_\Omega N\mathcal{G}c_2(\phi^3-\phi)\frac{\partial \phi}{\partial t}d\textbf{x}\\
&=\int_\Omega \nabla\mu_2\cdot \textbf{\textit{j}}_2d\textbf{x}-a_2\int_\Omega \mu_2\mathcal{R}d\textbf{x}+\int_\Omega\mathcal{G}\left(M\phi^2-\frac{N}{4}(\phi^2-1)^2\right)\textbf{\textit{u}}\cdot \nabla c_2d\textbf{x}\\
&+2M\int_{\Omega} \mathcal{G}c_2\phi\frac{\partial \phi}{\partial t}d\textbf{x}-\int_\Omega N\mathcal{G}c_2(\phi^3-\phi)\frac{\partial \phi}{\partial t}d\textbf{x},
\end{split}
\end{equation}

\begin{equation}\label{the2.13}
\begin{split}
\frac{dI_5}{dt}&=\int_{\Omega}\mathcal{G}\text{ln}\left(\frac{c_3}{c_0}\right)\frac{\partial c_3}{\partial t}d\textbf{x}+M\int_{\Omega}\frac{\partial c_3}{\partial t}\mathcal{G}\phi^2d\textbf{x}+2M\int_{\Omega} \mathcal{G}c_3\phi\frac{\partial \phi}{\partial t}d\textbf{x}\\
&-\frac{N}{4}\int_\Omega \frac{\partial c_3}{\partial t}\mathcal{G}(\phi^2-1)^2d\textbf{x}-\int_\Omega N\mathcal{G}c_3(\phi^3-\phi)\frac{\partial \phi}{\partial t}d\textbf{x}\\
&=\int_\Omega \nabla\mu_3\cdot \textbf{\textit{j}}_3d\textbf{x}+a_3\int_\Omega \mu_3\mathcal{R}d\textbf{x}-\int_\Omega\frac{\partial \lambda}{\partial c_3}\frac{\partial c_3}{\partial t}\left(F(\phi)+\frac{\epsilon^2}{2}\left|\nabla \phi\right|^2\right)d\textbf{x}\\
&-\int_\Omega \frac{\partial \lambda}{\partial c_3}\left(F(\phi)+\frac{\epsilon^2}{2}|\nabla \phi|^2\right)\textbf{\textit{u}}\cdot \nabla c_3d\textbf{x}+\int_\Omega\mathcal{G}\left(M\phi^2-\frac{N}{4}(\phi^2-1)^2\right)\textbf{\textit{u}}\cdot \nabla c_3d\textbf{x}\\
&+2M\int_{\Omega} \mathcal{G}c_3\phi\frac{\partial \phi}{\partial t}d\textbf{x}-\int_\Omega N\mathcal{G}c_3(\phi^3-\phi)\frac{\partial \phi}{\partial t}d\textbf{x}.
\end{split}
\end{equation}

Summing them up gives
\begin{equation}\label{the2.14}
	\begin{split}
\frac{dE_{\text{total}}}{dt}&=-\int_{\Omega} \sigma_{\nu}:\nabla\textbf{\textit{u}}d\textbf{x}-\int_\Omega \sigma_{\phi}: \nabla\textbf{\textit{u}}   d\textbf{x}-\int_{\Omega}p\nabla\cdot \textbf{\textit{u}}d\textbf{x}+\int_{\Omega} \nabla \mu_1\cdot\textbf{\textit{j}}_{1} d\textbf{x}-a_1\int_{\Omega}\mu_1\mathcal{R}  d\textbf{x}\\
&+\int_\Omega \nabla\mu_\phi \cdot \textbf{\textit{j}}_\phi d\textbf{x}-\int_\Omega \epsilon^2\left(\lambda\nabla\phi\otimes\nabla\phi\right): \nabla\textbf{\textit{u}}d\textbf{x}+\int_\Omega \nabla\mu_2\cdot \textbf{\textit{j}}_2d\textbf{x}-a_2\int_\Omega \mu_2\mathcal{R}d\textbf{x}\\
&+\int_\Omega \nabla\mu_3\cdot \textbf{\textit{j}}_3d\textbf{x}+a_3\int_\Omega \mu_3\mathcal{R}d\textbf{x}=-\Delta\leq 0.
	\end{split}
\end{equation}
 
By comparing with the predefined dissipation functional (\ref{the2.7}) yield
\begin{subequations}\label{the2.15}
\begin{align}
&\sigma_\nu=2\nu D_\nu-p\textbf{I}, \label{the2.15a}\\
&\sigma_\phi=-\lambda\epsilon^2(\nabla\phi\otimes\nabla\phi), \label{the2.15b}\\
&\textbf{\textit{j}}_\phi=-\mathcal{M}\nabla\mu_\phi, \label{the2.15c}\\
&\textbf{\textit{j}}_1=-\frac{D_1c_1}{\mathcal{G}}\nabla\mu_{1}, \label{the2.15d}\\
&\textbf{\textit{j}}_2=-\frac{D_2c_2}{\mathcal{G}}\nabla\mu_{2}, \label{the2.15e}\\
&\textbf{\textit{j}}_3=-\frac{D_3c_3}{\mathcal{G}}\nabla\mu_{3}, \label{the2.15f} \\
&\mathcal{R}=P(\phi)\mathcal{G}^{-1}(a_1\mu_1+a_2\mu_2-a_3\mu_3). \label{the2.15g}
\end{align}
\end{subequations}

Since the reaction only occurs at the interface, the following is a proper choice of the proliferation function $P(\phi)$:
\begin{equation}\label{the2.16}
\begin{split}
P(\phi)=k(\phi^2-1)^2 \ (-1\leq \phi \leq 1, \ k\geq 0), \quad P(\phi)=0 \ (\text{otherwise}).
\end{split}
\end{equation}

The proposed model can be summarized as follows:
\begin{subequations}\label{the2.17}
	\begin{align}
&\rho\left(\frac{\partial \textbf{\textit{u}}}{\partial t}+\left(\textbf{\textit{u}}\cdot \nabla\right)\textbf{\textit{u}}\right)=\nu\Delta \textbf{\textit{u}}-\nabla p-\epsilon^2\nabla\cdot(\lambda(\nabla\phi\otimes\nabla\phi)), \label{the2.17a}\\
&\nabla\cdot \textbf{\textit{u}}=0, \label{the2.17b}
\end{align}
\end{subequations}

\begin{subequations}\label{the2.18}
	\begin{align}
&\frac{\partial \phi}{\partial t}+\left(\textbf{\textit{u}}\cdot \nabla\right)\phi=\nabla\cdot(\mathcal{M}\nabla\mu_\phi), \label{the2.18a}\\
&\mu_\phi=-\epsilon^2\nabla \cdot \left(\lambda\nabla \phi\right)+\lambda F'(\phi)+2M\mathcal{G}c_2\phi-N\mathcal{G}c_2(\phi^3-\phi)+2M\mathcal{G}c_3\phi-N\mathcal{G}c_3(\phi^3-\phi), \label{the2.18b}
\end{align}
\end{subequations}

\begin{subequations}\label{the2.19}
	\begin{align}
&\frac{\partial c_1}{\partial t}+\left(\textbf{\textit{u}}\cdot \nabla\right)c_1=\nabla\cdot\left(D_1c_1\mathcal{G}^{-1}\nabla \mu_1\right)-a_1\mathcal{G}^{-1}P(\phi)(a_1\mu_1+a_2\mu_2-a_3\mu_3), \label{the2.19a}\\
&\mu_1=\mathcal{G}\text{ln}\left(\frac{c_1}{c_0}\right), \label{the2.19b}\\
&\frac{1}{D_1}=\frac{(\phi^2-1)^2}{q\epsilon}+\frac{1-\phi}{2D_1^{-}}+\frac{1+\phi}{2D_1^{+}},\label{the2.19c}
\end{align}
\end{subequations}

\begin{subequations}\label{the2.20}
\begin{align}
&\frac{\partial c_2}{\partial t}+\left(\textbf{\textit{u}}\cdot \nabla\right)c_2=D_2\mathcal{G}^{-1}\nabla\cdot\left(c_2\nabla \mu_2\right)-a_2\mathcal{G}^{-1}P(\phi)(a_1\mu_1+a_2\mu_2-a_3\mu_3), \label{the2.20a}\\
&\mu_2=\mathcal{G}\text{ln}\left(\frac{c_2}{c_0}\right)+M\mathcal{G}\phi^2-\frac{N}{4}\mathcal{G}(\phi^2-1)^2,\label{the2.20b}
\end{align}
\end{subequations}

\begin{subequations}\label{the2.21}
	\begin{align}
		&\frac{\partial c_3}{\partial t}+\left(\textbf{\textit{u}}\cdot \nabla\right)c_3=D_3\mathcal{G}^{-1}\nabla\cdot\left(c_3\nabla \mu_3\right)+a_3\mathcal{G}^{-1}P(\phi)(a_1\mu_1+a_2\mu_2-a_3\mu_3), \label{the2.21a}\\
		&\mu_3=\frac{\partial \lambda}{\partial c_3}\left(F(\phi)+\frac{\epsilon^2}{2}\left|\nabla \phi\right|^2\right)+\mathcal{G}\text{ln}\left(\frac{c_3}{c_0}\right)+M\mathcal{G}\phi^2-\frac{N}{4}\mathcal{G}(\phi^2-1)^2.
	\end{align}
\end{subequations}

By introducing the following dimensionless quantities:
\begin{equation}\label{the2.22}
	\begin{split}
&\textbf{\textit{x}}^*=\frac{\textbf{\textit{x}}}{L}, \quad  \textbf{\textit{u}}^*=\frac{\textbf{\textit{u}}}{U}, \quad t^*=\frac{t}{T}, \quad T=\frac{L}{U}, \quad D_1^{*}=\frac{D_1}{D_0}, \quad D_2^{*}=\frac{D_2}{D_0},\\
& D_3^{*}=\frac{D_3}{D_0}, \quad K^*=\frac{KL}{D_0}, \quad c_1^*=\frac{c_1}{c_0}, \quad c_2^*=\frac{c_2}{c_0}, \quad c_3^*=\frac{c_3}{c_0}, \quad \epsilon^*=\frac{\epsilon}{L},
\end{split}
\end{equation}
where $L, U, T, D_0$ are the characteristic length, velocity, time and diffusion coefficient. In addition, we can choose $c_0=\frac{\nu U}{\mathcal{G}L}$ to define the dimensionless variable $\lambda^*=\frac{\lambda}{\mathcal{G}c_0}=\frac{\lambda L}{\nu U}$.

Let $c_i+1=c_i^* \ (i=1,2,3)$ avoiding the singularity of $\text{ln}(c_i^*)$ caused by $c_i^*=0$. The dimensionless equations are as follows ($*$ is omitted for the convenience sake):
\begin{subequations}\label{the2.23}
\begin{align}
&Re\left(\frac{\partial \textbf{\textit{u}}}{\partial t}+\left(\textbf{\textit{u}}\cdot \nabla\right)\textbf{\textit{u}}\right)=\Delta \textbf{\textit{u}}-\nabla p-\phi\nabla\mu_\phi-c_1\nabla\mu_1-c_2\nabla\mu_2-c_3\nabla\mu_3, \label{the2.23a}\\
&\nabla\cdot \textbf{\textit{u}}=0, \label{the2.23b} \\
&\frac{\partial \phi}{\partial t}+\nabla\cdot\left(\textbf{\textit{u}}\phi\right)=\nabla\cdot(\mathcal{M}\nabla\mu_\phi), \label{the2.23c}\\
&\mu_\phi=-\epsilon^2\nabla \cdot \left(\lambda\nabla \phi\right)+\lambda F'(\phi)+2M(c_2+1)\phi-N(c_2+1)(\phi^3-\phi) \label{the2.23d}\\
&\quad +2M(c_3+1)\phi-N(c_3+1)(\phi^3-\phi), \notag\\
&\frac{\partial c_1}{\partial t}+\nabla\cdot\left(\textbf{\textit{u}}c_1\right)=\frac{1}{Pe}\nabla\cdot\left(D_1(c_1+1)\nabla \mu_1\right)-a_1P(\phi)(a_1\mu_1+a_2\mu_2-a_3\mu_3), \label{the2.23e}\\
&\mu_1=\text{ln}\left(c_1+1\right), \label{the2.23f}\\
&\frac{1}{D_1}=\frac{(\phi^2-1)^2}{q\epsilon}+\frac{1-\phi}{2D_1^{-}}+\frac{1+\phi}{2D_1^{+}},\label{the2.23g}\\
&\frac{\partial c_2}{\partial t}+\nabla\cdot\left(\textbf{\textit{u}}c_2\right)=\frac{1}{Pe}D_2\nabla\cdot\left((c_2+1)\nabla \mu_2\right)-a_2P(\phi)(a_1\mu_1+a_2\mu_2-a_3\mu_3), \label{the2.23h}\\
&\mu_2=\text{ln}\left(c_2+1\right)+M\phi^2-\frac{N}{4}(\phi^2-1)^2, \label{the2.23i}\\
&\frac{\partial c_3}{\partial t}+\nabla\cdot\left(\textbf{\textit{u}}c_3\right)=\frac{1}{Pe}D_3\nabla\cdot\left((c_3+1)\nabla \mu_3\right)+a_3P(\phi)(a_1\mu_1+a_2\mu_2-a_3\mu_3), \label{the2.23j}\\
&\mu_3=\frac{\partial \lambda}{\partial c_3}\left(F(\phi)+\frac{\epsilon^2}{2}\left|\nabla \phi\right|^2\right)+\text{ln}\left(c_3+1\right)+M\phi^2-\frac{N}{4}(\phi^2-1)^2, \label{the2.23k}
\end{align}
\end{subequations}
with the dimensionless parameters $Re=\frac{\rho UL}{\nu}$ and $Pe=\frac{UL}{D_0}$. We note the no-slip boundary for $\textbf{\textit{u}}$ and the homogeneous Neumann boundary for $\phi$, $\mu_\phi$, $c_1$, $\mu_1$, $c_2$, $\mu_2$, $c_3$, $\mu_3$. Here, we redefine the pressure
\begin{equation}\nonumber
\begin{split}
p&=p+\lambda\left(\frac{\epsilon^2}{2}|\nabla\phi|^2+F(\phi)\right)-\phi\mu_\phi
-c_1\mu_1-c_2\mu_2-c_3\mu_3-\phi(M(c_2+1)+M(c_3+1))\\
&+\frac{(\phi^2-1)^2}{4}(N(c_2+1)+N(c_3+1))+(c_1+1)\left(\text{ln}(c_1+1)-1\right)\\
&+(c_2+1)\left(\text{ln}(c_2+1)-1\right)+(c_3+1)\left(\text{ln}(c_3+1)-1\right).
\end{split}
\end{equation}

\subsection{Mass conservation and energy law}
We introduce the following theorem.
\begin{theorem}
\label{theorem2.1}
The total mass of system (\ref{the2.23}) is conserved
\begin{equation}\label{the2.24}
\begin{split}
\frac{d}{dt}\int_\Omega (c_1+c_2+c_3)d\textbf{x}=0,
\end{split}
\end{equation}
and satisfies the following energy dissipation law
\begin{equation}\label{the2.25}
	\begin{split}
\frac{d}{dt}\mathcal{E}_{\text{total}}&=-\left(\Vert\nabla \textbf{u}\Vert_{L^2}^2+\mathcal{M}\Vert\nabla\mu_\phi\Vert_{L^2}^2+\int_{\Omega}\frac{D_1(c_1+1)}{Pe}|\nabla\mu_1|^2 d\textbf{x}+\int_{\Omega}\frac{D_2(c_2+1)}{Pe}|\nabla\mu_2|^2 d\textbf{x} \right.\\
&\left. +\int_{\Omega}\frac{D_3(c_3+1)}{Pe}|\nabla\mu_3|^2 d\textbf{x}+\Vert\sqrt{P(\phi)}(a_1\mu_1+a_2\mu_2-a_3\mu_3)\Vert_{L^2}^2\right)\leq 0,
	\end{split}
\end{equation}
where
\begin{equation}\label{the2.26}
	\begin{split}
\mathcal{E}_{\text{total}}&=\int_\Omega \left(\frac{Re}{2}\textbf{u}^2+\frac{\lambda\epsilon^2}{2}|\nabla\phi|^2+\lambda F(\phi)+(c_1+1)\left(\text{ln}(c_1+1)-1\right)+(c_2+1)\left(\text{ln}(c_2+1)-1\right) \right.\\
&\left. +(c_3+1)\left(\text{ln}(c_3+1)-1\right) +\phi^2\left(M(c_2+1)+M(c_3+1)\right)-\frac{(\phi^2-1)^2}{4}\left(N(c_2+1)+N(c_3+1)\right)\right)d\textbf{x}.
	\end{split}
\end{equation}

\begin{proof}
Note that $a_1+a_2-a_3=0$, integrating (\ref{the2.23e}), (\ref{the2.23h}) and (\ref{the2.23j}) over $\Omega$, and summing that to obtain (\ref{the2.24}).
	
Taking the $L^2$ inner product of (\ref{the2.23a}), (\ref{the2.23c}) and (\ref{the2.23d}) with $\textbf{\textit{u}}$, $\mu_\phi$ and $\frac{\partial \phi}{\partial t}$, we obtain
\begin{equation}\label{the2.27}
\begin{split}
&\frac{Re}{2}\frac{d}{dt}\Vert \textbf{\textit{u}}\Vert_{L^2}^2+\int_\Omega \lambda\left(\frac{\epsilon^2}{2}\frac{\partial}{\partial t}|\nabla\phi|^2+\frac{\partial}{\partial t}F(\phi)\right)d\textbf{x}+\int_\Omega\left(M(c_2+1)+M(c_3+1)\right)\frac{\partial}{\partial t}\phi^2d\textbf{x}\\
&-\int_\Omega\left(N(c_2+1)+N(c_3+1)\right)\frac{\partial}{\partial t	}\frac{\left(\phi^2-1\right)^2}{4}d\textbf{x}\\
&=-\Vert\nabla \textbf{\textit{u}}\Vert_{L^2}^2-\mathcal{M}\Vert\nabla\mu_\phi\Vert_{L^2}^2-\int_{\Omega}(\phi\nabla\mu_\phi+c_1\nabla\mu_1+c_2\nabla\mu_2+c_3\nabla\mu_3)\cdot \textbf{\textit{u}}d\textbf{x}.
\end{split}
\end{equation}

Taking the $L^2$ inner product of (\ref{the2.23e}), (\ref{the2.23f}), (\ref{the2.23h}), (\ref{the2.23i}), (\ref{the2.23j}) and (\ref{the2.23k}) with $\mu_1$, $\frac{\partial c_1}{\partial t}$, $\mu_2$, $\frac{\partial c_2}{\partial t}$, $\mu_3$ and $\frac{\partial c_3}{\partial t}$, respectively. We then obtain
\begin{equation}\label{the2.28}
\begin{split}
&\int_\Omega \left[(c_1+1)\left(\text{ln}(c_1+1)-1\right)+(c_2+1)\left(\text{ln}(c_2+1)-1\right)+(c_3+1)\left(\text{ln}(c_3+1)-1\right)\right]d\textbf{x}\\
&+\int_\Omega \left(F(\phi)+\frac{\epsilon^2}{2}|\nabla\phi|^2\right)\frac{\partial \lambda}{\partial t}d\textbf{x}+\int_\Omega\phi^2\frac{\partial }{\partial t}\left(M(c_2+1)+M(c_3+1)\right)d\textbf{x}\\
&-\int_\Omega\frac{\left(\phi^2-1\right)^2}{4}\frac{\partial}{\partial t}\left(N(c_2+1)+N(c_3+1)\right)d\textbf{x}-\int_{\Omega}(\phi\nabla\mu_\phi+c_1\nabla\mu_1+c_2\nabla\mu_2+c_3\nabla\mu_3)\cdot \textbf{\textit{u}}d\textbf{x}\\
&=-\left(\int_{\Omega}\frac{D_1(c_1+1)}{Pe}|\nabla\mu_1|^2 d\textbf{x}+\int_{\Omega}\frac{D_2(c_2+1)}{Pe}|\nabla\mu_2|^2 d\textbf{x}+\int_{\Omega}\frac{D_3(c_3+1)}{Pe}|\nabla\mu_3|^2 d\textbf{x} \right.\\
&\left. +\Vert\sqrt{P(\phi)}(a_1\mu_1+a_2\mu_2-a_3\mu_3)\Vert_{L^2}^2\right).
\end{split}
\end{equation}

Summing up (\ref{the2.27}) and (\ref{the2.28}), and using the divergence-free condition, we obtain the desired result (\ref{the2.25}).
\end{proof}
\end{theorem}

\section{Structure-preserving numerical scheme}
\label{section3}
In this section, we first construct the time-discrete scheme with the first-order accuracy for the proposed phase field model, and then prove the mass conservation and energy stability for the numerical scheme. Finally, the fully discre finite element scheme is given.

Let $N$ be a positive integer and $T$ be the final time of computation. We set $\Delta t=T/N$ be the uniform time step. Let $(\cdot)^n$ be the numerical approximation of a specific variable at $t=n\Delta t$.

\subsection{Time-discrete scheme}
\label{Section3.1}
Following the stabilized method proposed in \cite{shen2010numerical}, we assume that the potential function $F(\phi)$ satisfies the condition: there exists a constant $L$ such that $\max\limits_{|\phi|\in \mathbb{R}}\left|F''(\phi)\right|\leq L$. We note that the double-well potential $F(\phi)=\frac{1}{4}\left(\phi^2-1\right)^2$ satisfies this condition by truncating it to quadratic growth outside of an interval $[-M, M]$ without affecting the solution if the maximum norm of the initial condition $\phi_0$ is bounded by $M$. It has been a common practice to deal with phase-field problems \cite{qin2022phase, shen2015decoupled, meng2025convergence}.

The first-order time-discrete scheme is constructed as follows:
\begin{subequations}\label{the3.1}
\begin{align}
&Re\left(\frac{\textbf{\textit{u}}^{n+1}-\textbf{\textit{u}}^{n}}{\Delta t}+\left(\textbf{\textit{u}}^{n}\cdot \nabla\right)\textbf{\textit{u}}^{n+1}\right)=\Delta\textbf{\textit{u}}^{n+1}-\nabla p^{n+1} \label{the3.1a}\\
&-\phi^{n}\nabla\mu_\phi^{n+1}-c_1^{n+1}\nabla\mu_1^{n+1}-c_2^{n+1}\nabla\mu_2^{n+1}-c_3^{n+1}\nabla\mu_3^{n+1}, \nonumber\\
&\nabla\cdot \textbf{\textit{u}}^{n+1}=0, \label{the3.1b}
	\end{align}
\end{subequations}

\begin{subequations}\label{the3.2}
\begin{align}
&\frac{\phi^{n+1}-\phi^{n}}{\Delta t}+\nabla\cdot(\textbf{\textit{u}}^{n+1} \phi^{n})=\nabla\cdot\left(\mathcal{M}\nabla \mu_\phi^{n+1}\right), \label{the3.2a}\\
&\mu_\phi^{n+1}=-\epsilon^2\nabla\cdot\left(\lambda^n\nabla\phi^{n+1}\right)+\lambda^nF'(\phi^n)+\lambda^nS(\phi^{n+1}-\phi^{n})+2\phi^{n+1}(M(c_2^{n+1}+1)+M(c_3^{n+1}+1)) \label{the3.2b}\\
&-\left((\phi^n)^3-\phi^n\right)\left(N(c_2^n+1)+N(c_3^n+1)\right)+S(\phi^{n+1}-\phi^{n})\left(N(c_2^n+1)+N(c_3^n+1)\right), \notag
\end{align}
\end{subequations}

\begin{subequations}\label{the3.3}
\begin{align}
&\frac{c_1^{n+1}-c_1^{n}}{\Delta t}+\nabla\cdot\left(\textbf{\textit{u}}^{n+1}c_1^{n+1}\right)=\frac{1}{Pe}\nabla\cdot\left(D_1^n(c_1^{n+1}+1)\nabla \mu_1^{n+1}\right)-a_1P(\phi^n)(a_1\mu_1^{n+1}+a_2\mu_2^{n+1}-a_3\mu_3^{n+1}), \label{the3.3a}\\
&\mu_1^{n+1}=\text{ln}\left(c_1^{n+1}+1\right), \label{the3.3b}\\
&\frac{1}{D_1^n}=\frac{((\phi^n)^2-1)^2}{q^n\epsilon}+\frac{1-\phi^n}{2D_1^{-}}+\frac{1+\phi^n}{2D_1^{+}},\label{the3.3c}
\end{align}
\end{subequations}

\begin{subequations}\label{the3.4}
\begin{align}
&\frac{c_2^{n+1}-c_2^{n}}{\Delta t}+\nabla\cdot\left(\textbf{\textit{u}}^{n+1}c_2^{n+1}\right)=\frac{D_2}{Pe}\nabla\cdot\left((c_2^{n+1}+1)\nabla \mu_2^{n+1}\right)-a_2P(\phi^n)(a_1\mu_1^{n+1}+a_2\mu_2^{n+1}-a_3\mu_3^{n+1}), \label{the3.4a}\\
&\mu_2^{n+1}= \text{ln}\left(c_2^{n+1}+1\right)+M(\phi^n)^2-\frac{N}{4}((\phi^{n+1})^2-1)^2, \label{the3.4b}
\end{align}
\end{subequations}

\begin{subequations}\label{the3.5}
\begin{align}
&\frac{c_3^{n+1}-c_3^{n}}{\Delta t}+\nabla\cdot\left(\textbf{\textit{u}}^{n+1}c_3^{n+1}\right)=\frac{D_3}{Pe}\nabla\cdot\left((c_3^{n+1}+1)\nabla \mu_3^{n+1}\right)+a_3P(\phi^n)(a_1\mu_1^{n+1}+a_2\mu_2^{n+1}-a_3\mu_3^{n+1}), \label{the3.5a}\\
&\mu_3^{n+1}=\frac{\lambda^{n+1}-\lambda^{n+1}}{c_3^{n+1}-c_3^{n}}\left(F(\phi^{n+1})+\frac{\epsilon^2}{2}\left|\nabla \phi^{n+1}\right|^2\right)+\text{ln}\left(c_3^{n+1}+1\right)+M(\phi^n)^2-\frac{N}{4}((\phi^{n+1})^2-1)^2. \label{the3.5b}
\end{align}
\end{subequations}
Here, $S\geq 0$ in (\ref{the3.2b}) is the stabilization parameter. The following homogeneous boundary conditions on the boundary:
\begin{equation}\nonumber
\begin{split}
&\nabla c_1^{n+1}\cdot \textbf{\textit{n}}|_{\partial \Omega}=\nabla \mu_1^{n+1}\cdot \textbf{\textit{n}}|_{\partial \Omega}=\nabla c_2^{n+1}\cdot \textbf{\textit{n}}|_{\partial \Omega}=\nabla \mu_2^{n+1}\cdot \textbf{\textit{n}}|_{\partial \Omega}=\nabla c_3^{n+1}\cdot \textbf{\textit{n}}|_{\partial \Omega}=\nabla \mu_3^{n+1}\cdot \textbf{\textit{n}}|_{\partial \Omega}=0,\\
&\textbf{\textit{u}}^{n+1}|_{\partial \Omega}=0, \ \nabla \phi^{n+1}\cdot \textbf{\textit{n}}|_{\partial \Omega}=\nabla \mu_\phi^{n+1}\cdot \textbf{\textit{n}}|_{\partial \Omega}=0.
\end{split}
\end{equation}

We next prove the mass conservation and energy stability for the proposed numerical scheme. The discrete energy is defined as:
\begin{equation}\label{the3.6}
\begin{split}
\mathcal{E}^{n}&=\frac{Re}{2}\Vert \textbf{\textit{u}}^n\Vert_{L^2}^2+\int_\Omega \left(\frac{\epsilon^2}{2}\lambda^{n}|\nabla\phi^n|^2+\lambda^nF(\phi^n)\right)d\mathbf{x}+\int_\Omega (c_1^n+1)\left(\text{ln}(c_1^n+1)-1\right)d\mathbf{x}\\
&+\int_\Omega (c_2^n+1)\left(\text{ln}(c_2^n+1)-1\right)d\mathbf{x}+\int_\Omega (c_3^n+1)\left(\text{ln}(c_3^n+1)-1\right)d\mathbf{x}+\int_\Omega(\phi^n)^2\left(M(c_2^n+1)+M(c_3^n+1)\right)d\textbf{x}\\
&-\int_\Omega\frac{((\phi^n)^2-1)^2}{4}\left(N(c_2^n+1)+N(c_3^n+1)\right)d\textbf{x}.
\end{split}
\end{equation}

\begin{theorem}
	\label{theorem3.1}
The time-discrete scheme (\ref{the3.1})-(\ref{the3.5}) satisfies the mass conservation in the sense that
\begin{equation}\label{the3.7}
\begin{split}
\int_\Omega	(c_1^n+c_2^n+c_3^n)d\textbf{x}=\int_\Omega	(c_1^0+c_2^0+c_3^0)d\textbf{x}.
\end{split}
\end{equation}	
Let $S\geq \frac{L}{2}$, it is energy stable in the sense that the following discrete energy law holds:
\begin{equation}\label{the3.8}
\begin{split}
&\mathcal{E}^{n+1}-\mathcal{E}^{n}+\mathcal{T}^{n+1,n}=  -\Delta t\left(\mathcal{M}\Vert\nabla\mu_\phi^{n+1}\Vert_{L^2}^2+\int_\Omega   \frac{D_1^n(c_1^{n+1}+1)}{Pe}\left|\nabla\mu_1^{n+1}\right|^2+   \frac{D_2(c_2^{n+1}+1)}{Pe}\left|\nabla\mu_2^{n+1}\right|^2 \right.\\
&\left. +   \frac{D_3(c_3^{n+1}+1)}{Pe}\left|\nabla\mu_3^{n+1}\right|^2d\textbf{x}+\Vert\sqrt{P(\phi^n)}(a_1\mu_1^{n+1}+a_2\mu_2^{n+1}-a_3\mu_3^{n+1})\Vert_{L^2}^2\right)\leq 0.
\end{split}
\end{equation}	
Here, $\mathcal{T}^{n+1,n}$ represents non-negative terms that are independent of the energy.
\begin{proof}
Integrating (\ref{the3.3a}), (\ref{the3.4a}) and (\ref{the3.5a}), and using integration by parts, we can obtain (\ref{the3.7}).
	
Taking the inner product of (\ref{the3.1a}), (\ref{the3.2a}) and (\ref{the3.2b}) with $\Delta t\textbf{\textit{u}}^{n+1}$, $\Delta t\mu_\phi^{n+1}$ and $\phi^{n+1}-\phi^{n}$, respectively, we have
\begin{equation}\label{the3.9}
\begin{split}
&\frac{Re}{2}\left(\Vert\textbf{\textit{u}}^{n+1}\Vert_{L^2}^2-\Vert\textbf{\textit{u}}^{n}\Vert_{L^2}^2+\Vert\textbf{\textit{u}}^{n+1}-\textbf{\textit{u}}^{n}\Vert_{L^2}^2\right)+\Delta t\Vert\nabla\textbf{\textit{u}}^{n+1}\Vert_{L^2}^2\\
&=-\Delta t\left(\phi^{n}\nabla\mu_\phi^{n+1}+c_1^{n+1}\nabla\mu_1^{n+1}+c_2^{n+1}\nabla\mu_2^{n+1}+c_3^{n+1}\nabla\mu_3^{n+1},\textbf{\textit{u}}^{n+1}\right),
\end{split}
\end{equation}	
\begin{equation}\label{the3.10}
\begin{split}
&\frac{\epsilon^2}{2}\int_\Omega\lambda^n\left(\left|\nabla\phi^{n+1}\right|^2+\left|\nabla\phi^{n}\right|^2+\left|\nabla(\phi^{n+1}-\phi^{n})\right|^2\right)d\textbf{x}+\int_\Omega\lambda^n \left(F(\phi^{n+1})-F(\phi^{n})\right)d\textbf{x}\\
&-\int_\Omega\frac{\lambda^nF''(\xi_\phi^n)}{2}(\phi^{n+1}-\phi^{n})^2d\textbf{x}+\int_\Omega\lambda^n S(\phi^{n+1}-\phi^{n})^2d\textbf{x}\\
&+\int_\Omega \left((\phi^{n+1})^2-(\phi^{n})^2+\left(\phi^{n+1}-\phi^n\right)^2\right)\left(M(c_2^{n+1}+1)+M(c_3^{n+1}+1)\right) d\textbf{x}\\
&-\int_\Omega\left(\frac{((\phi^{n+1})^2-1)^2}{4}-\frac{((\phi^n)^2-1)^2}{4}\right)\left(N(c_2^n+1)+N(c_3^n+1)\right)d\textbf{x}\\
&+\int_\Omega\frac{F''(\xi_\phi^n)}{2}(\phi^{n+1}-\phi^{n})^2\left(N(c_2^n+1)+N(c_3^n+1)\right)d\textbf{x}+\int_\Omega S(\phi^{n+1}-\phi^{n})^2\left(N(c_2^n+1)+N(c_3^n+1)\right)d\textbf{x}\\
&=-\mathcal{M}\Delta t\Vert\nabla\mu_\phi^{n+1}\Vert_{L^2}^2+\Delta t(\phi^n\nabla \mu_\phi^{n+1},\textbf{\textit{u}}^{n+1}),
\end{split}
\end{equation}	
where we use the identity $(a,a-b)=\frac{1}{2}(a^2-b^2)+\frac{1}{2}(a-b)^2$ and the expansion:
\begin{equation}\label{the3.11}
	\begin{split}
F(\phi^{n+1})-F(\phi^{n})=F'(\phi^{n})\left(\phi^{n+1}-\phi^{n}\right)+\frac{F''(\xi_\phi^n)}{2}\left(\phi^{n+1}-\phi^{n}\right)^2.
	\end{split}
\end{equation}	
	
Taking the inner product of (\ref{the3.3a}) and (\ref{the3.3b}) with $\Delta t\mu_1^{n+1}$ and $c_1^{n+1}-c_1^{n}$, respectively. We thus obtain
\begin{equation}\label{the3.12}
\begin{split}
&\int_\Omega\left(c_1^{n+1}+1\right)\left(\text{ln}\left(c_1^{n+1}+1\right)-1\right)d\textbf{x}-\int_\Omega\left(c_1^{n}+1\right)\left(\text{ln}\left(c_1^{n}+1\right)-1\right)d\textbf{x}\\
&=-\int_\Omega   \frac{\Delta tD_1^n(c_1^{n+1}+1)}{Pe}\left|\nabla\mu_1^{n+1}\right|^2d\textbf{x}+\Delta t\left(c_1^{n+1}\nabla\mu_1^{n+1},\textbf{\textit{u}}^{n+1}\right)-\int_\Omega\frac{1}{2(\xi_{c_1}^n+1)}\left(c_1^{n+1}-c_1^n\right)^2 d\textbf{x}\\
&-\Delta t\left(a_1P(\phi^n)(a_1\mu_1^{n+1}+a_2\mu_2^{n+1}-a_3\mu_3^{n+1}), \mu_1^{n+1}\right),
\end{split}
\end{equation}	
where we use the following fact:
\begin{equation}\label{the3.13}
\begin{split}
&\left(c_1^{n+1}+1\right)\left(\text{ln}\left(c_1^{n+1}+1\right)-1\right)-\left(c_1^{n}+1\right)\left(\text{ln}\left(c_1^{n}+1\right)-1\right)\\
&=\text{ln}\left(c_1^{n+1}+1\right)(c_1^{n+1}-c_1^{n})-\frac{1}{2(\xi_{c_1}^n+1)}(c_1^{n+1}-c_1^{n})^2.
\end{split}
\end{equation}	
	
Similarly, we take the inner product of (\ref{the3.4a}), (\ref{the3.4b}), (\ref{the3.5a}), and (\ref{the3.5b}) with $\Delta t\mu_2^{n+1}$, $c_2^{n+1}-c_2^{n}$, $\Delta t\mu_3^{n+1}$ and $c_3^{n+1}-c_3^{n}$, respectively, we obtain
\begin{equation}\label{the3.14}
\begin{split}
&\int_\Omega\left(c_2^{n+1}+1\right)\left(\text{ln}\left(c_2^{n+1}+1\right)-1\right)d\textbf{x}-\int_\Omega\left(c_2^{n}+1\right)\left(\text{ln}\left(c_2^{n}+1\right)-1\right)d\textbf{x}\\
&+M\int_{\Omega}(\phi^n)^2\left((c_2^{n+1}+1)-(c_2^{n}+1)\right)d\textbf{x}-\frac{N}{4}\int_{\Omega}\left((\phi^{n+1})^2-1\right)^2\left((c_2^{n+1}+1)-(c_2^{n}+1)\right)d\textbf{x}\\
&=-\int_\Omega   \frac{\Delta tD_2(c_2^{n+1}+1)}{Pe}\left|\nabla\mu_2^{n+1}\right|^2d\textbf{x}+\Delta t\left(c_2^{n+1}\nabla\mu_2^{n+1},\textbf{\textit{u}}^{n+1}\right)-\int_\Omega\frac{1}{\xi_{c_2}^n+1}\left(c_2^{n+1}-c_2^n\right)^2 d\textbf{x}\\
&-\Delta t\left(a_2P(\phi^n)(a_1\mu_1^{n+1}+a_2\mu_2^{n+1}-a_3\mu_3^{n+1}), \mu_2^{n+1}\right),
\end{split}
\end{equation}	
and
\begin{equation}\label{the3.15}
\begin{split}
&\int_\Omega\left(c_3^{n+1}+1\right)\left(\text{ln}\left(c_3^{n+1}+1\right)-1\right)d\textbf{x}-\int_\Omega\left(c_3^{n}+1\right)\left(\text{ln}\left(c_3^{n}+1\right)-1\right)d\textbf{x}\\
&+\int_\Omega \left(\lambda^{n+1}-\lambda^{n}\right)\left(F(\phi^{n+1})+\frac{\epsilon^2}{2}\left|\nabla \phi^{n+1}\right|^2\right)d\textbf{x}\\
&+M\int_{\Omega}(\phi^n)^2\left((c_3^{n+1}+1)-(c_3^{n}+1)\right)d\textbf{x}-\frac{N}{4}\int_{\Omega}\left((\phi^{n+1})^2-1\right)^2\left((c_3^{n+1}+1)-(c_3^{n}+1)\right)d\textbf{x}\\
&=-\int_\Omega   \frac{\Delta tD_3(c_3^{n+1}+1)}{Pe}\left|\nabla\mu_3^{n+1}\right|^2d\textbf{x}+\Delta t\left(c_3^{n+1}\nabla\mu_3^{n+1},\textbf{\textit{u}}^{n+1}\right)-\int_\Omega\frac{1}{\xi_{c_3}^n+1}\left(c_3^{n+1}-c_3^n\right)^2 d\textbf{x}\\
&+\Delta t\left(a_3P(\phi^n)(a_1\mu_1^{n+1}+a_2\mu_2^{n+1}-a_3\mu_3^{n+1}), \mu_3^{n+1}\right).
\end{split}
\end{equation}	

Combining (\ref{the3.9}), (\ref{the3.10}), (\ref{the3.12}), (\ref{the3.14}) and (\ref{the3.15}), one obtains immediately (\ref{the3.8}).

The proof is completed.
\end{proof}
\end{theorem}

\begin{remark}
We can see that the governing system is highly nonlinear, and the constructed numerical scheme preserves the decay of the original energy. To solve the system, we use the Newton iteration method, which may increase the computational cost. The recently developed SAV approach \cite{shen2019new, shen2018scalar} may be used to design higher-order decoupled schemes for this problem, which will be considered in future work.
\end{remark}

\subsection{Fully discrete finite element scheme}
The fully discrete finite element scheme is developed in this subsection. Let $\mathcal{T}_h$ be a regular triangulation of $\Omega\subset\mathbb{R}^2$ with mesh size $h$. We use $\mathcal{P}_l$ to denote the space of continuous piecewise polynomials of total degree at most $l$. Several continuous finite element spaces are introduced as follows:
\begin{equation}\label{the3.16}
\begin{split}
&\mathbf{V}_h=\left\{\mathbf{v}_h\in C(\Omega): \mathbf{v}_h|_K\in (\mathcal{P}_2)^2, \ \forall K\in\mathcal{T}_h,  \mathbf{v}_h|_{\partial\Omega}=0\right\},\\
& Q_h=\left\{q_h\in C(\Omega): q_h|_K\in \mathcal{P}_1, \ \forall K\in\mathcal{T}_h, q_h|_{\partial\Omega}=0 \right\},\\
&X_h=\left\{(w_1)_h\in C(\Omega): (w_1)_h|_K\in \mathcal{P}_2, \ \forall K\in\mathcal{T}_h \right\}, \ Y_h=\left\{(w_2)_h\in C(\Omega): (w_2)_h|_K\in \mathcal{P}_2, \ \forall K\in\mathcal{T}_h \right\}.
\end{split}
\end{equation}	
We assume the pair of spaces $\left(\mathbf{V}_h, Q_h\right)$ satisfy the \textit{Inf-Sup} condition \cite{girault2012finite}. Let
\begin{equation}\nonumber
	\begin{split}
\mathbf{W}_h=\mathbf{V}_h\times Q_h\times X_h\times X_h\times X_h\times Y_h\times X_h\times Y_h\times X_h\times Y_h,
	\end{split}
\end{equation}	

The fully-discrete numerical scheme of system (\ref{the2.23}) reads as follows: find
\begin{equation}\nonumber
	\begin{split}
\left(\textbf{\textit{u}}^{n+1}_h, p^{n+1}_h, \phi_h^{n+1}, (\mu_\phi)_h^{n+1}, (c_1)_h^{n+1}, (\mu_1)_h^{n+1}, (c_2)_h^{n+1}, (\mu_2)_h^{n+1}, (c_3)_h^{n+1}, (\mu_3)_h^{n+1}\right)\in \mathbf{W}_h,
	\end{split}
\end{equation}	
such that for all $\left(\textbf{\textit{v}}_h, q_h, \psi_h, w_h, (\psi_1)_h, (w_1)_h, (\psi_2)_h, (w_2)_h, (\psi_3)_h, (w_3)_h\right)\in \mathbf{W}_h$, there hold
\begin{subequations}\label{the3.17}
\begin{align}
&Re\left(\left(\frac{\textbf{\textit{u}}_h^{n+1}-\textbf{\textit{u}}_h^{n}}{\Delta t},\textbf{\textit{v}}_h \right)+\left(\left(\textbf{\textit{u}}^{n}_h\cdot \nabla\right)\textbf{\textit{u}}^{n+1}_h, \textbf{\textit{v}}_h\right)\right)+\left(D_\nu(\textbf{\textit{u}}_h^{n+1}),D_\nu(\textbf{\textit{v}}_h)\right)+\left(\nabla p^{n+1}_h, \textbf{\textit{v}}_h\right) \label{the3.17a}\\
&=-\left(\phi^{n}_h\nabla(\mu_\phi)^{n+1}_h+(c_1)^{n+1}_h\nabla(\mu_1)^{n+1}_h+(c_2)_h^{n+1}\nabla(\mu_2)_h^{n+1}+(c_3)_h^{n+1}\nabla(\mu_3)_h^{n+1},\textbf{\textit{v}}_h\right), \nonumber\\
&\left(\nabla\cdot \textbf{\textit{u}}^{n+1}_h,q_h\right)=0, \label{the3.17b}\\
&\left(\frac{\phi^{n+1}_h-\phi^{n}_h}{\Delta t}, \psi_h\right)-\left(\textbf{\textit{u}}^{n+1}_h \phi^{n}_h, \nabla\psi_h\right)+\mathcal{M}\left(\nabla (\mu_\phi)^{n+1}_h,\nabla\psi_h\right)=0, \label{the3.17c}\\
&\left((\mu_\phi)_h^{n+1}, w_h \right)=\epsilon^2\left(\lambda^n_h\nabla\phi^{n+1}_h, \nabla w_h\right)+\left(\lambda^n_hF'(\phi^n_h), w_h\right)+S\left(\lambda^n_h(\phi^{n+1}_h-\phi^{n}_h), w_h\right)\\
&+2\left(\phi^{n+1}_h(M((c_2)^{n+1}_h+1)+M((c_3)^{n+1}_h+1)),w_h\right) \notag\\
&+S\left((\phi^{n+1}_h-\phi^{n}_h)\left(N((c_2)^n_h+1)+N((c_3)^n_h+1)\right),w_h\right) \notag\\
&-\left(\left((\phi^n_h)^3-\phi^n_h\right)\left(N((c_2)^n_h+1)+N((c_3)^n_h+1)\right),w_h\right), \notag\\
&\left(\frac{(c_1)^{n+1}_h-(c_1)^{n}_h}{\Delta t},(\psi_1)_h\right)-\left(\textbf{\textit{u}}^{n+1}_h(c_1)^{n+1}_h,\nabla(\psi_1)_h\right)+\frac{1}{Pe}\left(D_1^n((c_1)^{n+1}_h+1)\nabla (\mu_1)^{n+1}_h,\nabla(\psi_1)_h\right)\\
&=-\left(a_1P(\phi^n_h)(a_1(\mu_1)^{n+1}_h+a_2(\mu_2)^{n+1}_h-a_3(\mu_3)^{n+1}_h),(\psi_1)_h\right), \notag\\
&\left((\mu_1)^{n+1}_h,(w_1)_h\right)=\left(\text{ln}\left((c_1)^{n+1}_h+1\right),(w_1)_h\right), \\
&\frac{1}{D_1^n}=\frac{((\phi^n_h)^2-1)^2}{q^n_h\epsilon}+\frac{1-\phi^n_h}{2D_1^{-}}+\frac{1+\phi^n_h}{2D_1^{+}}, \notag\\
&\left(\frac{(c_2)^{n+1}_h-(c_2)^{n}_h}{\Delta t},(\psi_2)_h\right)-\left(\textbf{\textit{u}}^{n+1}_h(c_2)^{n+1}_h,\nabla(\psi_2)_h\right)+\frac{D_2}{Pe}\left(((c_2)^{n+1}_h+1)\nabla (\mu_2)^{n+1}_h,\nabla(\psi_2)_h\right) \\
&=-\left(a_2P(\phi^n_h)(a_1(\mu_1)^{n+1}_h+a_2(\mu_2)^{n+1}_h-a_3(\mu_3)^{n+1}_h),(\psi_2)_h\right), \notag\\
&\left((\mu_2)^{n+1}_h,(w_2)_h\right)=\left(\text{ln}\left((c_2)^{n+1}_h+1\right),(w_2)_h\right)+M\left((\phi^n_h)^2,(w_2)_h\right)-\frac{N}{4}\left(((\phi^n_h)^2-1)^2,(w_2)_h\right), \\
&\left(\frac{(c_3)^{n+1}_h-(c_3)^{n}_h}{\Delta t},(\psi_3)_h\right)-\left(\textbf{\textit{u}}^{n+1}_h(c_3)^{n+1}_h,\nabla(\psi_3)_h\right)+\frac{D_3}{Pe}\left(((c_3)^{n+1}_h+1)\nabla (\mu_3)^{n+1}_h,\nabla(\psi_3)_h\right) \\
&=\left(a_3P(\phi^n_h)(a_1(\mu_1)^{n+1}_h+a_2(\mu_2)^{n+1}_h-a_3(\mu_3)^{n+1}_h),(\psi_3)_h\right), \notag\\
&\left((\mu_3)^{n+1}_h,(w_3)_h\right)=\left(\frac{\lambda^{n+1}_h-\lambda^n_h}{(c_3)^{n+1}_h-(c_3)^{n}_h}\left(F(\phi^{n+1}_h)+\frac{\epsilon^2}{2}\left|\nabla \phi^{n+1}_h\right|^2\right),(w_3)_h\right)\\
&+\left(\text{ln}\left((c_3)^{n+1}_h+1\right),(w_3)_h\right)+M\left((\phi^n_h)^2,(w_3)_h\right)-\frac{N}{4}\left(((\phi^n_h)^2-1)^2,(w_3)_h\right). \notag
\end{align}
\end{subequations}

By following the process of proving Theorem \ref{theorem3.1}, it is straightforward to verify that the above fully discrete scheme preserves mass conservation and the energy dissipation law.

\section{Error estimates}
\label{section4}
In this section, we present a rigorous error estimate for the time-discrete scheme applied to a simplified problem, which is formulated as a nonlinear system of the Cahn-Hilliard-Navier-Stokes equations coupled with a reaction-diffusion equation. The first-order discrete scheme corresponding to subsection \ref{Section3.1} is as follows:
\begin{subequations}\label{the4.1}
\begin{align}
&\frac{\textbf{\textit{u}}^{n+1}-\textbf{\textit{u}}^{n}}{\Delta t}+\left(\textbf{\textit{u}}^{n}\cdot \nabla\right)\textbf{\textit{u}}^{n+1}=\frac{1}{Re}\Delta\textbf{\textit{u}}^{n+1}-\nabla p^{n+1}-\phi^{n}\nabla\mu_\phi^{n+1}-\nabla c^{n+1}, \label{the4.1a}\\
&\nabla\cdot \textbf{\textit{u}}^{n+1}=0, \label{the4.1b} \\
&\frac{\phi^{n+1}-\phi^{n}}{\Delta t}+\textbf{\textit{u}}^{n+1}\cdot \nabla\phi^{n}=\mathcal{M} \Delta\mu_\phi^{n+1}, \label{the4.1c}\\
&\mu_\phi^{n+1}=-\epsilon^2\Delta\phi^{n+1}+ F'(\phi^n)+S(\phi^{n+1}-\phi^{n}), \label{the4.1d}\\
&\frac{c^{n+1}-c^{n}}{\Delta t}+\textbf{\textit{u}}^{n+1}\cdot \nabla c^{n+1}=D\Delta c^{n+1}-P(\phi^n)\text{ln}\left(c^{n+1}+1\right). \label{the4.1e}
\end{align}
\end{subequations}
Here, $D$ is a constant. The boundary condition of $\textbf{\textit{u}}$ is no-slip, and other variables follow the homogeneous Neumann boundary condtion.
\begin{remark}
It can be seen that the above system is a simplified version of (\ref{the3.1}), and the error analysis performed is a preliminary attempt as an illustrative example. 
\end{remark}

To derive the error estimate, we first give some necessary regularity assumptions for the exact solution.
\begin{assumption}
\label{assumption4.1}
We assume that the exact solution of the system (\ref{the4.1}) satisfies the following regularity conditions:
\begin{equation}\nonumber
\begin{split}
&\textbf{\textit{u}}\in L^\infty(0,T;W^{1,\infty})\cap W^{1,2}(0,T;L^2)\cap W^{2,2}(0,T;L^2),\\
&\phi\in L^\infty(0,T;W^{3,\infty})\cap W^{1,2}(0,T;L^2)\cap W^{2,2}(0,T;H^{-1}),\\
&c\in L^{\infty}(0,T; W^{1,\infty})\cap W^{2,2}(0,T;L^{2}).
\end{split}
\end{equation}	
\end{assumption}

For notational simplicity, we set $e^n_{\textbf{\textit{u}}}=\textbf{\textit{u}}(t^n)-\textbf{\textit{u}}^n$, $e^n_{p}=p(t^n)-p^n$, $e^n_{\phi}=\phi(t^n)-\phi^n$, $e^n_{\mu}=\mu_\phi(t^n)-\mu_\phi^n$, $e^n_c=c(t^n)-c^n$. The main results are stated in the following theorem.

\begin{theorem}
We assume that $S> 2L$ and assumption \ref{assumption4.1} holds. Then for the discrete scheme (\ref{the4.1a})-(\ref{the4.1e}), there exists a positive constant $C_{T43}$ independent of $\Delta t$ such that 
\begin{equation}\label{the4.2}
\begin{split}
&\Vert e_{\textbf{\textit{u}}}^{n}\Vert_{H^1}^2+\Vert e_\phi^{n}\Vert_{H^1}^2+\Vert e_c^{n}\Vert_{L^2}^2+\Delta t\sum_{k=0}^{n} \left(\Vert e_{\textbf{\textit{u}}}^{k+1}\Vert_{H^2}^2+\Vert e_p^{k+1}\Vert_{H^1}^2+\Vert \nabla e_\mu^{k}\Vert_{L^2}^2+\Vert \nabla e_c^{k}\Vert_{L^2}^2\right)\leq C_{T43}(\Delta t)^2,
\end{split}
\end{equation}	
where $n\leq N$.
\begin{proof}
The truncation form of the system (\ref{the4.1}) is as follows:	
\begin{subequations}\label{the4.3}
\begin{align}
&\frac{\textbf{\textit{u}}(t^{n+1})-\textbf{\textit{u}}(t^{n})}{\Delta t}+\left(\textbf{\textit{u}}(t^{n})\cdot \nabla\right)\textbf{\textit{u}}(t^{n+1}) \label{the4.3a}\\
&=\frac{1}{Re}\Delta\textbf{\textit{u}}(t^{n+1})-\nabla p(t^{n+1})-\phi(t^{n})\nabla\mu_\phi(t^{n+1})-\nabla c(t^{n+1})+T_{\textbf{\textit{u}}}^{n+1}, \notag\\
&\nabla\cdot \textbf{\textit{u}}(t^{n+1})=0, \label{the4.3b} \\
&\frac{\phi(t^{n+1})-\phi(t^{n})}{\Delta t}+\textbf{\textit{u}}(t^{n+1})\cdot \nabla\phi(t^{n})=\mathcal{M} \Delta\mu_\phi(t^{n+1})+T_{\phi}^{n+1}, \label{the4.3c}\\
&\mu_\phi(t^{n+1})=-\epsilon^2\Delta\phi(t^{n+1})+ F'(\phi(t^{n}))+S(\phi(t^{n+1})-\phi(t^{n}))+T_{\mu}^{n+1}, \label{the4.3d}\\
&\frac{c(t^{n+1})-c(t^{n})}{\Delta t}+\textbf{\textit{u}}(t^{n+1})\cdot \nabla c(t^{n+1})=D\Delta c(t^{n+1})-P(\phi(t^n))\text{ln}(c(t^{n+1})+1)+T_{c}^{n+1}, \label{the4.3e}
\end{align}
\end{subequations}
where $T_{\textbf{\textit{u}}}^{n+1}$, $T_{\phi}^{n+1}$, $T_{\mu}^{n+1}$ and $T_{c}^{n+1}$ are the truncation errors:
\begin{equation}\label{the4.4}
\begin{split}
&T_{\textbf{\textit{u}}}^{n+1}=\frac{\partial {\textbf{\textit{u}}}(t^{n+1})}{\partial t}-\frac{\textbf{\textit{u}}\left(t^{n+1}\right)-\textbf{\textit{u}}\left(t^{n}\right)}{\Delta t}+\left(\textbf{\textit{u}}\left(t^{n+1}\right)-\textbf{\textit{u}}\left(t^{n}\right)\right)\cdot\nabla \textbf{\textit{u}}\left(t^{n+1}\right)\\
&-\left(\phi(t^{n+1})-\phi(t^{n+1})\right)\nabla\mu_\phi(t^{n+1})\\
&=\frac{1}{\Delta t}\int_{t^n}^{t^{n+1}}(t^n-s)\frac{\partial^2 {\textbf{\textit{u}}}}{\partial s^2}ds+\nabla \textbf{\textit{u}}\left(t^{n+1}\right)\int_{t^n}^{t^{n+1}} \frac{\partial \textbf{\textit{u}}}{\partial s}ds-\nabla\mu_\phi(t^{n+1})\int_{t^n}^{t^{n+1}} \frac{\partial \phi}{\partial s}ds,
\end{split}
\end{equation}
\begin{equation}\label{the4.5}
\begin{split}
&T_{\phi}^{n+1}=\frac{\partial \phi(t^{n+1})}{\partial t}-\frac{\phi(^{n+1})-\phi(^{n})}{\Delta t}+\left(\nabla\phi(t^{n+1})-\nabla\phi(t^{n})\right)\cdot\textbf{\textit{u}}^{n+1}\\
&=\frac{1}{\Delta t}\int_{t^n}^{t^{n+1}}(t^n-s)\frac{\partial^2 \phi}{\partial s^2}ds+\nabla\phi(t^{n+1})\int_{t^n}^{t^{n+1}} \nabla\frac{\partial \phi}{\partial s}ds,
\end{split}
\end{equation}
\begin{equation}\label{the4.6}
\begin{split}
T_{\mu}^{n+1}=F'(\phi(t^{n+1}))-F'(\phi(t^{n}))+S\int_{t^n}^{t^{n+1}} \frac{\partial \phi}{\partial s}ds,
\end{split}
\end{equation}
\begin{equation}\label{the4.7}
\begin{split}
T_{c}^{n+1}=\frac{1}{\Delta t}\int_{t^n}^{t^{n+1}}(t^n-s)\frac{\partial^2 c}{\partial s^2}ds-\left(P(\phi(t^{n+1}))-P(\phi(t^{n}))\right)\text{ln}(c(t^{n+1})+1).
\end{split}
\end{equation}

Combining (\ref{the4.3}) and (\ref{the4.4})-(\ref{the4.7}), we thus obtain the error equations:
\begin{subequations}\label{the4.8}
\begin{align}
&\frac{e_{\textbf{\textit{u}}}^{n+1}-e_{\textbf{\textit{u}}}^{n}}{\Delta t}-\frac{1}{Re}\Delta e_{\textbf{\textit{u}}}^{n+1}+e_{\textbf{\textit{u}}}^{n}\cdot\nabla\textbf{\textit{u}}(t^{n+1})+\textbf{\textit{u}}^n\cdot\nabla e_{\textbf{\textit{u}}}^{n+1} \label{the4.8a}\\
&=-\nabla e_p^{n+1}-e_{\phi}^n\nabla\mu_\phi(t^{n+1})-\phi^n\nabla e_{\mu}^{n+1}-\nabla e_c^{n+1}+T_{\textbf{\textit{u}}}^{n+1}, \notag\\
&\nabla\cdot e_\textbf{\textit{u}}^{n+1}=0, \label{the4.8b} \\
&\frac{e_\phi^{n+1}-e_\phi^{n}}{\Delta t}+{\textbf{\textit{u}}}(t^{n+1})\cdot\nabla e_\phi^n+e_{\textbf{\textit{u}}}^{n+1}\cdot\nabla\phi^n-\mathcal{M} \Delta e_\mu^{n+1}=T_{\phi}^{n+1}, \label{the4.8c}\\
&e_\mu^{n+1}=-\epsilon^2\Delta e_\phi^{n+1}+ F'(\phi(t^{n}))-F'(\phi^n)+Se_\phi^{n+1}-Se_\phi^{n}+T_{\mu}^{n+1}, \label{the4.8d}\\
&\frac{e_c^{n+1}-e_c^{n}}{\Delta t}+e_\textbf{\textit{u}}^{n+1}\cdot\nabla c(t^{n+1})+\textbf{\textit{u}}^{n+1}\cdot\nabla e_c^{n+1}-D\Delta e_c^{n+1} \notag\\
&=-\left(P(\phi(t^{n}))-P(\phi^n)\right)\text{ln}(c(t^{n+1})+1)-P(\phi^n)\left(\text{ln}(c(t^{n+1})+1)-\text{ln}(c^{n+1}+1)\right)+T_{c}^{n+1}. \label{the4.8e}
\end{align}
\end{subequations}

Taking the inner product of (\ref{the4.8a}) with $2\Delta t e_\textbf{\textit{u}}^{n+1}$ and using the regularity assumption \ref{assumption4.1}, we have
\begin{equation}\label{the4.9}
\begin{split}
&\Vert e_{\textbf{\textit{u}}}^{n+1}\Vert_{L^2}^2-\Vert e_{\textbf{\textit{u}}}^{n}\Vert_{L^2}^2+\Vert e_{\textbf{\textit{u}}}^{n+1}-e_{\textbf{\textit{u}}}^{n}\Vert_{L^2}^2+\frac{2\Delta t}{Re}\Vert \nabla e_{\textbf{\textit{u}}}^{n+1}\Vert_{L^2}^2=-2\Delta t\left(e_{\textbf{\textit{u}}}^{n}\cdot\nabla\textbf{\textit{u}}(t^{n+1}), e_\textbf{\textit{u}}^{n+1}\right)\\
&-2\Delta t\left(e_{\phi}^n\nabla\mu_\phi(t^{n+1}), e_\textbf{\textit{u}}^{n+1}\right)-2\Delta t\left(\phi^n\nabla e_{\mu}^{n+1}, e_\textbf{\textit{u}}^{n+1}\right)+2\Delta t(T_{\textbf{\textit{u}}}^{n+1},e_\textbf{\textit{u}}^{n+1})\\
&\leq -2\Delta t\left(\phi^n\nabla e_{\mu}^{n+1}, e_\textbf{\textit{u}}^{n+1}\right)+C_\varepsilon\Delta t \Vert e_\textbf{\textit{u}}^{n+1}\Vert_{L^2}^2+C\Delta t \Vert e_\textbf{\textit{u}}^{n}\Vert_{L^2}^2+C\Delta t\Vert e_\phi^{n}\Vert_{L^2}^2+C\Delta t\Vert T_{\textbf{\textit{u}}}^{n+1}\Vert_{L^2}^2,
\end{split}
\end{equation}
where $C_{\varepsilon}$ is a sufficiently small positive costant.

Taking the inner product of (\ref{the4.8c}) and (\ref{the4.8d}) with $2\Delta te_\mu^{n+1}$ and $2(e_\phi^{n+1}-e_\phi^{n})$, we then obtain
\begin{equation}\label{the4.10}
\begin{split}
&(S-2L)\left(\Vert e_\phi^{n+1}\Vert_{L^2}^2-\Vert e_\phi^{n}\Vert_{L^2}^2+\Vert e_\phi^{n+1}-e_\phi^{n}\Vert_{L^2}^2\right)+{\epsilon^2}\left(\Vert \nabla e_\phi^{n+1}\Vert_{L^2}^2-\Vert \nabla e_\phi^{n}\Vert_{L^2}^2+\Vert\nabla e_\phi^{n+1}-\nabla e_\phi^{n}\Vert_{L^2}^2\right)\\
&+2\mathcal{M}\Delta t\Vert\nabla e_\mu^{n+1}\Vert_{L^2}^2\leq 2\Delta t\left(\phi^n\nabla e_{\mu}^{n+1}, e_\textbf{\textit{u}}^{n+1}\right)+C_\varepsilon\Vert e_\phi^{n+1}-e_\phi^{n}\Vert_{L^2}^2+C_\varepsilon\Delta t\Vert \nabla e_\mu^{n+1}\Vert_{L^2}^2+C\Delta t\Vert \nabla e_\phi^{n}\Vert_{L^2}^2\\
&+C\Delta t\Vert T_\phi^{n+1}\Vert_{H^{-1}}^2+2(-Se_\phi^n+T_\mu^{n+1},e_\phi^{n+1}-e_\phi^{n}),
\end{split}
\end{equation}
where we use the fact $\int_\Omega e_\mu^{n+1}e_{\textbf{\textit{u}}}^{n+1}\cdot\nabla \phi^nd\textbf{x}=-\int_\Omega \phi^n\nabla e_\mu^n\cdot e_{\textbf{\textit{u}}}^{n+1}d\textbf{x}$. 

Let $K_1=1+\max\limits_{0\leq t\leq T}\Vert \phi(t)\Vert_{L^{\infty}}$. We assume that $\Vert \phi^n\Vert_{L^\infty}\leq K_1$ for $n=0,1,2,...,N-1$, which will be verified later in the proof.

The last term on the right hand side of (\ref{the4.10}) can be estimated by
\begin{equation}\label{the4.11}
\begin{split}
&2(-Se_\phi^n+T_\mu^{n+1},e_\phi^{n+1}-e_\phi^{n})\leq  C_\varepsilon \Delta t\Vert e_{\textbf{\textit{u}}}^{n+1}\Vert_{L^2}^2+C_\varepsilon \Delta t\Vert\nabla e_\mu^{n+1}\Vert_{L^2}^2+C_{K_1}\Delta t\Vert\nabla e_\phi^n\Vert_{L^2}^2\\
&+ C\Delta t\Vert T_\phi^{n+1}\Vert_{H^{-1}}^2+C\Delta t\Vert T_\mu^{n+1}\Vert_{H^1}^2.
\end{split}
\end{equation}

Taking the inner product of (\ref{the4.8e}) with $2\Delta te_c^{n+1}$, we obtain
\begin{equation}\label{the4.12}
	\begin{split}
&\Vert e_c^{n+1}\Vert_{L^2}^2-\Vert e_c^{n}\Vert_{L^2}^2+\Vert e_c^{n+1}-e_c^n\Vert_{L^2}^2+D\Delta t\left(\Vert \nabla e_c^{n+1}\Vert_{L^2}^2-\Vert \nabla e_c^{n}\Vert_{L^2}^2+\Vert \nabla e_c^{n+1}-\nabla e_c^n\Vert_{L^2}^2\right)\\
&\leq -\frac{C_{K_1}\Delta t}{\xi_c^{n+1}+1}\Vert e_c^{n+1}\Vert_{L^2}^2+C_1 \Delta t\Vert e_{\textbf{\textit{u}}}^{n+1}\Vert_{L^2}^2+C_\varepsilon\Delta t\Vert  e_c^{n+1}\Vert_{H^1}^2+C_{K_1}\Delta t\Vert e_\phi^{n}\Vert_{L^2}^2+C\Delta t\Vert T_c^{n+1}\Vert_{L^2}^2,
	\end{split}
\end{equation}
where $\xi_c^{n+1}\in\left(c^{n+1},c(t^{n+1})\right)$.

We restrict the time step $\Delta t$ to satisfy $\Delta t<\left(C_1\right)^{-1}$. By combining with equations (\ref{the4.9})-(\ref{the4.12}), summing over $n$ and using the H\"older inequality, we have
\begin{equation}\label{the4.13}
\begin{split}
&\Vert e_{\textbf{\textit{u}}}^{n+1}\Vert_{L^2}^2+\Vert e_\phi^{n+1}\Vert_{L^2}^2+\Vert \nabla e_\phi^{n+1}\Vert_{L^2}^2+\Vert e_c^{n+1}\Vert_{L^2}^2+\sum_{k=0}^{n}\Vert  e_{\textbf{\textit{u}}}^{k+1}-e_{\textbf{\textit{u}}}^{k}\Vert_{L^2}^2\\
&+\Delta t\sum_{k=0}^{n}\Vert \nabla e_{\textbf{\textit{u}}}^{k+1}\Vert_{L^2}^2+\Delta t\sum_{k=0}^{n}\Vert \nabla e_\mu^{k+1}\Vert_{L^2}^2+\Delta t\sum_{k=0}^{n}\Vert \nabla e_c^{k+1}\Vert_{L^2}^2\\
&\leq C\Delta t\left(\Vert e_\textbf{\textit{u}}^{n}\Vert_{L^2}^2+\Vert e_\phi^{n}\Vert_{L^2}^2+\Vert \nabla e_\phi^{n}\Vert_{L^2}^2+\Vert e_c^{n}\Vert_{L^2}^2\right)\\
&+C\Delta t\sum_{k=0}^{n}\left(\Vert T_\textbf{\textit{u}}^{k+1}\Vert_{L^2}^2+\Vert T_\phi^{k+1}\Vert_{H^{-1}}^2+\Vert T_\mu^{k+1}\Vert_{H^{1}}^2+\Vert T_c^{k+1}\Vert_{L^2}^2\right)\\
&\leq C\Delta t\left(\Vert e_\textbf{\textit{u}}^{n}\Vert_{L^2}^2+\Vert e_\phi^{n}\Vert_{L^2}^2+\Vert \nabla e_\phi^{n}\Vert_{L^2}^2+\Vert e_c^{n}\Vert_{L^2}^2\right)\\
&+C\left(\Delta t\right)^2\left(\Vert {\textbf{\textit{u}}}\Vert_{W^{2,2}(0,T;L^2)}^2+\Vert {\textbf{\textit{u}}}\Vert_{W^{1,2}(0,T;L^2)}^2+\Vert \phi\Vert_{W^{1,2}(0,T;L^2)}^2+\Vert \phi\Vert_{W^{2,2}(0,T;H^{-1})}^2+\Vert c\Vert_{W^{2,2}(0,T;L^{2})}^2\right).
\end{split}
\end{equation}

By using the regularity assumption \ref{assumption4.1} and the discrete Gronwall inequality \cite{quarteroni2008numerical}, we can obtain
\begin{equation}\label{the4.14}
\begin{split}
&\Vert e_{\textbf{\textit{u}}}^{n+1}\Vert_{L^2}^2+\Vert e_\phi^{n+1}\Vert_{L^2}^2+\Vert \nabla e_\phi^{n+1}\Vert_{L^2}^2+\Vert e_c^{n+1}\Vert_{L^2}^2+\sum_{k=0}^{n}\Vert e_{\textbf{\textit{u}}}^{k+1}-e_{\textbf{\textit{u}}}^{k}\Vert_{L^2}^2\\
&+\Delta t\sum_{k=0}^{n}\Vert \nabla e_{\textbf{\textit{u}}}^{k+1}\Vert_{L^2}^2+\Delta t\sum_{k=0}^{n}\Vert \nabla e_\mu^{k+1}\Vert_{L^2}^2+\Delta t\sum_{k=0}^{n}\Vert \nabla e_c^{k+1}\Vert_{L^2}^2\leq C(\Delta t)^2.
\end{split}
\end{equation}
	
We next verify the assumption $\Vert \phi^n\Vert_{L^\infty}\leq K_1$. The proof below is based on mathematical induction.

For $k=0$, $\Vert\phi^0\Vert_{L^2}\leq 1\leq K_1$ is obvious. We assume that $\Vert\phi^k\Vert_{L^\infty}\leq K_1$ holds true for $1\leq k\leq n$, we derive $\Vert\phi^{n+1}\Vert_{L^\infty}\leq K_1$ is also true. It can be seen from (\ref{the4.13}) that $\Vert T_\mu^{n+1}\Vert_{L^2}^2\leq C\Delta t$. Applying the $H^2$-elliptic regularity of (\ref{the4.8d}), we have
\begin{equation}\label{the4.15}
	\begin{split}
\Vert e_\phi^{n+1}\Vert_{H^2}^2\leq \Vert e_\phi^{n+1}\Vert_{H^1}^2+\Vert e_\mu^{n+1}\Vert_{L^2}^2+\Vert T_\mu^{n+1}\Vert_{L^2}^2\leq C\Delta t.
	\end{split}
\end{equation}
	
Since
\begin{equation}\label{the4.16}
	\begin{split}
\Vert\phi^{n+1}\Vert_{L^\infty}\leq \Vert e_\phi^{n+1}\Vert_{L^\infty}+\Vert \phi(t^{n+1})\Vert_{L^\infty}\leq C_2\left(\Delta t\right)^{\frac{1}{2}}+\Vert \phi(t^{n+1})\Vert_{L^\infty},
	\end{split}
\end{equation}
where the Sobolev embedding theorem $H^2 \hookrightarrow L^{\infty}$ is used. We set the time step $\Delta t$ such that $\Delta t\leq (C_2)^{-2}$, leading to $\Vert\phi^{n+1}\Vert_{L^\infty}\leq K_1$.

Next, we also use induction to estimate $e_{\textbf{\textit{u}}}^{{n}}$ in $H^2$-norm . Let $K_2=\Vert {\textbf{\textit{u}}}(0)\Vert_{L^\infty}+1+\max\limits_{0\leq t\leq T}\Vert \textbf{\textit{u}}(t)\Vert_{L^{\infty}}$. Assuming $\Vert{\textbf{\textit{u}}}^{{n}}\Vert_{L^\infty}\leq K_2$ for $n=0,1,2,...,N-1$.

Taking the inner product of (\ref{the4.8a}) with $\Delta t\left(e_{\textbf{\textit{u}}}^{n+1}-e_{\textbf{\textit{u}}}^{n}\right)$, we obtain
\begin{equation}\label{the4.17}
\begin{split}
&\Vert e_{\textbf{\textit{u}}}^{n+1}-e_{\textbf{\textit{u}}}^{n}\Vert_{L^2}^2+\Delta t\left(\Vert\nabla e_{\textbf{\textit{u}}}^{n+1}\Vert_{L^2}^2-\Vert\nabla e_{\textbf{\textit{u}}}^{n}\Vert_{L^2}^2+\Vert\nabla e_{\textbf{\textit{u}}}^{n+1}-\nabla e_{\textbf{\textit{u}}}^{n}\Vert_{L^2}^2\right)\leq C_\varepsilon\Vert e_{\textbf{\textit{u}}}^{n+1}-e_{\textbf{\textit{u}}}^{n}\Vert_{L^2}^2\\
& +C(\Delta t)^2\Vert e_{\textbf{\textit{u}}}^{n} \Vert_{L^2}^2+C_{K_2}(\Delta t)^2\Vert\nabla e_{\textbf{\textit{u}}}^{n+1}\Vert_{L^2}^2+C(\Delta t)^2\Vert e_\phi^n\Vert_{L^2}^2+C(\Delta t)^2\Vert \nabla e_\mu^{n+1}\Vert_{L^2}^2+C(\Delta t)^2\Vert T_{\textbf{\textit{u}}}^{n+1}\Vert_{L^2}^2
	\end{split}
\end{equation}

Summing (\ref{the4.17}) over $n$ and using (\ref{the4.14}), we arrive at
\begin{equation}\label{the4.18}
\begin{split}
\Delta t\Vert\nabla e_{\textbf{\textit{u}}}^{n+1}\Vert_{L^2}^2+\sum_{k=0}^{n}\Vert e_{\textbf{\textit{u}}}^{k+1}-e_{\textbf{\textit{u}}}^{k}\Vert_{L^2}^2\leq C\left(\Delta t\right)^3.
\end{split}
\end{equation}

By the regularity of the solution to the Stokes problem (see \cite{temam2024navier}), we can get
\begin{equation}\label{the4.19}
\begin{split}
&\Delta t\sum_{k=0}^{n} \left(\Vert e_{\textbf{\textit{u}}}^{k+1}\Vert_{H^2}^2+\Vert e_p^{k+1}\Vert_{H^1}^2\right)\\
&\leq C\Delta t\sum_{k=0}^{n}\left(\left\Vert \frac{e_{\textbf{\textit{u}}}^{k+1}- e_{\textbf{\textit{u}}}^{k}}{\Delta t} \right\Vert_{L^2}^2+\Vert e_{\textbf{\textit{u}}}^{k}\Vert_{H^1}^2+\Vert  e_\phi^{k}\Vert_{L^2}^2+\Vert  \nabla e_\mu^{k+1}\Vert_{L^2}^2+\Vert  \nabla e_c^{k+1}\Vert_{L^2}^2+\Vert T_{\textbf{\textit{u}}}^{k+1}\Vert_{L^2}^2\right)\\
&\leq C\left(\Delta t\right)^2.
\end{split}
\end{equation}

For $k=0$, we have $\Vert {\textbf{\textit{u}}}^0\Vert_{L^\infty}\leq K_2$. For $\Vert {\textbf{\textit{u}}}^k\Vert_{L^\infty}\leq K_2$ holds true for $1\leq k\leq n$, we can see from (\ref{the4.19}) that $\Vert e_{{\textbf{\textit{u}}}}^{n+1}\Vert_{H^2}\leq C\left(\Delta t\right)^{\frac{1}{2}}$. Therefore, 
\begin{equation}\label{the4.20}
	\begin{split}
\Vert {\textbf{\textit{u}}}^{n+1}\Vert_{L^\infty}^2\leq C_3\left(\Delta t\right)^{\frac{1}{2}}+\Vert {\textbf{\textit{u}}}(t^{n+1})\Vert_{L^\infty}.
	\end{split}
\end{equation}
We complete the induction by setting $\Delta t\leq \left(C_3\right)^{-2}(\Vert {\textbf{\textit{u}}}(0)\Vert_{L^\infty}+1)^2$.

By setting $\Delta t< \text{min}\{\left(C_1\right)^{-1}, \left(C_2\right)^{-2}, \left(C_3\right)^{-2}(\Vert {\textbf{\textit{u}}}(0)\Vert_{L^\infty}+1)^2\}$, we thus complete the proof.
\end{proof}
\end{theorem}

\section{Numerical simulations}
\label{section5}
In this section, we first perform a numerical test to verify the accuracy and stability of the scheme. Then, we present two numerical examples to illustrate the performance of the proposed model. Finally, we conducted simulations on straight and bifurcated vascular structures to investigate risks associated with microaneurysm development. Without specific needs, we set $Re=Pe=1$, $\mathcal{M}=0.05$, $D_2=D_3=1$, $M=N=1$, $k=1$, $\epsilon=0.04$, $\lambda=1+\lambda_0c_3$, $S=4$, $a_1=a_2=1$ and $a_3=2$ for simulations. The mesh size is fixed $h=\frac{1}{64}$.

\subsection{Convergence and stability tests}
In this test, we first verify the accuracy in time of the proposed scheme by refining the time step. We set the final time $T=0.1$, the domain $\Omega=[0,1]^2$, $\lambda_0=0.05$, $Q=q_0c_1$ and $q=\frac{dQ}{dc_1}=q_0$ with $q_0=1$. The following initial conditions are used: 
\begin{equation}\label{the5.1}
	\begin{split}
&\textbf{\textit{u}}(0,\textbf{x})=0, \quad \phi(0,\textbf{x})=\text{tanh}\left(\frac{0.20-\sqrt{(x-0.5)^2+(y-0.5)^2}}{\sqrt{2}\epsilon}\right),\\
&c_1(0,\textbf{x})=2.0, \quad c_2(0,\textbf{x})=2.0, \quad c_3(0,\textbf{x})=0.2.
	\end{split}
\end{equation}
All variables along $x$-direction are periodic. On the top and bottom boundaries, zero Neumann and no-slip boundary conditions are used for scalar variables and velocities, respectively. 

\begin{figure}[H]
	\centering
	\includegraphics[scale=0.58]{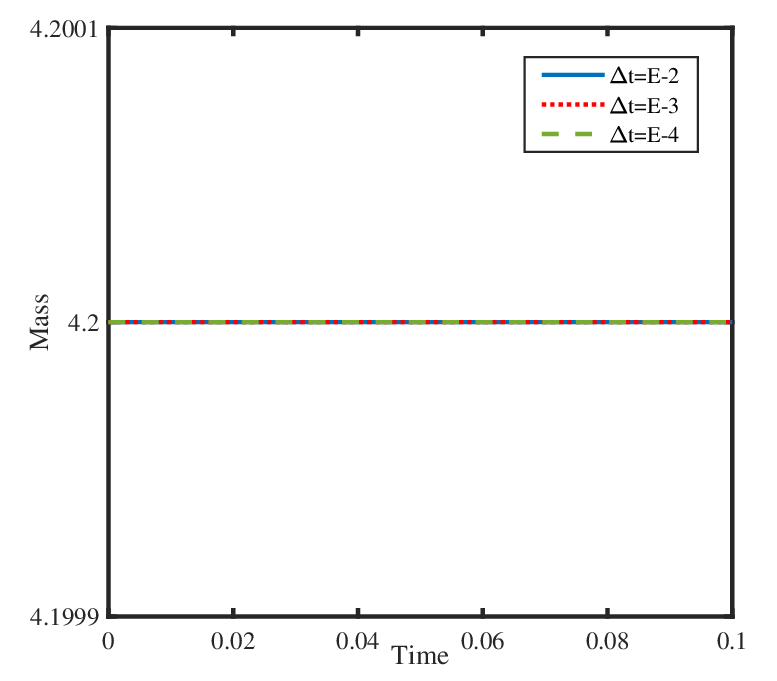}
	\includegraphics[scale=0.58]{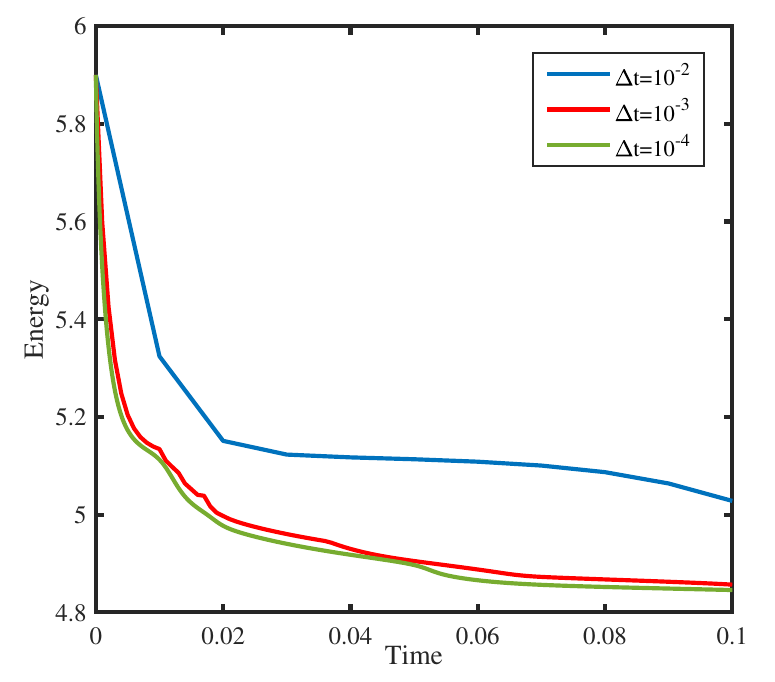}
	\caption{Evolution of total mass and total energy with $\Delta t=10^{-2}$, $\Delta t=10^{-3}$ and $\Delta t=10^{-4}$.}
	\label{MassEnergy}
\end{figure}

We use a finer time step $\Delta t^r=10^{-4}$ to perform the test until $T=0.1$ as the reference solution. A set of decreasing time steps is used to test the convergence rate. As presented in Table $\ref{table1}$, the results confirm the temporal accuracy of the proposed scheme. We further test the stability of the scheme. Figure \ref{MassEnergy} shows the temporal evolution of the total mass and total energy of the system under different time steps, which demonstrates the proposed scheme is mass conserved and satisfies the energy dissipation law.

\begin{table}[H]
\captionsetup{justification=centering} 
\caption{Convergence rates of scalar variables and velocities at $T=0.1$.}\label{table1}
\centering
\scalebox{0.98}{
\resizebox{\linewidth}{!}{
\begin{tabular}{c c c c c c c c c c c} \hline  		
\scriptsize   $\Delta t$   & \scriptsize $\Vert e_\textbf{\textit{u}}\Vert_{L^2}$ & \scriptsize   Order & \scriptsize  $\Vert e_\phi\Vert_{L^2}$ & \scriptsize  Order & \scriptsize $\Vert e_{c_1}\Vert_{L^2}$ & \scriptsize  Order & \scriptsize $\Vert e_{c_2}\Vert_{L^2}$ & \scriptsize Order & \scriptsize $\Vert e_{c_3}\Vert_{L^2}$ & \scriptsize Order\\ \hline
\scriptsize  $4\cdot 10^{-3}$		& \scriptsize  2.57E-2 &  & \scriptsize 8.48E-3  & &\scriptsize 1.80E-2 & &\scriptsize 1.44E-2 & & \scriptsize 1.40E-2 &\\

\scriptsize  $2\cdot 10^{-3}$		& \scriptsize  1.65E-2 & \scriptsize 0.64 & \scriptsize 4.87E-3  & \scriptsize 0.80 &\scriptsize  1.15E-2 & \scriptsize 0.65 &\scriptsize 8.99E-3 & \scriptsize 0.68 & \scriptsize  8.86E-3 & \scriptsize 0.66\\

\scriptsize  $1\cdot 10^{-3}$		& \scriptsize   9.36E-3  & \scriptsize 0.82 & \scriptsize 2.54E-3  & \scriptsize 0.94 &\scriptsize 6.54E-3 &  \scriptsize 0.82 &\scriptsize 5.02E-3 & \scriptsize 0.84 & \scriptsize 5.05E-3 & \scriptsize 0.81\\

\scriptsize  $5\cdot 10^{-4}$		& \scriptsize  4.68E-3 & \scriptsize 1.00 & \scriptsize 1.25E-3  & \scriptsize 1.02 &\scriptsize 3.26E-3 & \scriptsize 1.00 &\scriptsize 2.53E-3 & \scriptsize 0.99 & \scriptsize 2.54E-3 & \scriptsize 0.99 \\  \hline
		\end{tabular}}
	}
\end{table}

\subsection{Adsorption test under shear flow}
\begin{figure}[H]
	\centering
	\subfigure[$t=0.03$]{
		\includegraphics[scale=0.10]{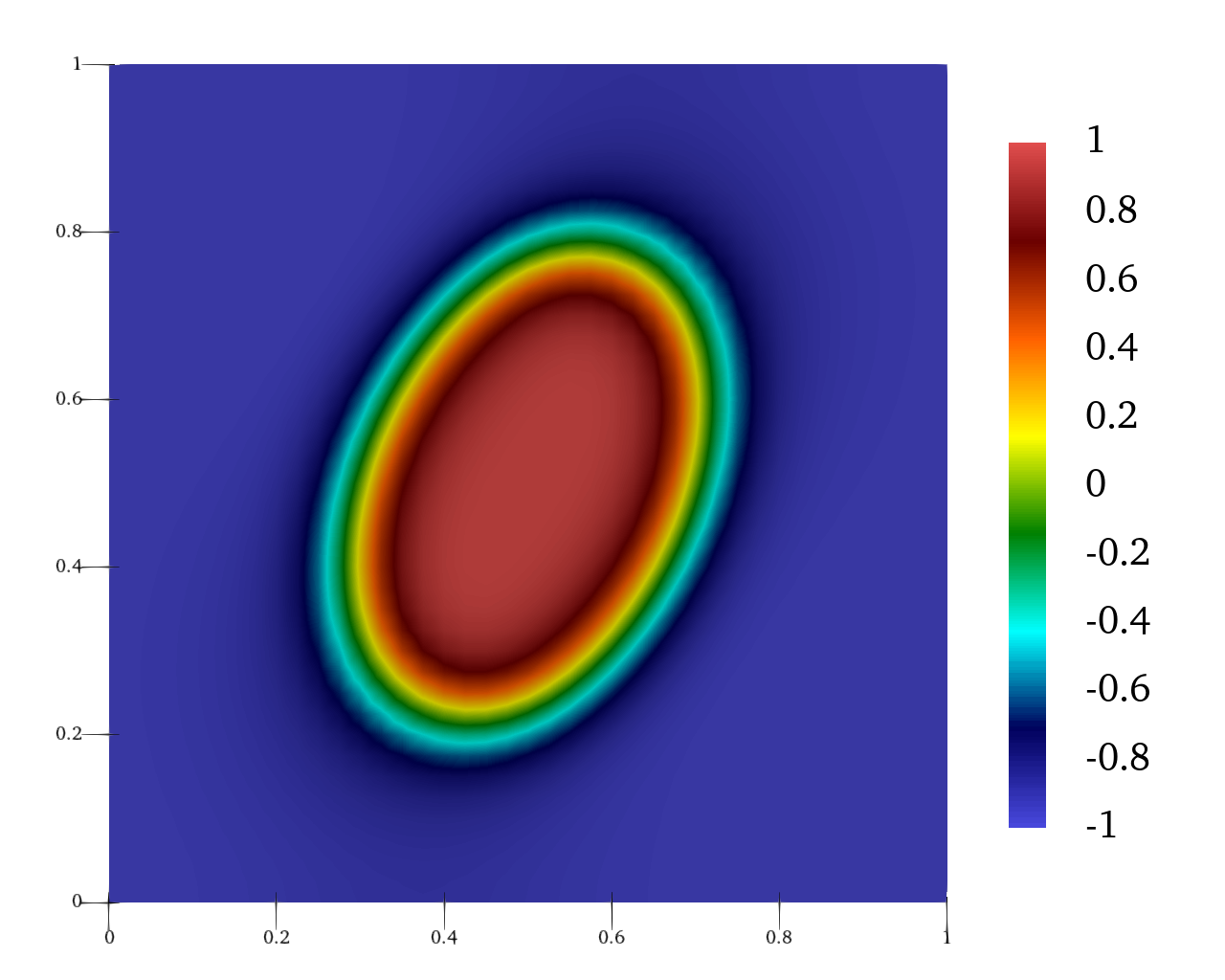}
	}
	\subfigure[$t=0.07$]{
		\includegraphics[scale=0.10]{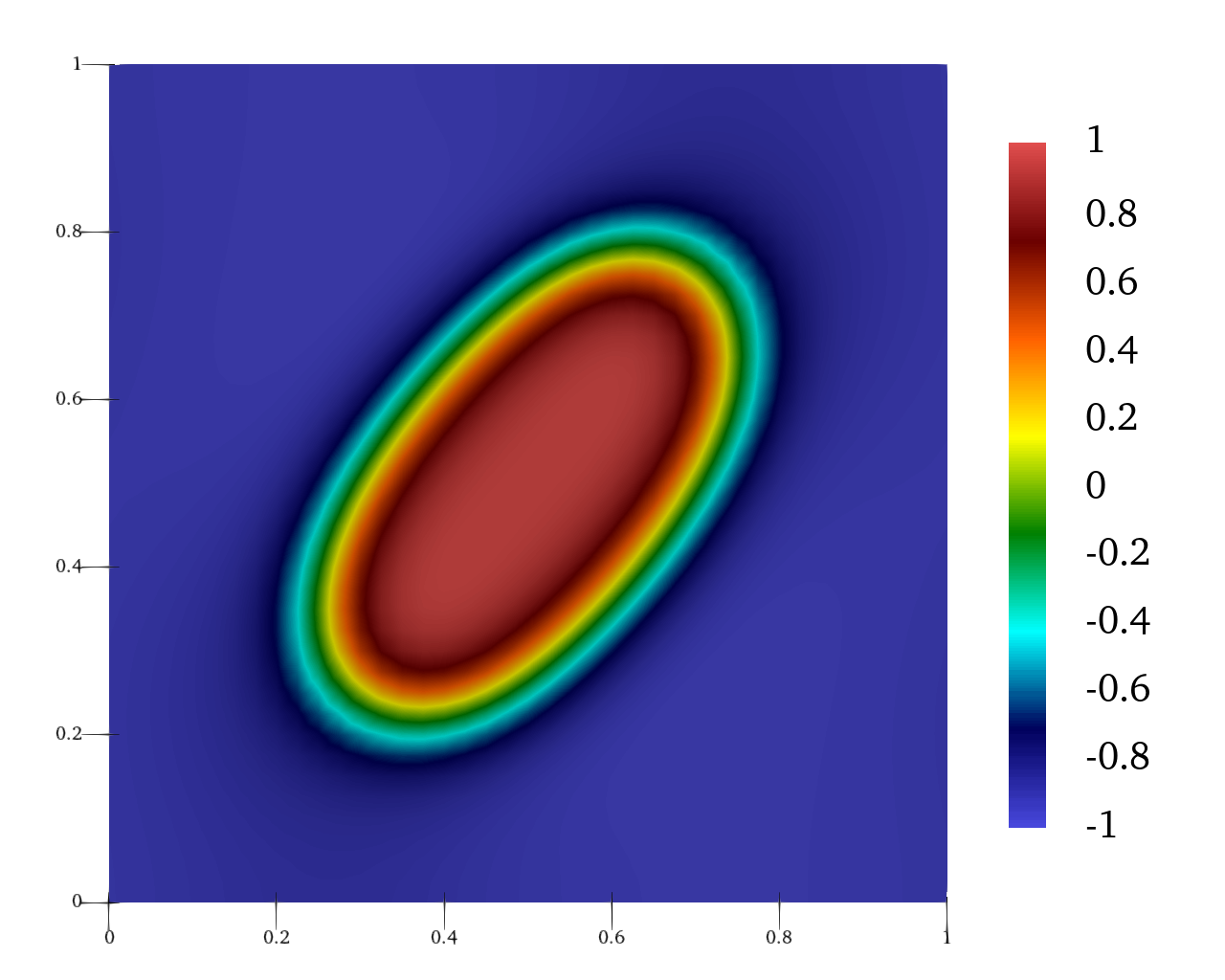}
	}
	\subfigure[$t=0.10$]{
		\includegraphics[scale=0.10]{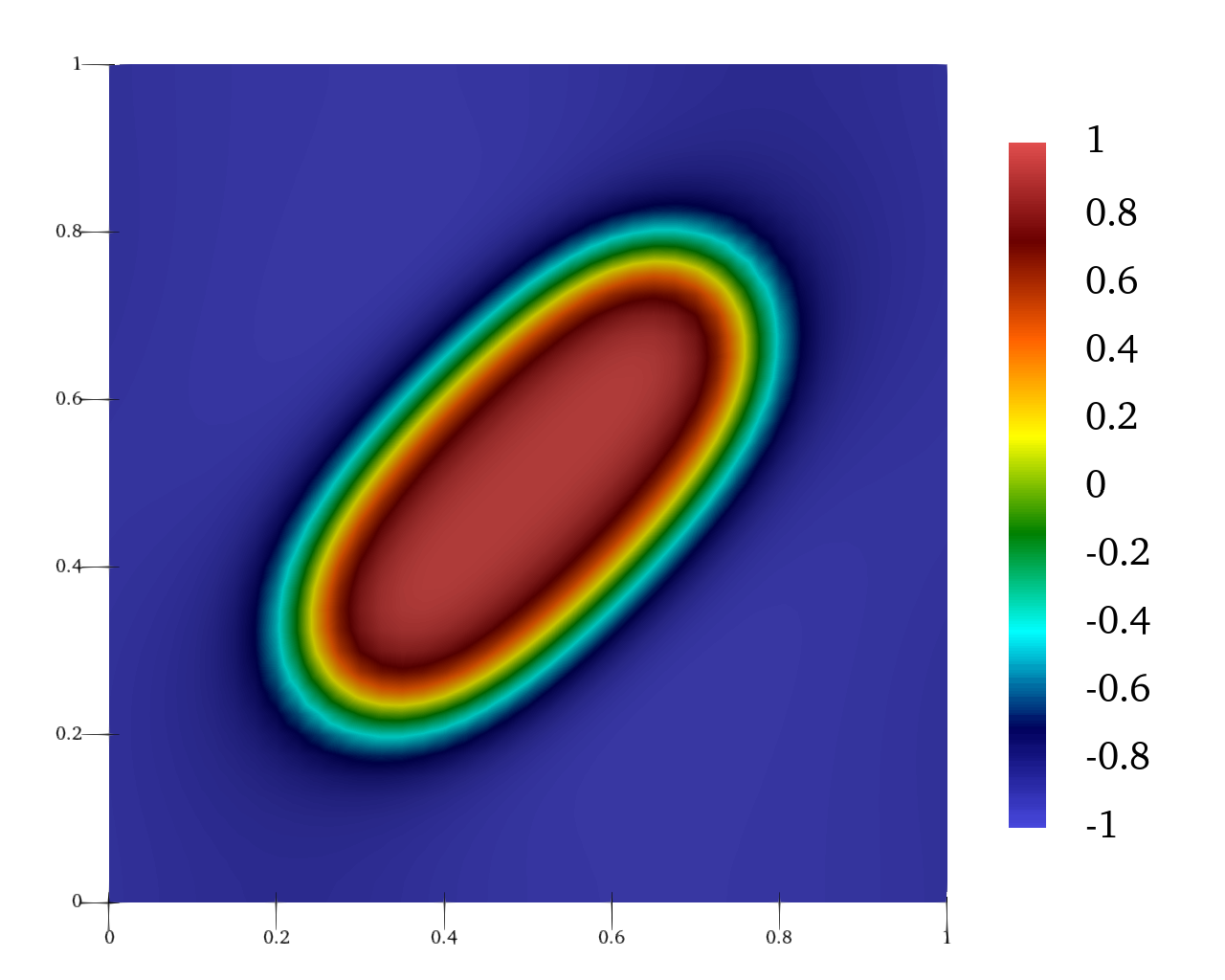}
	}
	
	\caption{Snapshots of $\phi$ at different times with $M=N=0$.}
	\label{figure2}
\end{figure}

\begin{figure}[H]
	\centering
	\subfigure[$t=0.03$]{
		\includegraphics[scale=0.101]{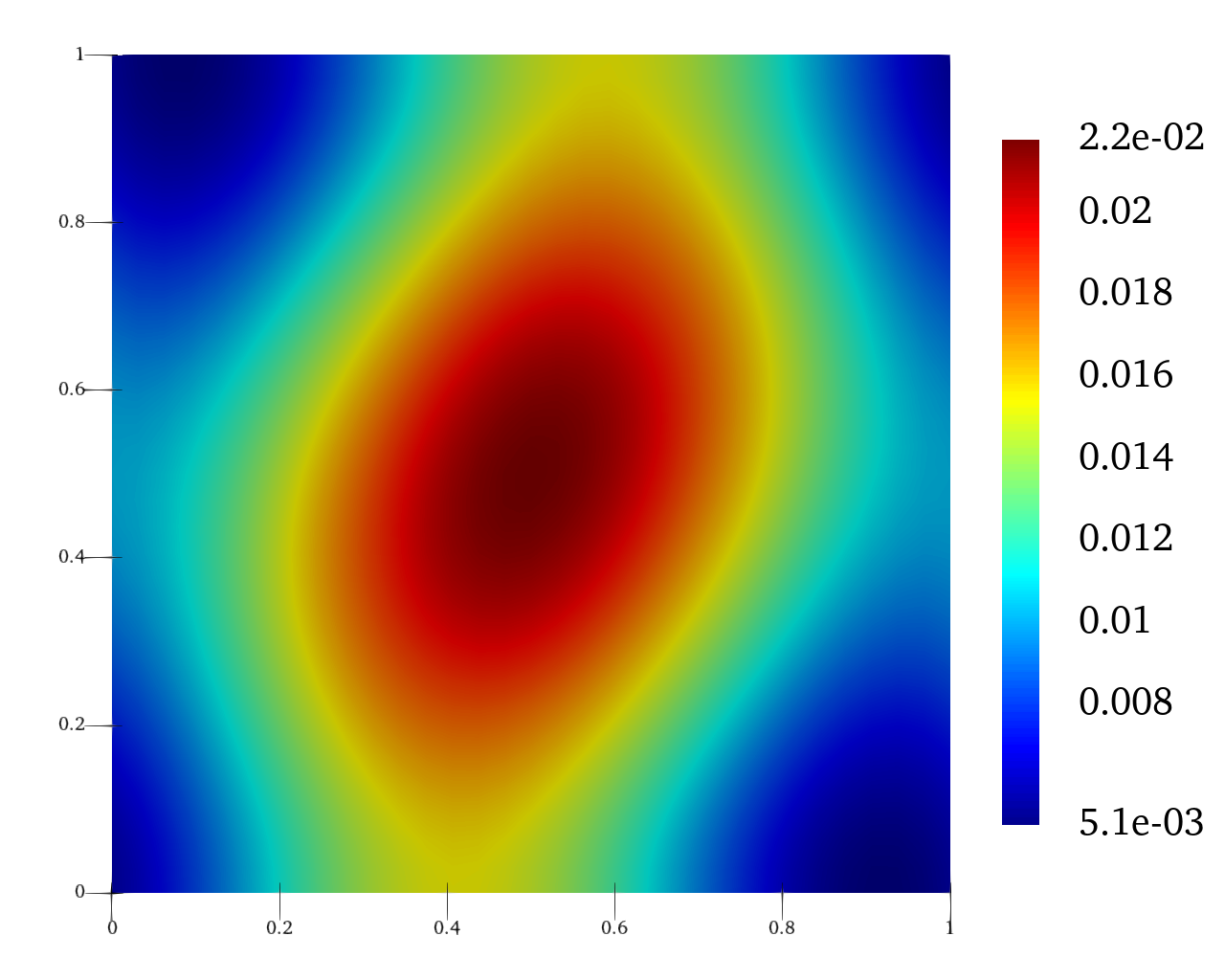}
	}
	\subfigure[$t=0.07$]{
		\includegraphics[scale=0.10]{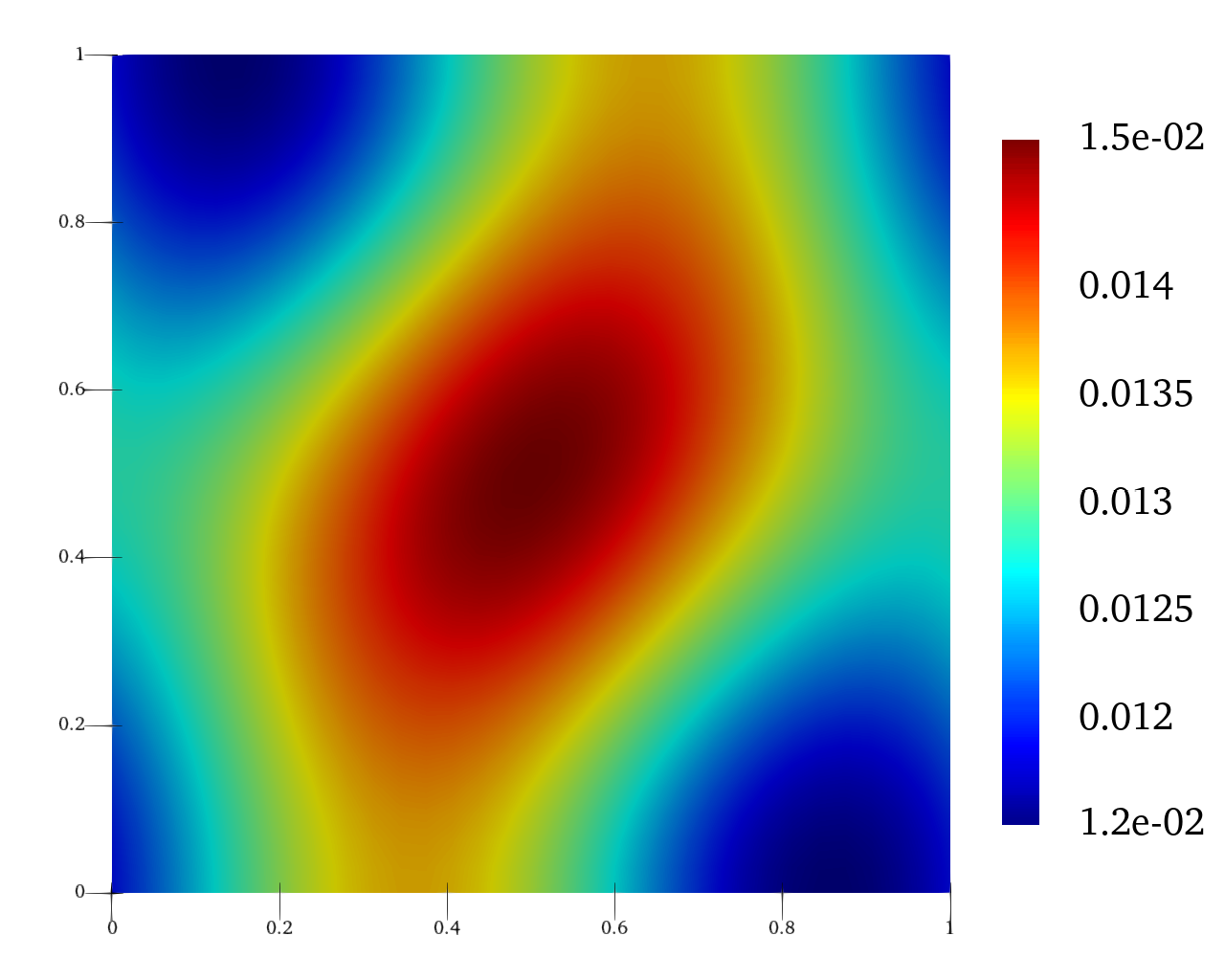}
	}
	\subfigure[$t=0.10$]{
		\includegraphics[scale=0.10]{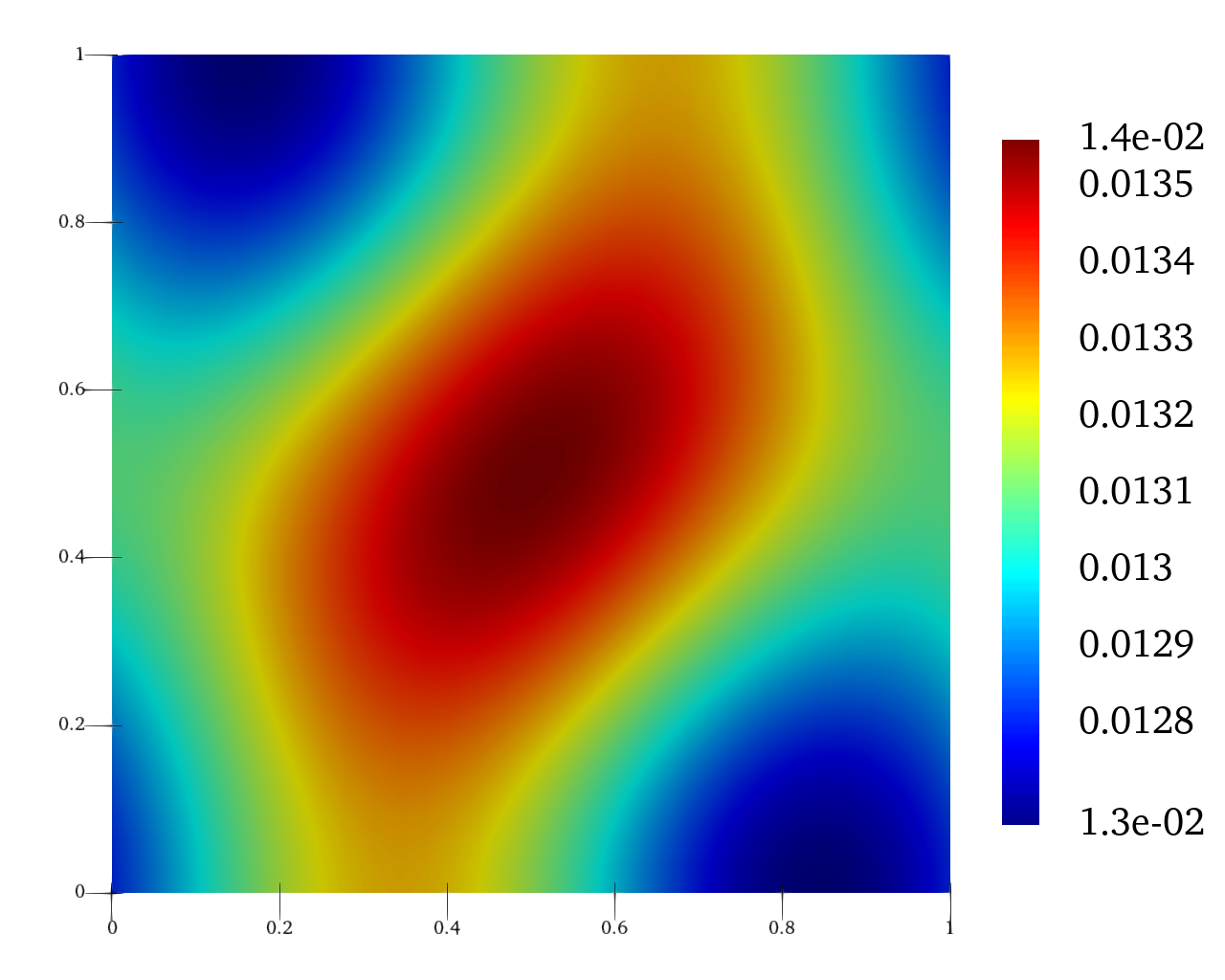}
	}
	
	\caption{Snapshots of $Y_2$ at different times with $M=N=0$.}
	\label{figure3}
\end{figure}

\begin{figure}[H]
	\centering
	\subfigure[$t=0.03$]{
		\includegraphics[scale=0.101]{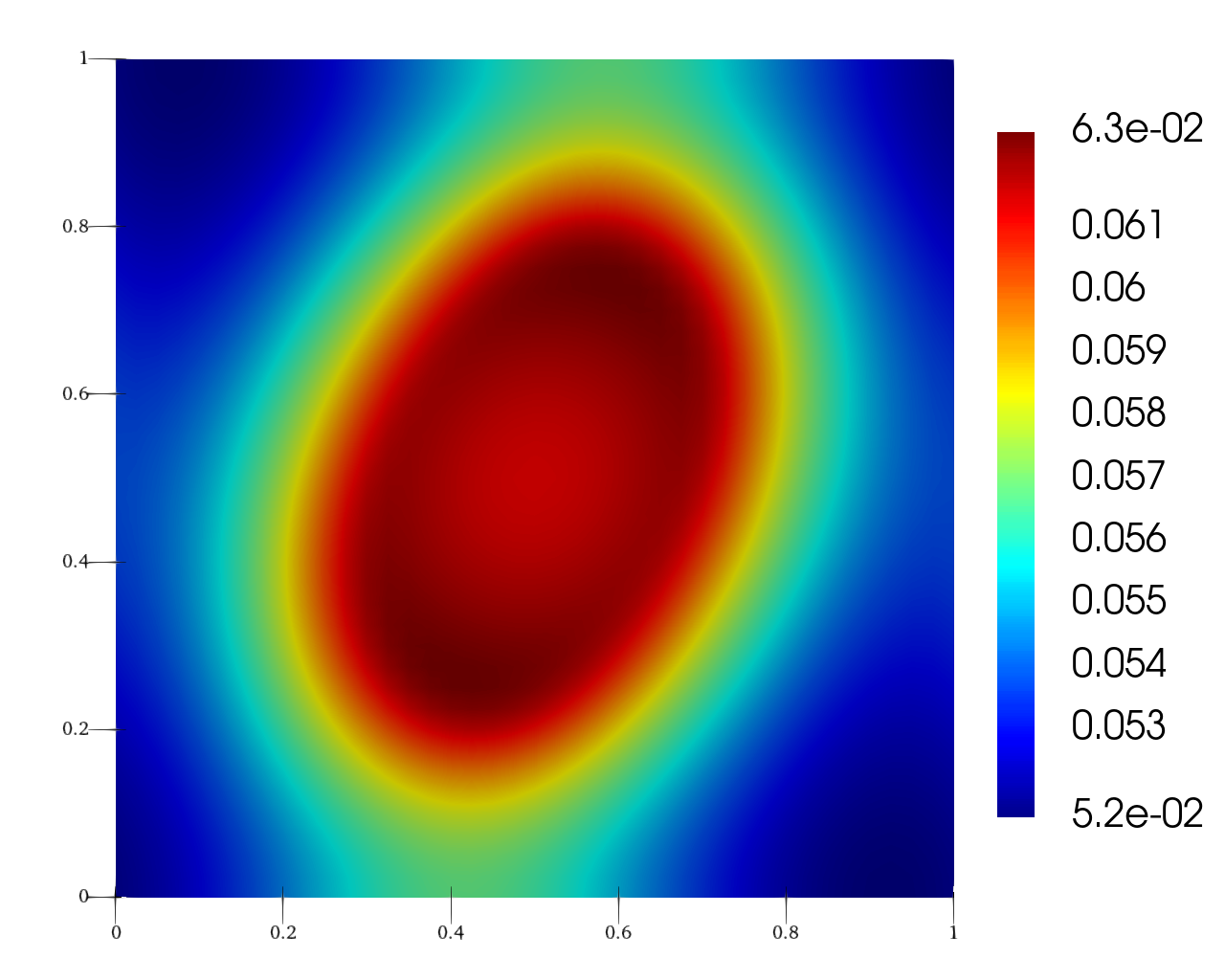}
	}
	\subfigure[$t=0.07$]{
		\includegraphics[scale=0.10]{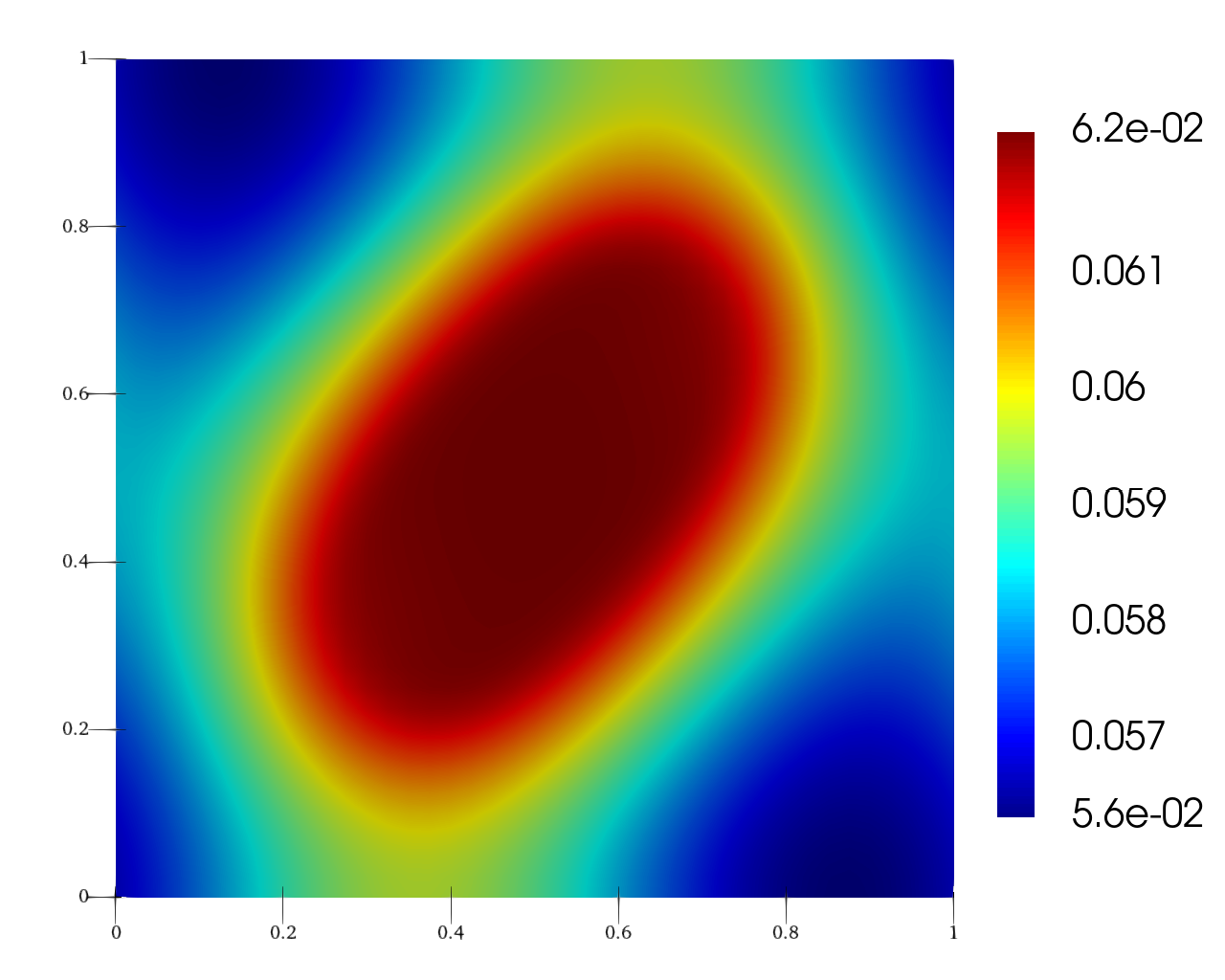}
	}
	\subfigure[$t=0.10$]{
		\includegraphics[scale=0.10]{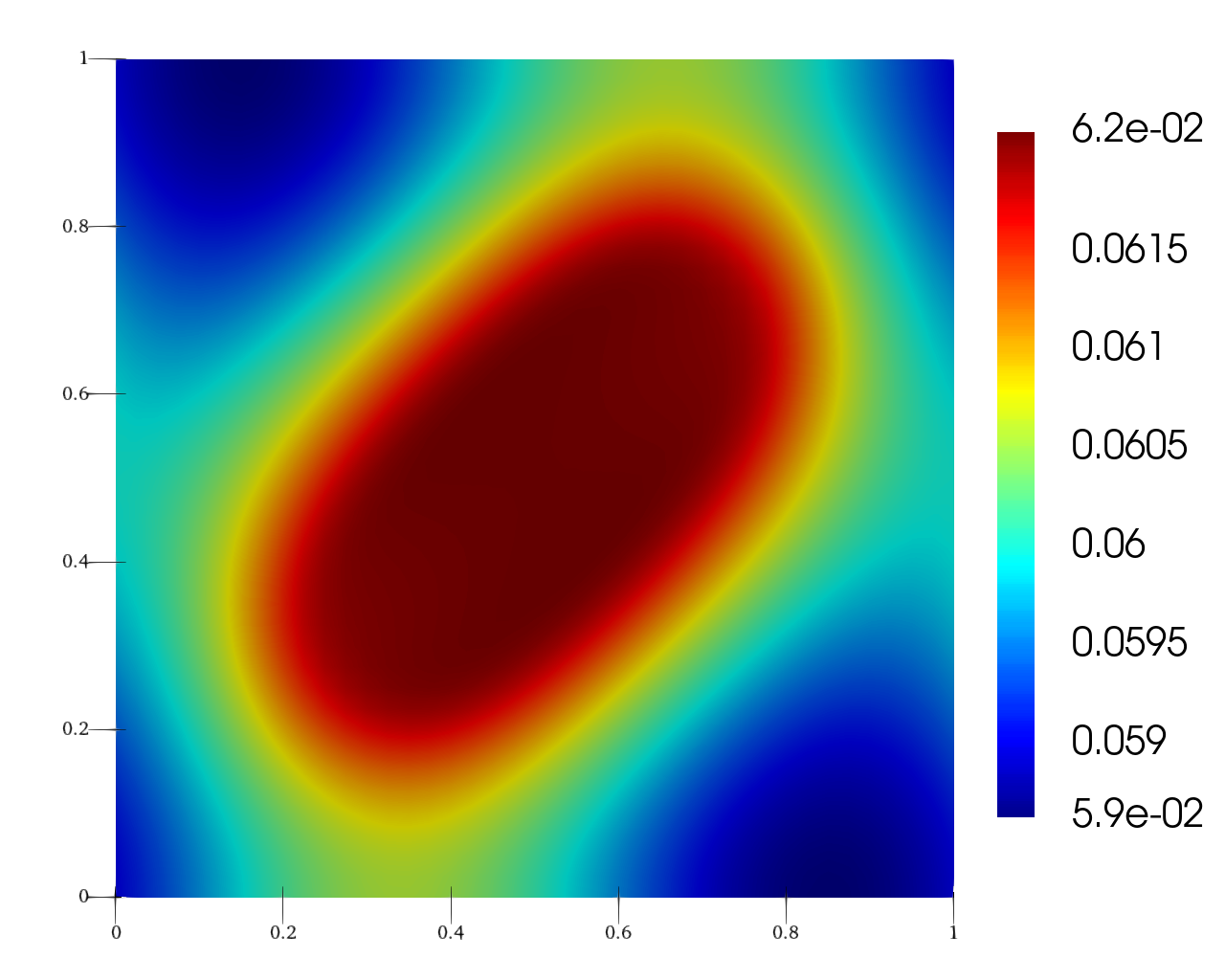}
	}
	
	\caption{Snapshots of $Y_3$ at different times with $M=N=0$.}
	\label{figure4}
\end{figure}

\begin{figure}[H]
	\centering
	\subfigure[$t=0.03$]{
		\includegraphics[scale=0.10]{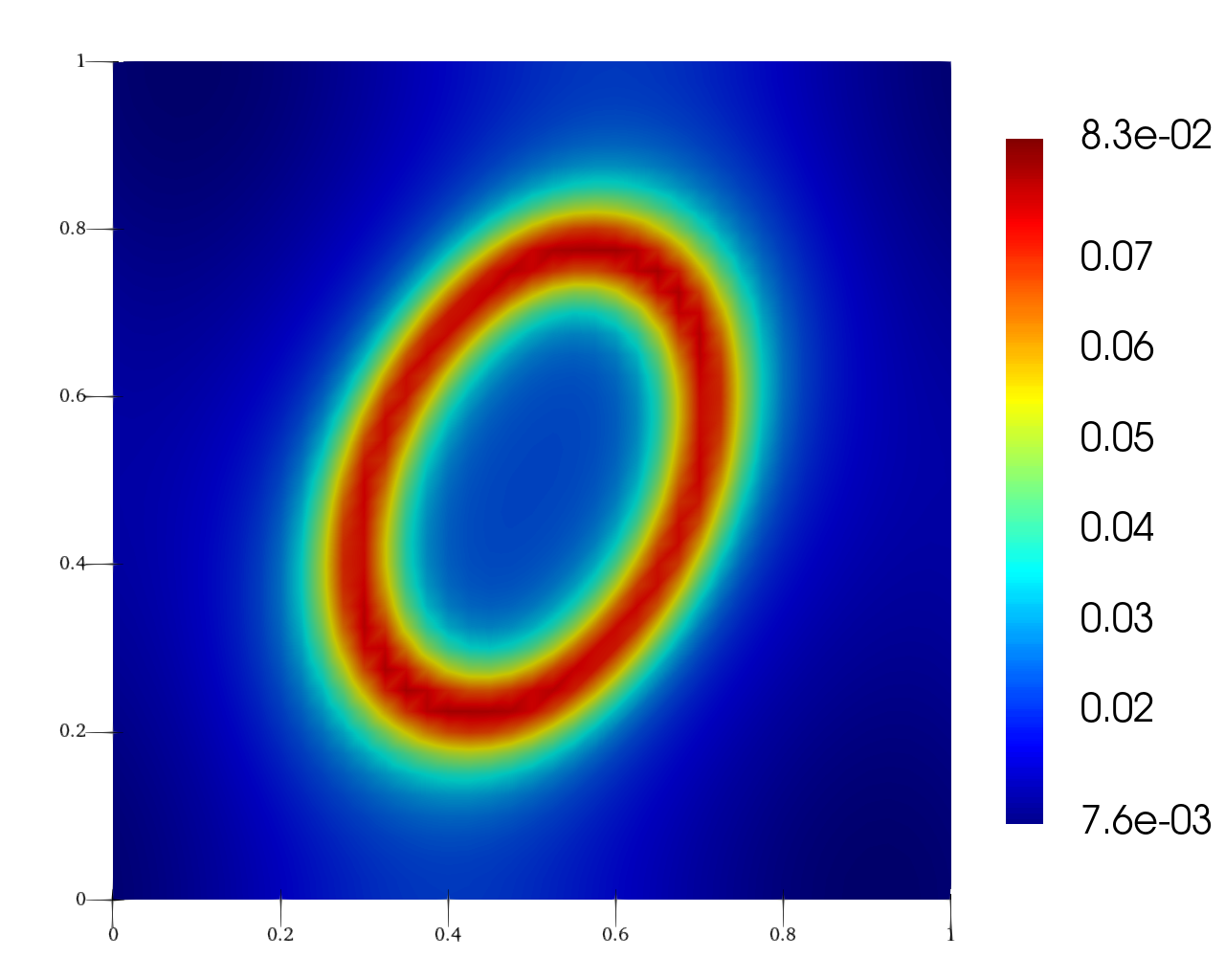}
	}
	\subfigure[$t=0.07$]{
		\includegraphics[scale=0.10]{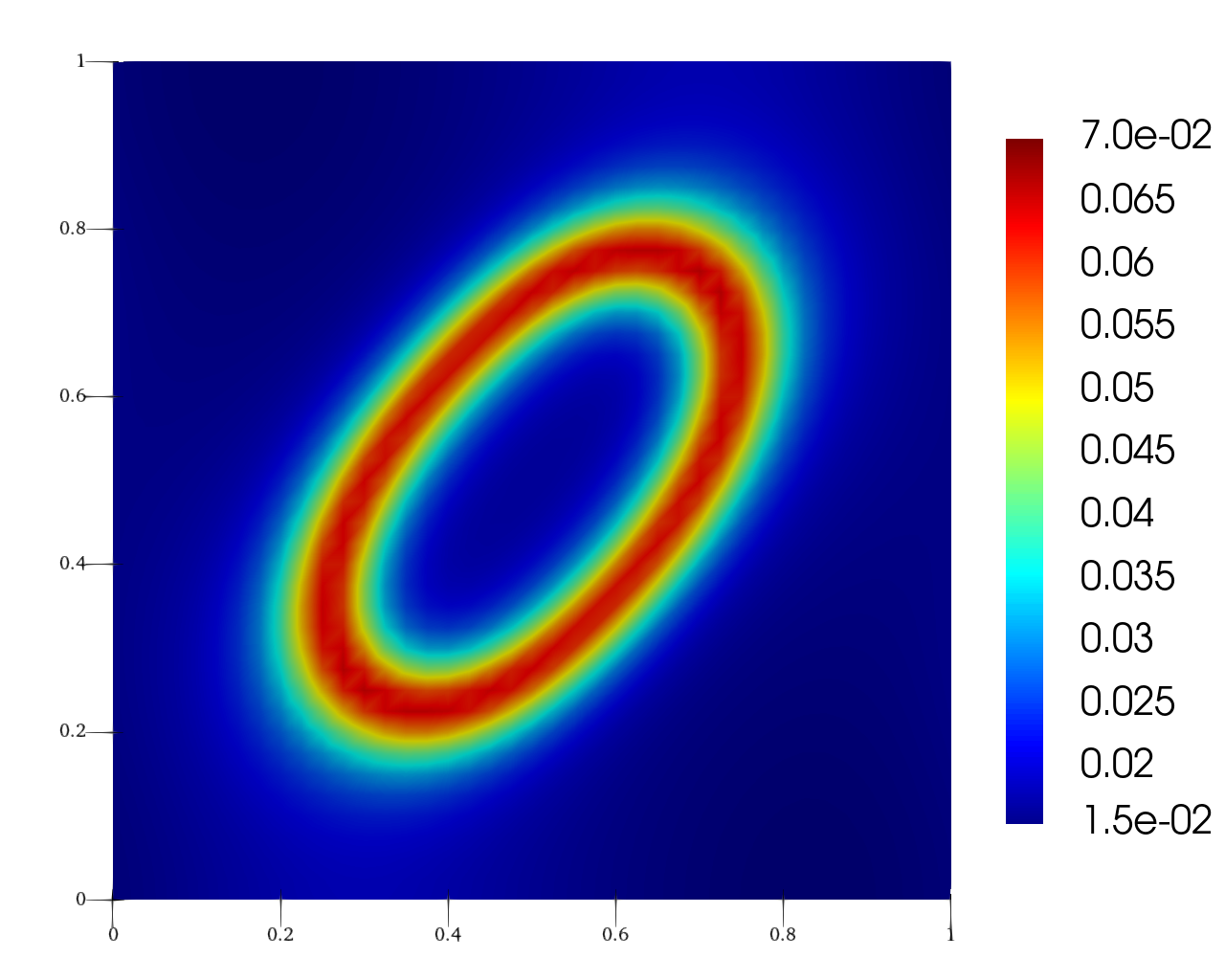}
	}
	\subfigure[$t=0.10$]{
		\includegraphics[scale=0.10]{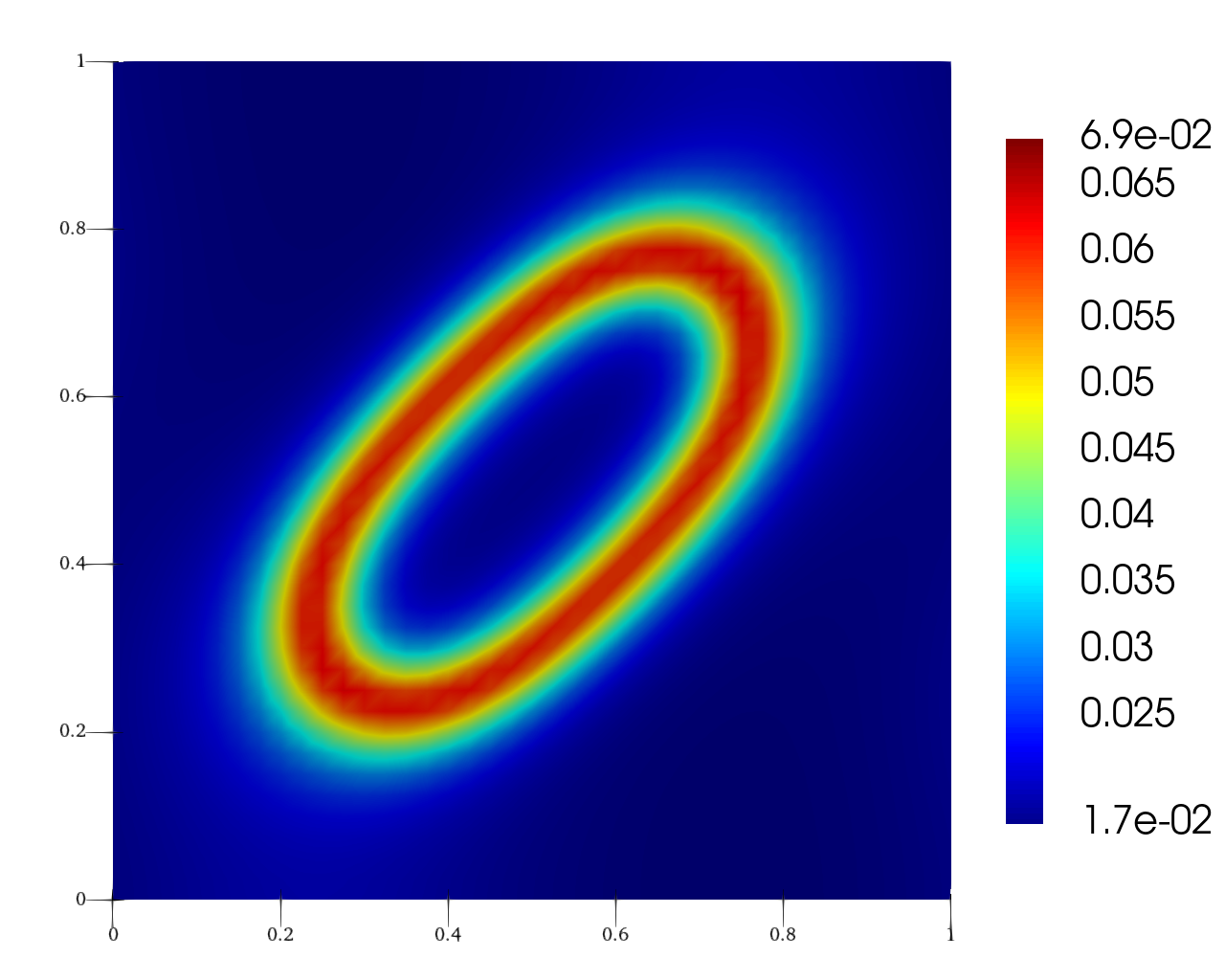}
	}
	
	\caption{Snapshots of $Y_2$ at different times with $M=N=1$.}
	\label{figure5}
\end{figure}

\begin{figure}[H]
	\centering
	\subfigure[$t=0.03$]{
		\includegraphics[scale=0.10]{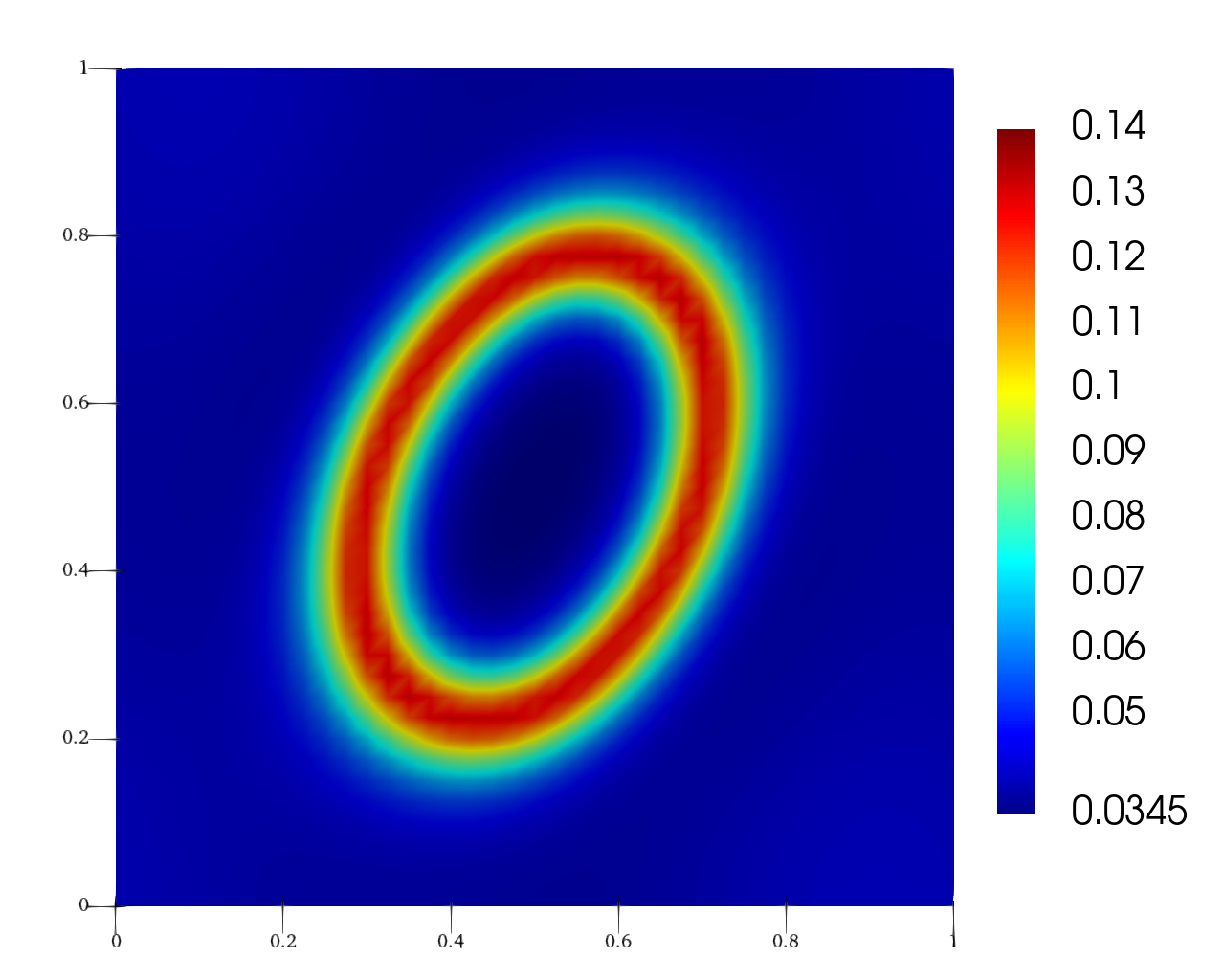}
	}
	\subfigure[$t=0.07$]{
		\includegraphics[scale=0.10]{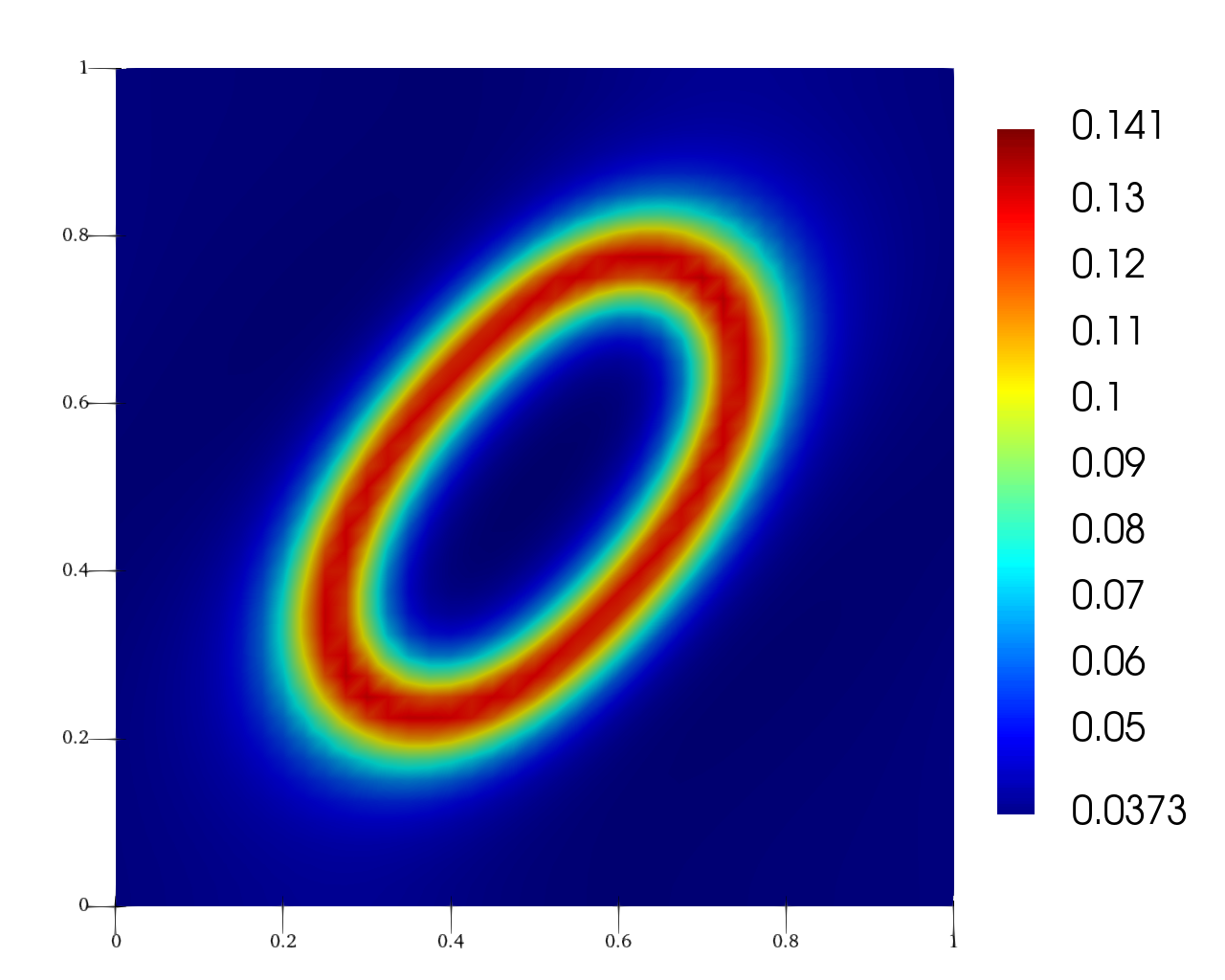}
	}
	\subfigure[$t=0.10$]{
		\includegraphics[scale=0.10]{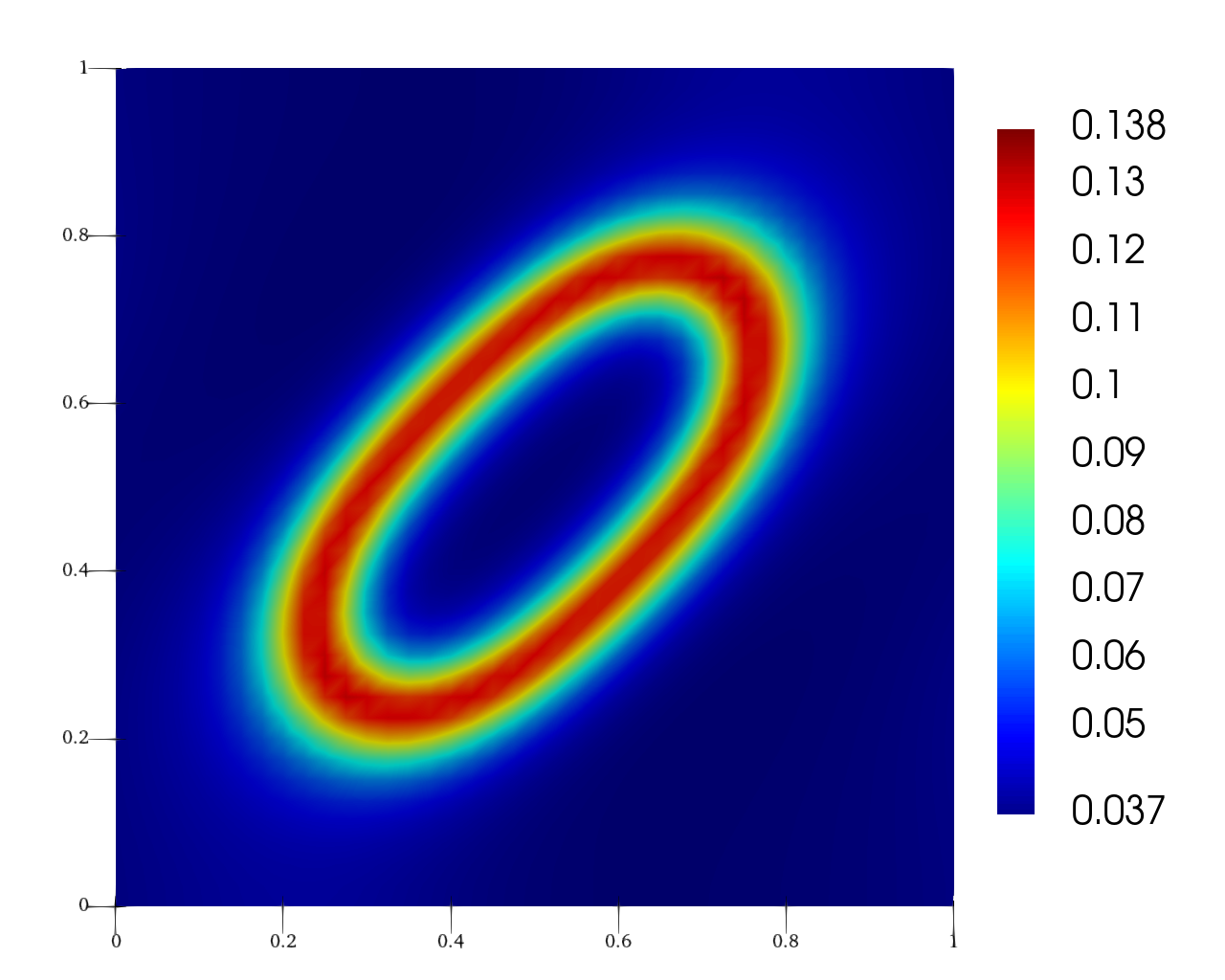}
	}
	
	\caption{Snapshots of $Y_3$ at different times with $M=N=1$.}
	\label{figure6}
\end{figure}

To investigate the effects of parameters $M$ and $N$ on adsorption characteristics, we set $M=N=0$ and $M=N=1$ to compare the adsorption performance of the interface on $Y_2$ and $Y_3$. The initial conditions are defined to be
\begin{equation}\label{the5.2}
\begin{split}
&\textbf{\textit{u}}(0,\textbf{x})=\left(10(y-0.5),0\right)^{\dagger}, \quad \phi(0,\textbf{x})=\text{tanh}\left(\frac{0.20-\sqrt{(x-0.5)^2+\frac{(y-0.5)^2}{2}}}{\sqrt{2}\epsilon}\right),\\
&c_1(0,\textbf{x})=0.2, \quad c_2(0,\textbf{x})=0.2\left(1-\left|\phi(0,\textbf{x})\right|\right), \quad c_3(0,\textbf{x})=0.
\end{split}
\end{equation}
We use $\Delta t=10^{-3}$, $\lambda_0=0$ and $q_0=1$ for the simulation. In the absence of adhesion effects, Figure \ref{figure2} shows that the initial elliptical shape is stretched by the shear flow until a balance is achieved between the surface tension and the applied velocity.
The corresponding distribution patterns of $Y_2$ and $Y_3$ are shown in Figure \ref{figure3} and \ref{figure4}, both exhibiting diffusive characteristics. When $M=N=1$, it can be seen from Figure
\ref{figure5} and \ref{figure6} that the interface exhibits adhesion effects, with $Y_2$ and the generated $Y_3$ both confined to the interface.


\subsection{Arterial microaneurysm}
In this subsection, we apply the proposed model to the problem of microaneurysm formation. The phase-field function $\phi$ is regarded as the order parameter to distinguish between the two phases inside and outside the blood vessel. The diffuse interface, in which $\phi\in(-1,1)$, represents the vascular wall. $Y_1$, $Y_2$ and $Y_3$ can be represented as glucose, structural proteins in the vascular wall, and advanced glycation end products (AGEs), respectively. While early glycation might be reversible, sustained AGEs accumulation leads to vascular stiffening, loss of elasticity, and structural degradation, which contribute to microaneurysm formation. In the following, we simulate microvessels with straight and bifurcated structures to illustrate the risk of microaneurysm formation.

Research has shown that AGEs disrupt endothelial junction integrity by downregulating tight junction proteins, increasing vascular permeability and promoting glucose extravasation into interstitial tissue, particularly under hyperglycemic conditions \cite{hudson2018targeting}. Therefore, we set $q=1+2 c_3$ in the following simulations. The settings for the other parameters are $D^{+}=1$, $D^{-}=0.5$ and $\lambda_0=0.5$.

\subsubsection{Straight microvessel} 
Straight-structured vessels are commonly observed in the microvascular system. We set the final time $T=3$ and the domain $\Omega=[0,2]^2$. The following initial and boundary conditions for $\phi$ and $\textbf{\textit{u}}$ are used:
\begin{equation}\nonumber
	\begin{split}
&\phi(0,\textbf{x})=0.95 \ (0\leq x\leq 2, 0.8\leq y\leq 1.2), \ \phi(0,\textbf{x})=-0.95 \ (\text{otherwise}), \\
&\nabla \phi\cdot \textbf{\textit{n}}|_{\partial \Omega}=\nabla \mu_\phi\cdot \textbf{\textit{n}}|_{\partial \Omega}=0.\\
&\textbf{\textit{u}}(0,\textbf{x})=0, \ \textbf{\textit{u}}|_{\{x=0, \ 0.8\leq y\leq 1.2\}}=25\left(y-0.8\right)\left(y-1.2\right), \ p|_{\{x=2, \ 0.8\leq y\leq 1.2\}}=0.
	\end{split}
\end{equation}
The other boundaries of $\textbf{\textit{u}}$ are no-slip. For $Y_1$ and $Y_2$, we set a local high concentration at the initial for the reaction to generate $Y_3$:
\begin{equation}\nonumber
	\begin{split}
&c_1(0,\textbf{x})=c_2(0,\textbf{x})=\text{tanh}\left(\frac{0.08-\sqrt{(x-1)^2+(y-1.2)^2}}{\sqrt{2}\epsilon}\right)+1.2, \  c_3(0,\textbf{x})=0, \\
& \nabla c_2\cdot \textbf{\textit{n}}|_{\partial \Omega}=\nabla c_3\cdot \textbf{\textit{n}}|_{\partial \Omega}=0.
	\end{split}
\end{equation}
On the left boundary $x=0$, the boundary condition of $c_1$ is set to $c_1|_{\{x=0, \ 0.8\leq y\leq 1.2\}}=2$, and on other boundaries it is the homogeneous Neumann condition.

\begin{figure}[H]
	\centering
	\subfigure[$t=1.0$]{
		\includegraphics[scale=0.10]{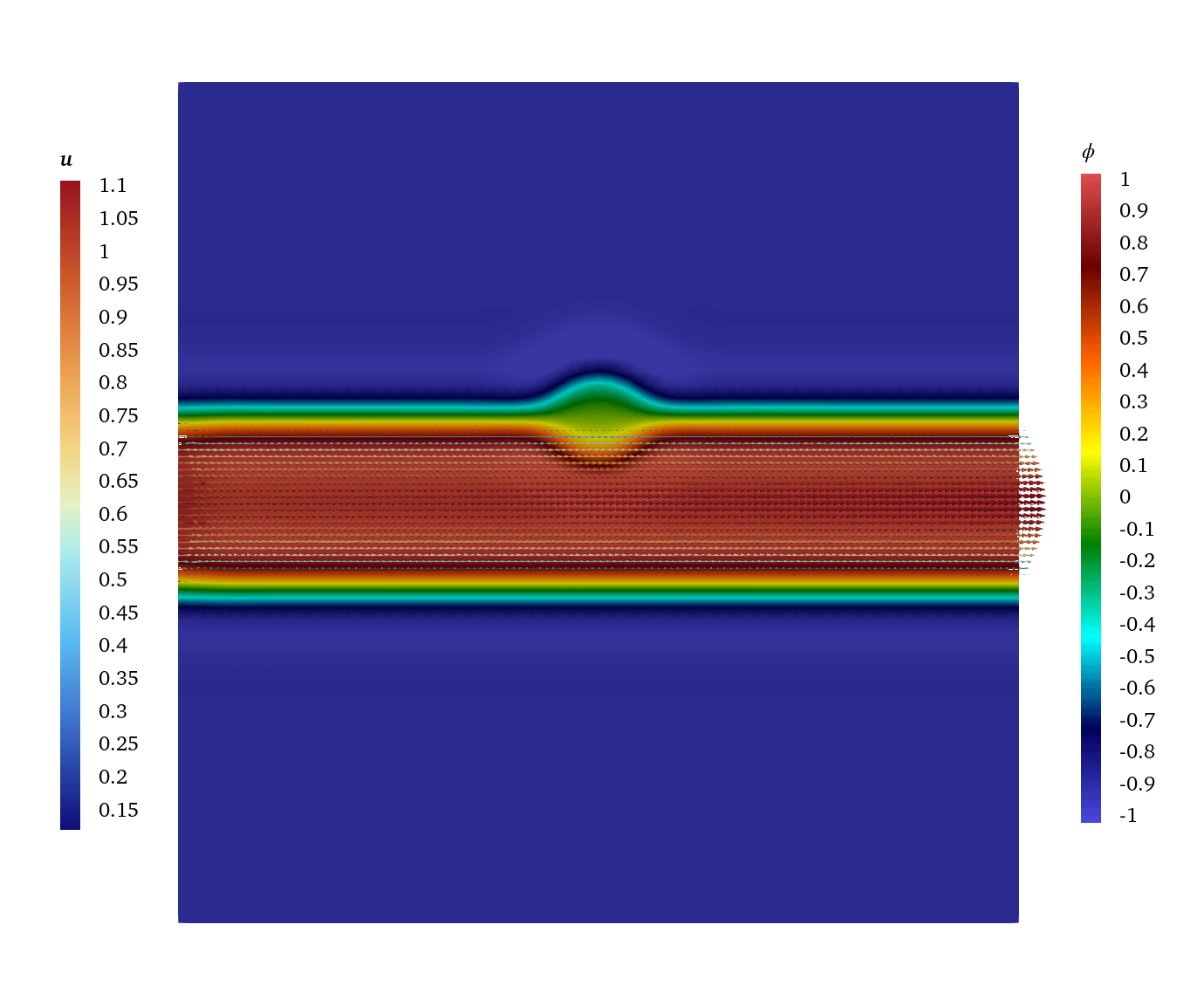}
	}
	\subfigure[$t=2.0$]{
		\includegraphics[scale=0.10]{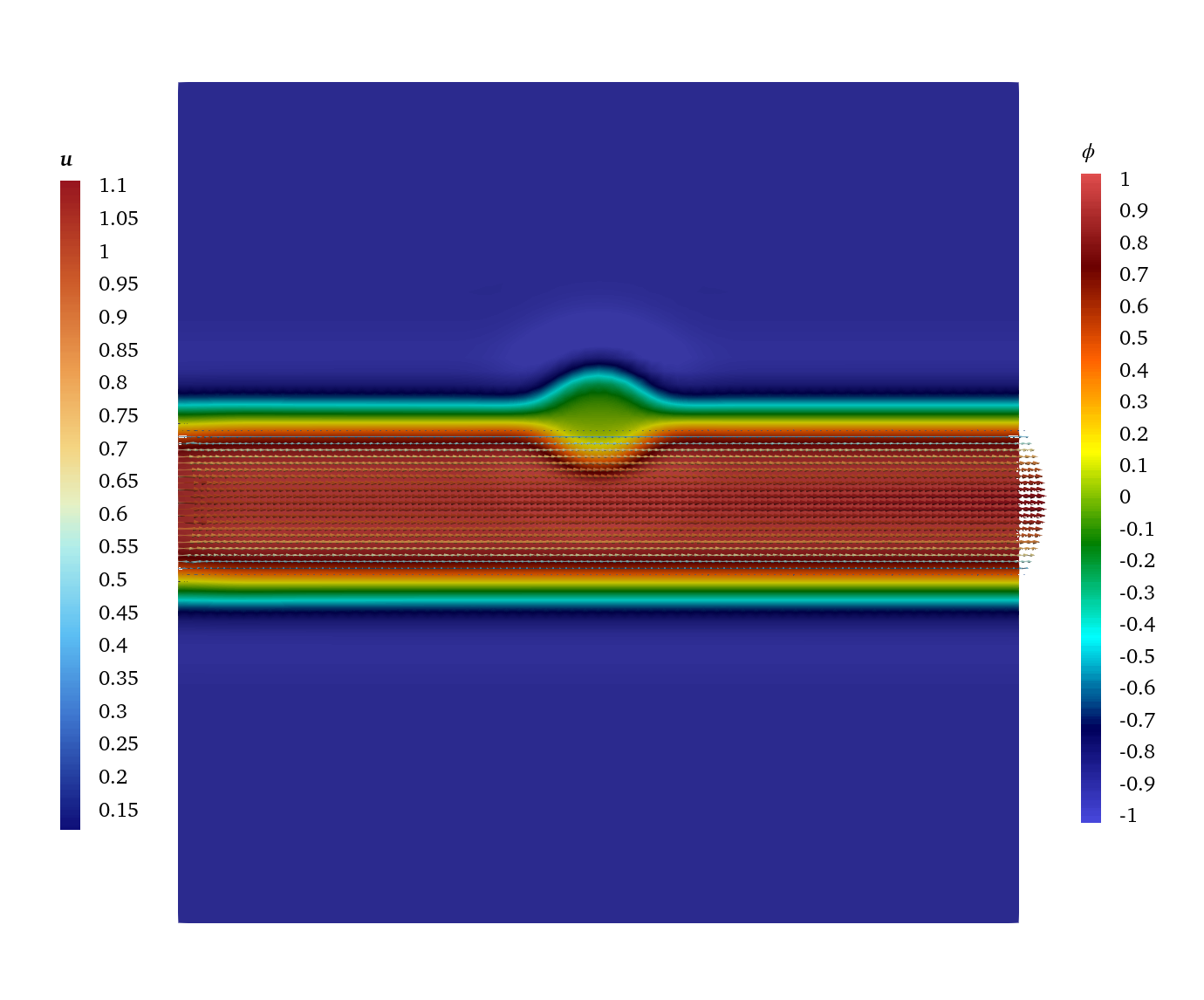}
	}
	\subfigure[$t=3.0$]{
		\includegraphics[scale=0.10]{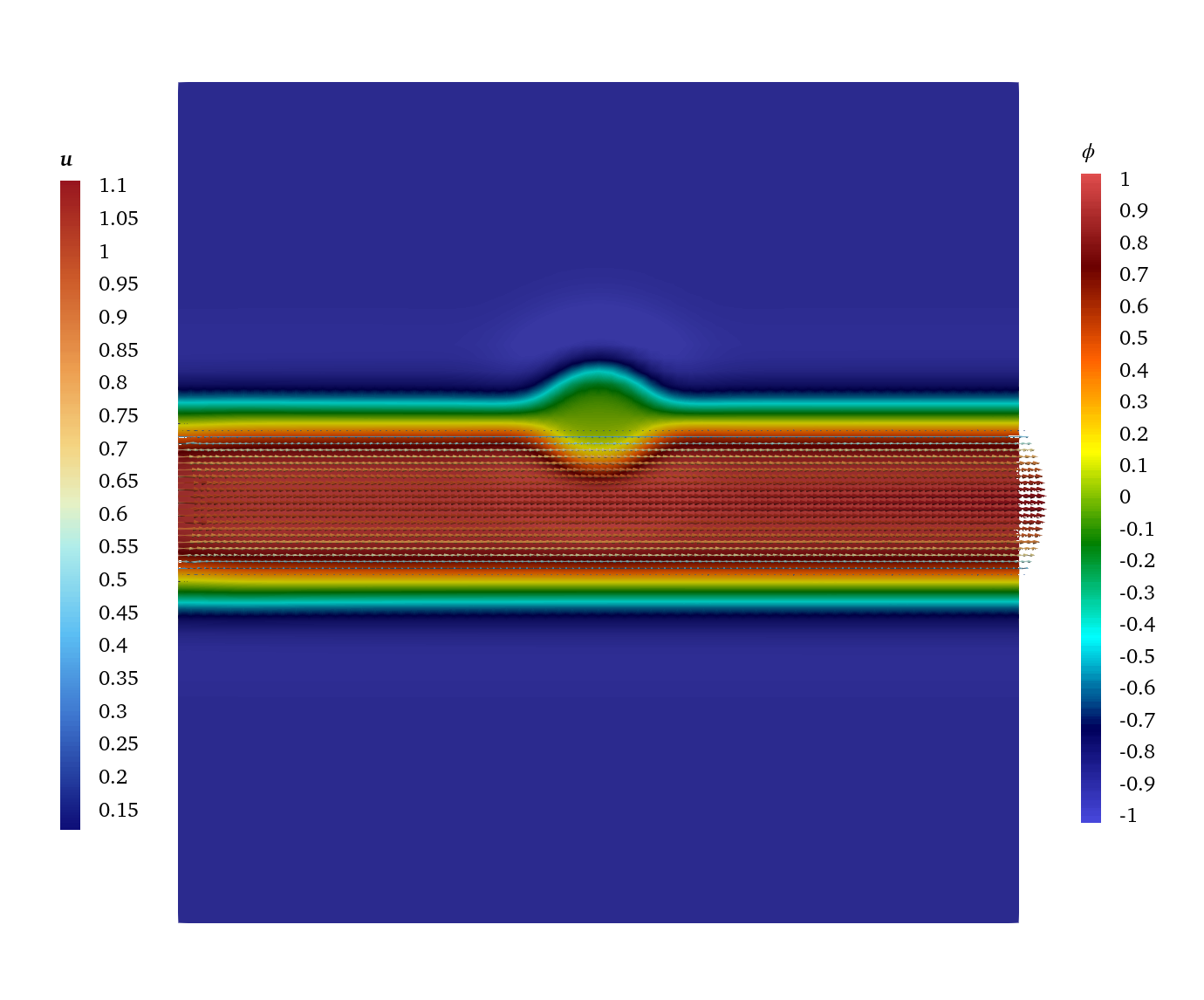}
	}
	
	\caption{Blood vessel snapshots at different times.}
	\label{figure7}
\end{figure}

\begin{figure}[H]
	\centering
	\subfigure[$t=1.0$]{
		\includegraphics[scale=0.10]{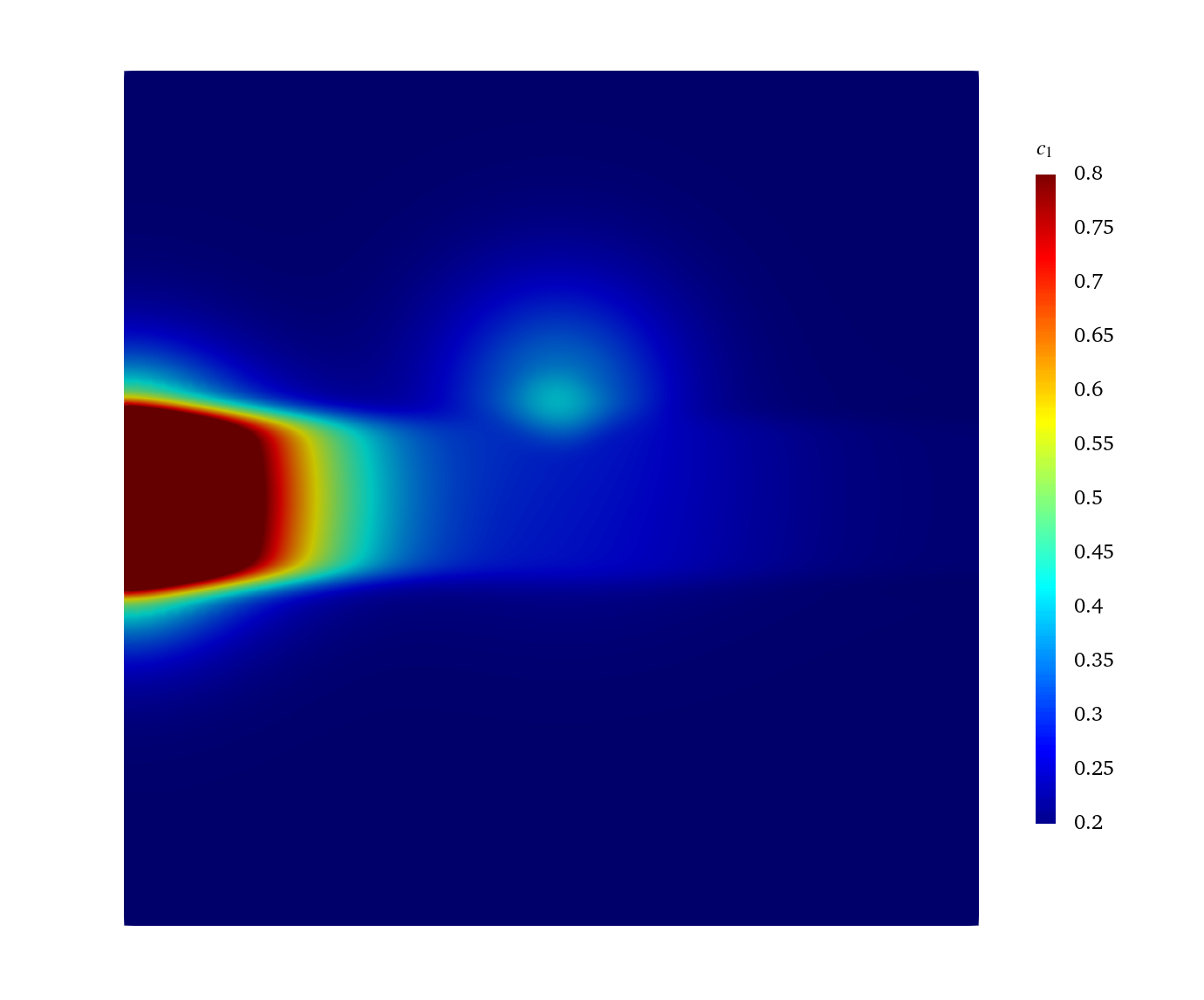}
	}
	\subfigure[$t=2.0$]{
		\includegraphics[scale=0.10]{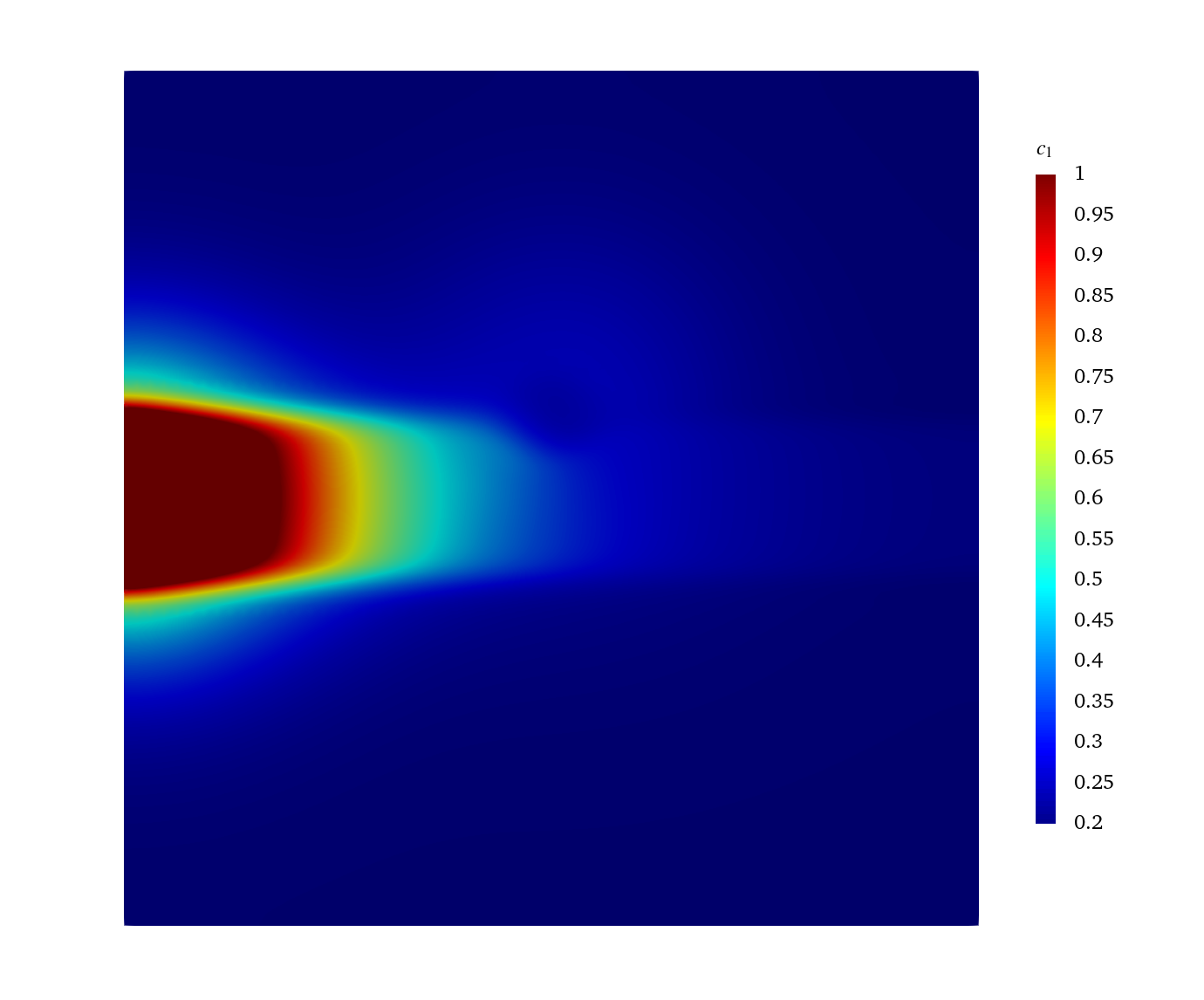}
	}
	\subfigure[$t=3.0$]{
		\includegraphics[scale=0.10]{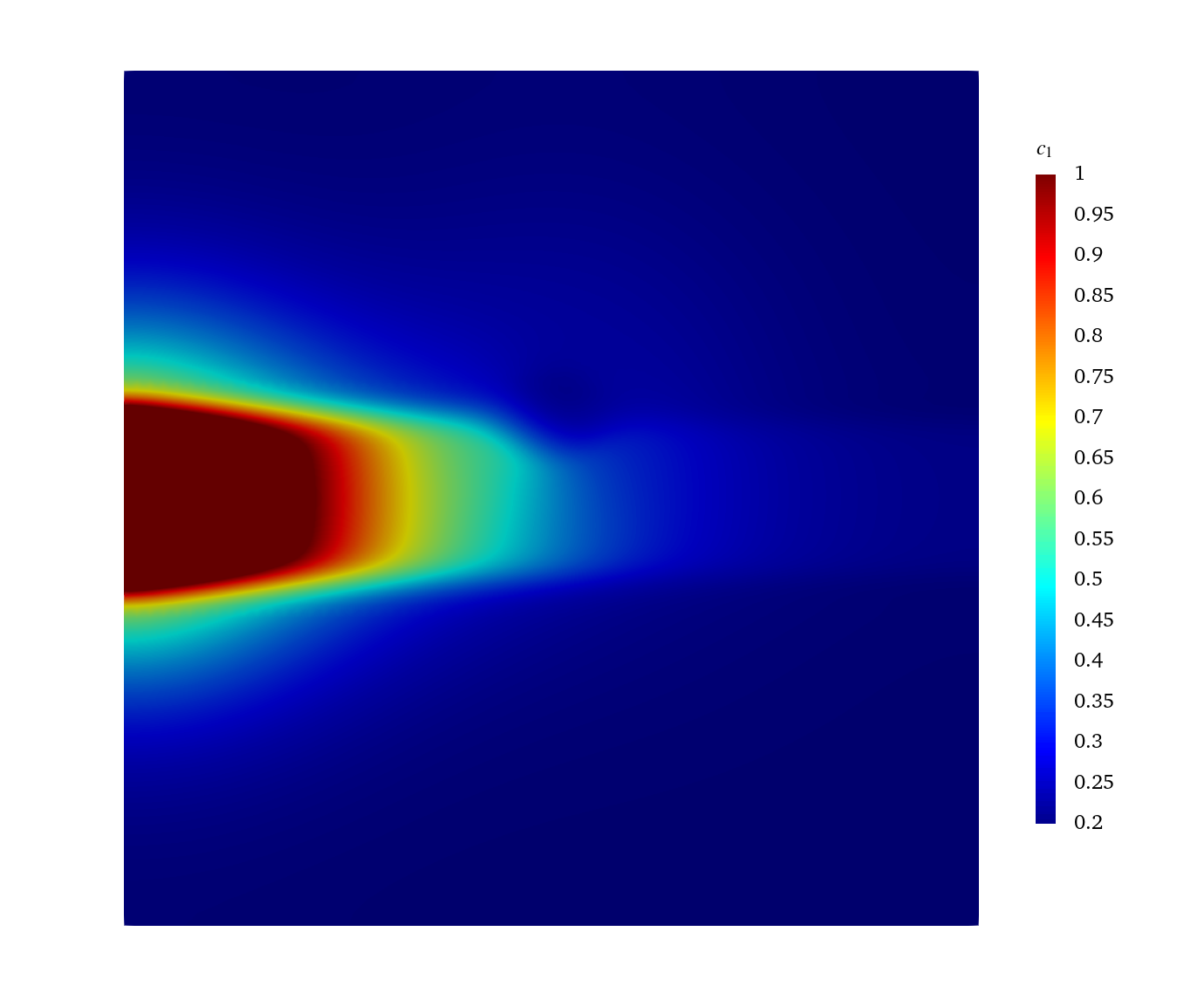}
	}
	
	\caption{Snapshots of glucose $(Y_1)$ at different times.}
	\label{figure8}
\end{figure}

\begin{figure}[H]
	\centering
	\subfigure[$t=1.0$]{
		\includegraphics[scale=0.10]{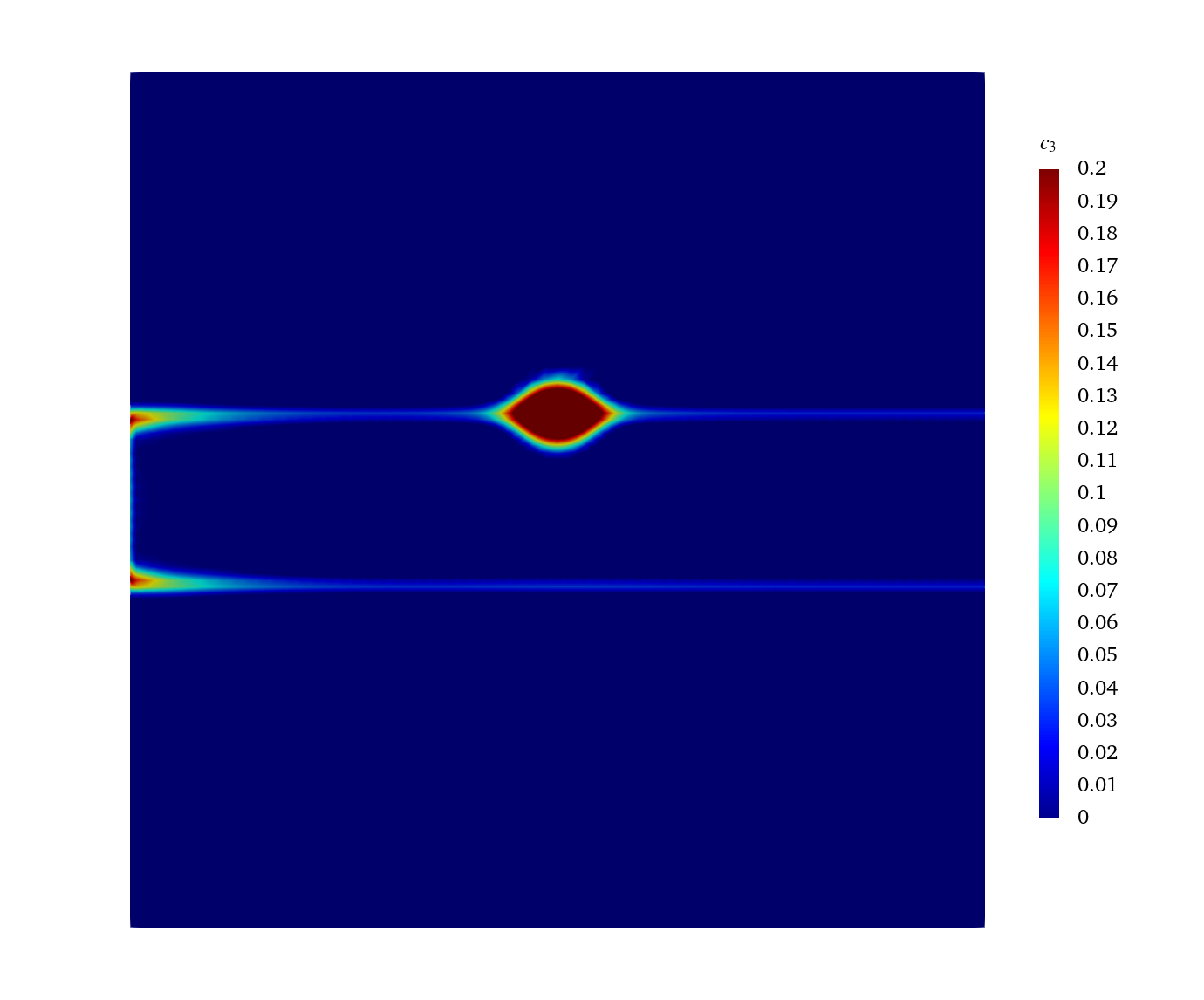}
	}
	\subfigure[$t=2.0$]{
		\includegraphics[scale=0.10]{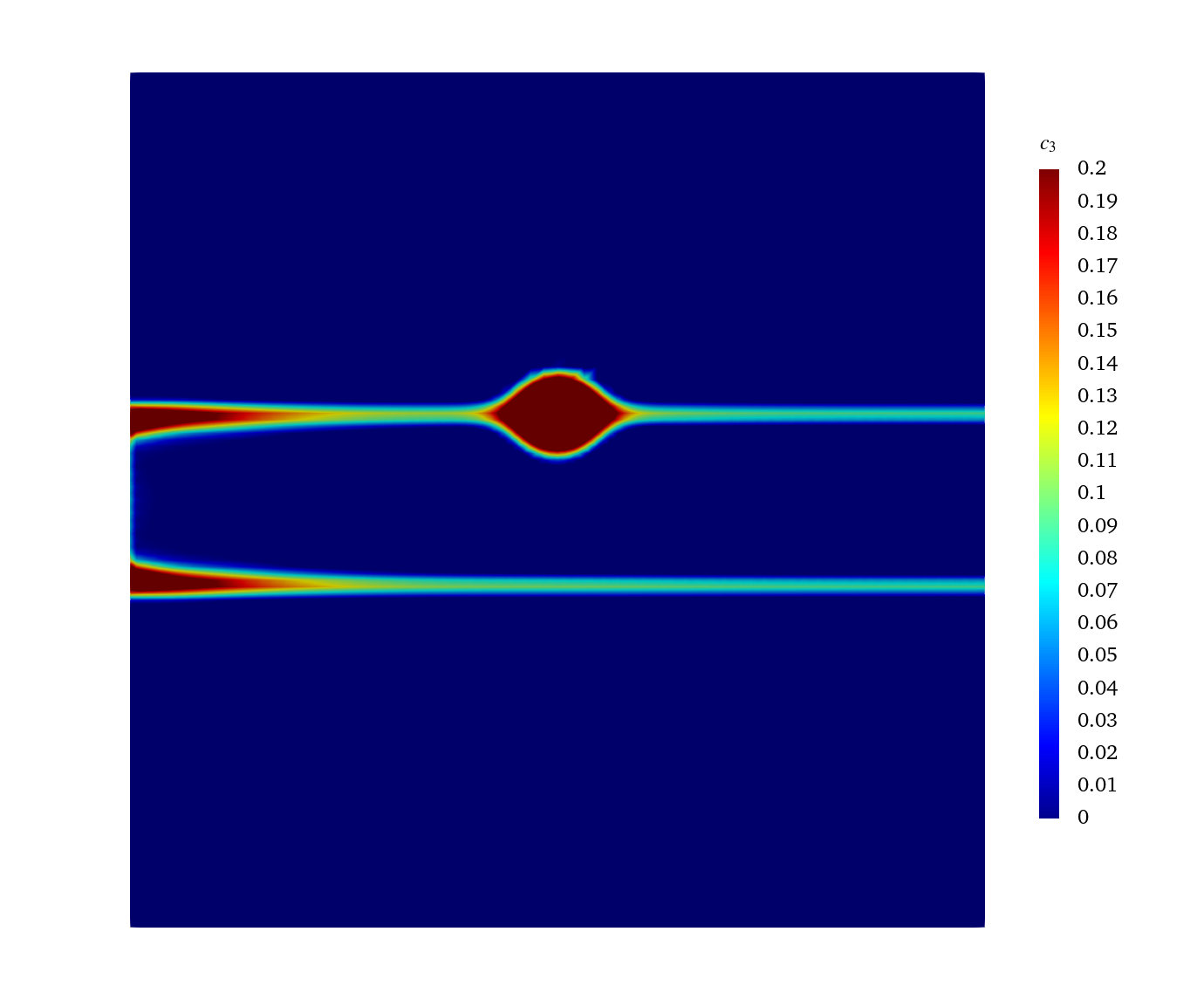}
	}
	\subfigure[$t=3.0$]{
		\includegraphics[scale=0.10]{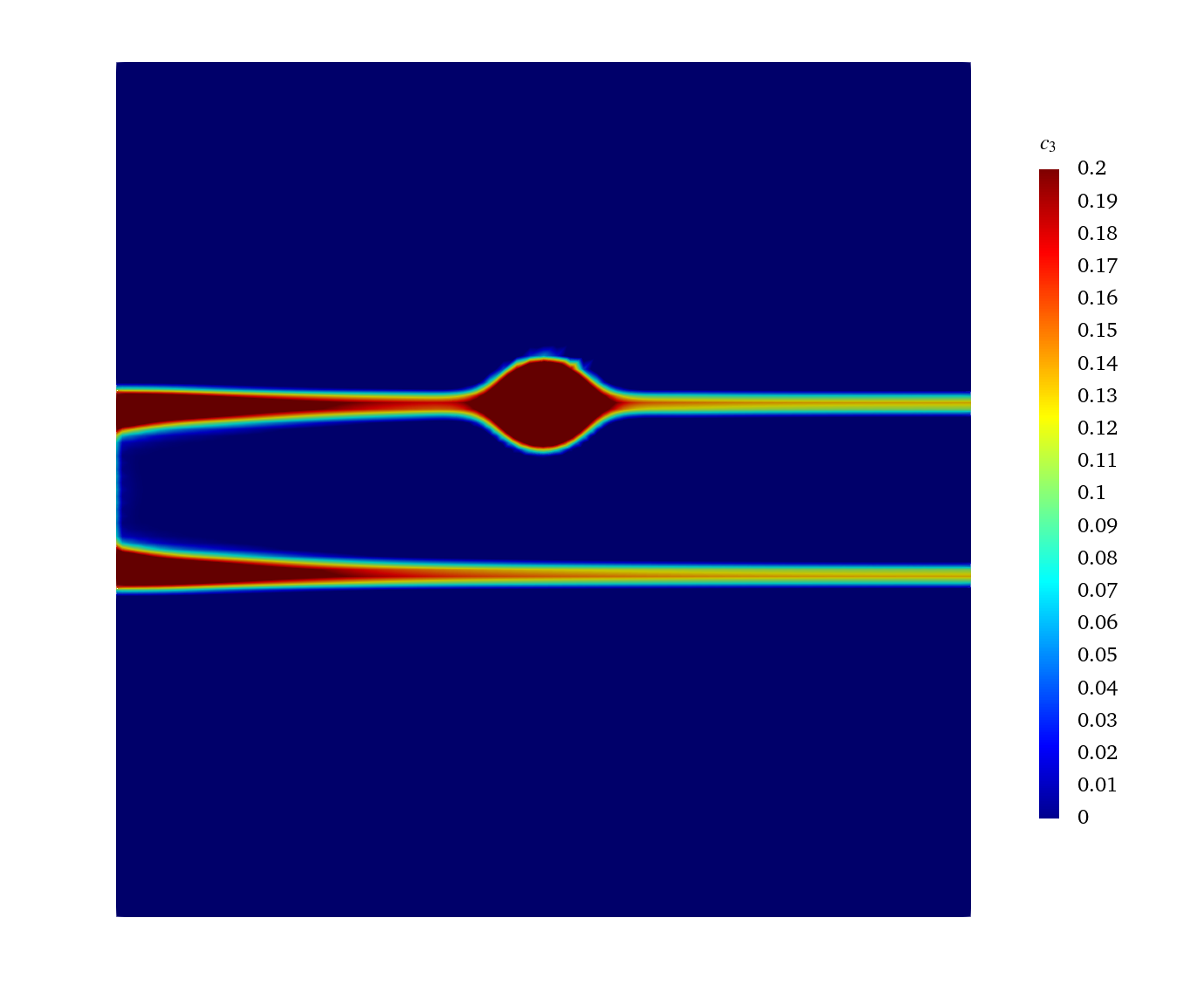}
	}
	
	\caption{Snapshots of AGEs $(Y_3)$ at different times.}
	\label{figure9}
\end{figure}

Figure \ref{figure7} shows the temporal evolution of a microaneurysm within a blood vessel. We can observe that the initial stage of aneurysm formation is visible as a localized bulge, which gradually evolves into a larger size. The temporal evolution profiles of glucose $(Y_1)$ and AGEs $(Y_3)$ are respectively presented in Figures \ref{figure8} and \ref{figure9}. The continuous consumption of glucose to produce AGEs correlates with the progressive aneurysm growth observed in Figure \ref{figure7}. Consequently, as the mixing energy density function $\lambda(c_3)$ increases, the interface tends to expand laterally to reduce the phase mixing energy of the system. In addition, it can be seen that as AGEs accumulate on the vascular wall near the inlet, glucose diffuses across the vascular wall into the interstitial fluid.


\subsubsection{Bifurcated microvessel} 
\begin{figure}[H]
	\centering
	\subfigure[$t=1.0$]{
		\includegraphics[scale=0.11]{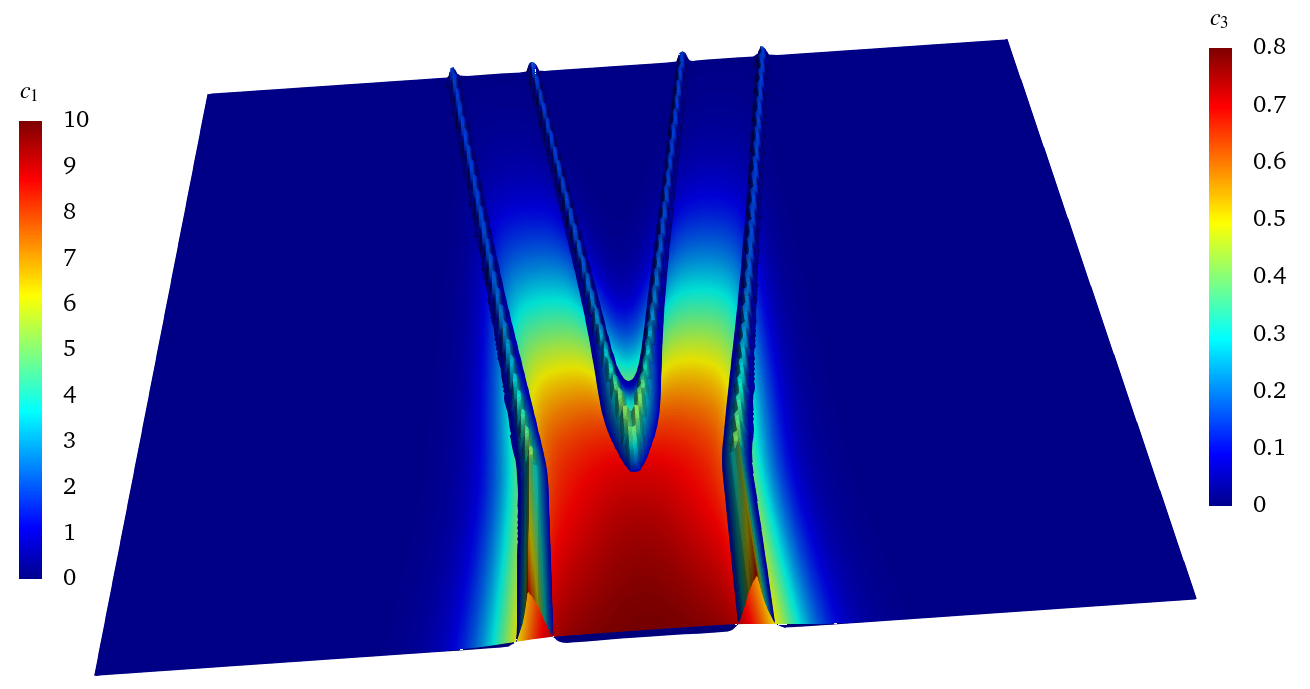}
	}
	\subfigure[$t=2.0$]{
		\includegraphics[scale=0.22]{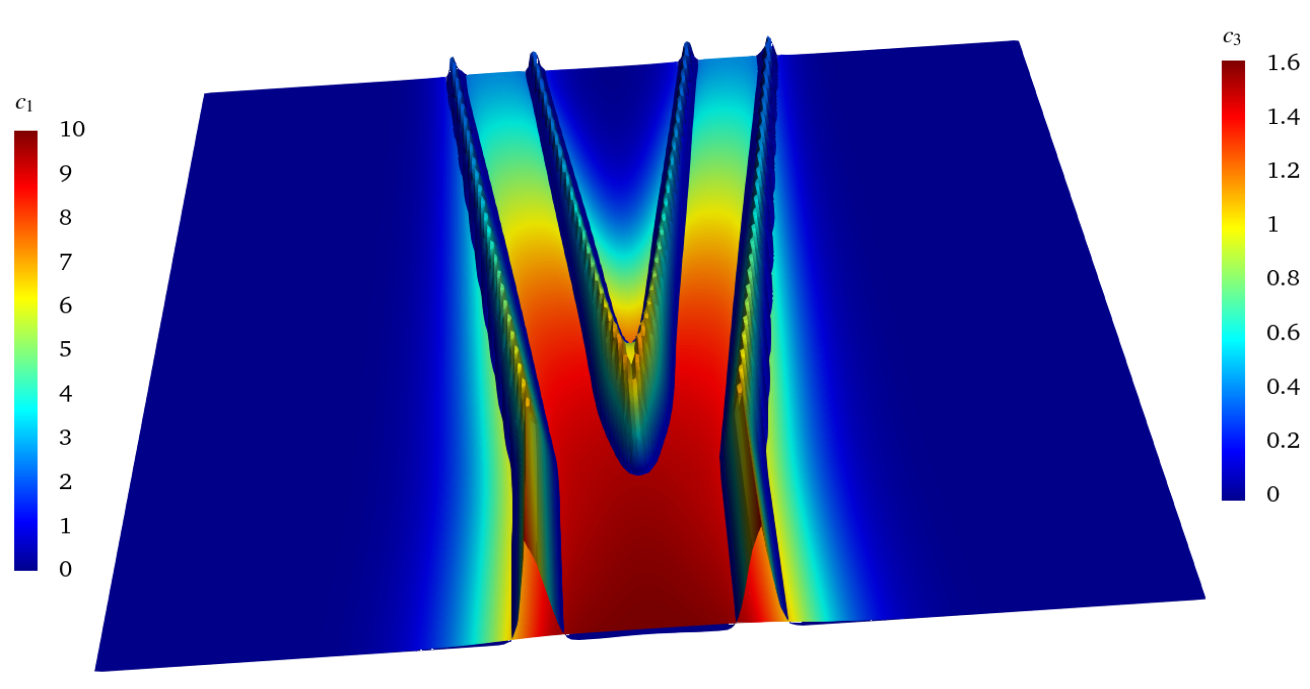}
	}
	\subfigure[$t=3.0$]{
		\includegraphics[scale=0.22]{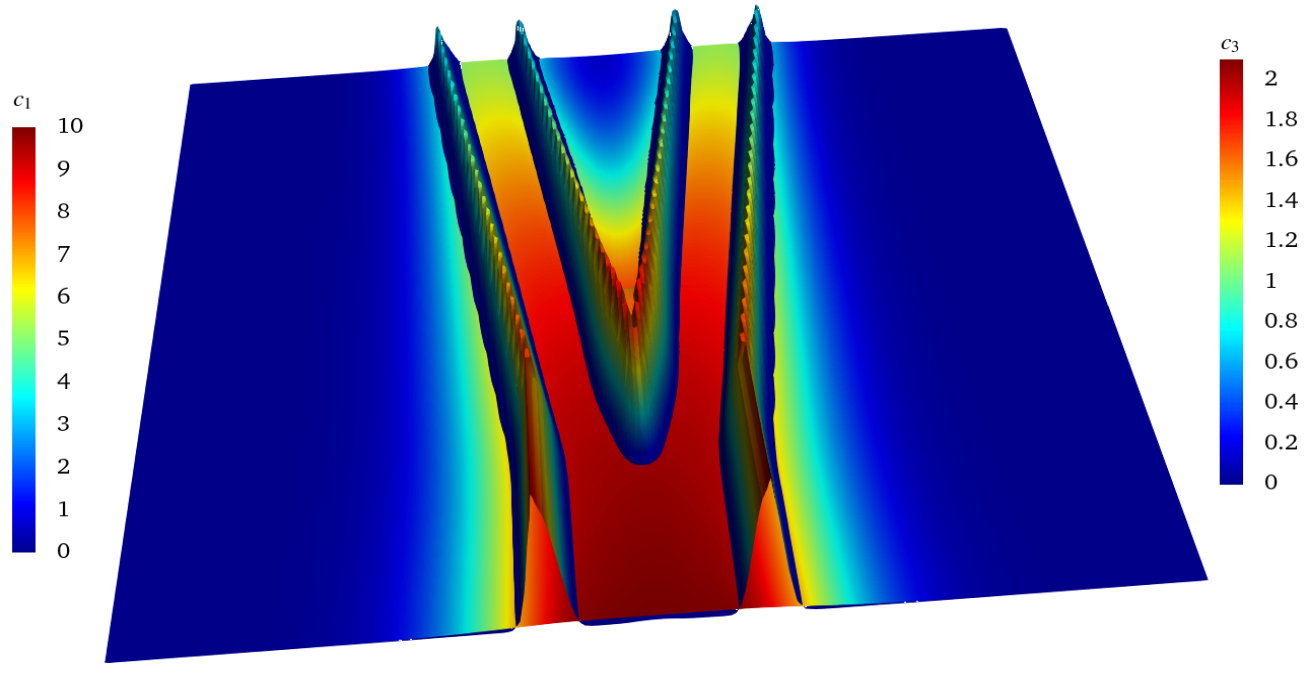}
	}
\caption{Three-dimensional visualization distribution of glucose $(Y_1)$ and AGEs $(Y_3)$ at different times. The base displays the distribution of glucose and the vertical protrusions represent
the concentration of AGEs. The vascular structure is symmetrical, with a bifurcation angle of $14^\circ$, and the bifurcation point is located at $(0.5,1)$.}
	\label{figure10}
\end{figure}

We simulate the Y-shaped bifurcated vessel without changing the basic parameter settings to investigate the spatial distribution of AGEs. Figure \ref{figure10} presents the distribution of glucose and AGEs in a bifurcated vessel structure. The influx of substantial glucose from the inlet into the vasculature, and its reaction with structural proteins on the vessel wall, results in progressive accumulation of AGEs.

\begin{figure}[H]
	\centering
	\includegraphics[scale=0.58]{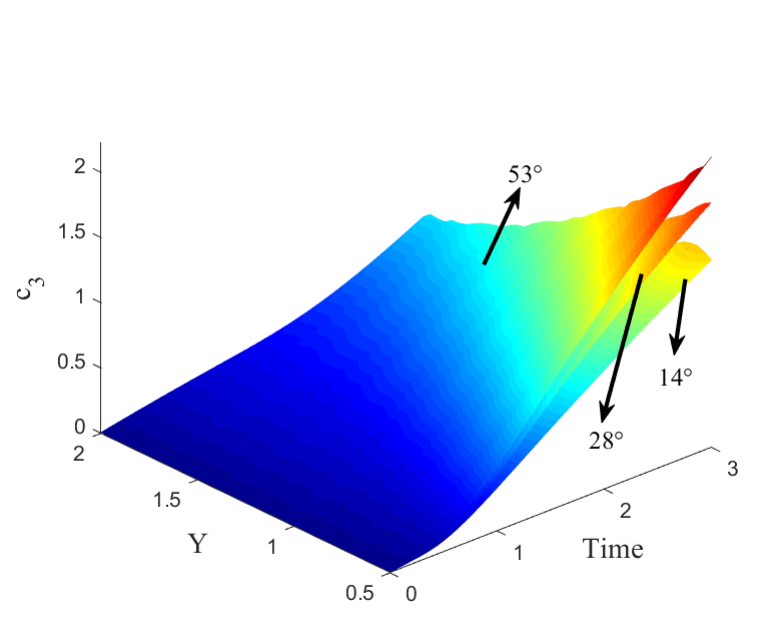}
	\caption{The distribution of AGEs on the vascular wall at the bifurcation point under different bifurcation angles.}
	\label{figure11}
\end{figure}

Near vascular bifurcations, blood flow is prone to flow separation, generating vortex structures or stagnation zones. These hemodynamic disturbances appear to facilitate the accumulation of AGEs, which can lead to the formation of microaneurysms. We simulate vessels with larger bifurcation angles and extract the concentration of AGEs on the vascular wall at the bifurcation point for plotting. As shown in Figure \ref{figure11}, larger bifurcation angles are associated with higher concentrations of AGEs near the bifurcation. This phenomenon may be attributed to the hemodynamic environment, where larger bifurcation angles induce more pronounced flow separation and recirculation zones, leading to prolonged residence time, which in turn facilitates the glycation reaction.


\section{Conclusions and remarks}
\label{section6}
In this paper, we develop a thermodynamically consistent phase-field model to describe mass transport, interfacial reactions, and deformation in fluid environments. The proposed model satisfies both mass conservation and energy dissipation laws. A structure-preserving numerical scheme is designed and a rigorous error analysis is performed for a simplified case. Several numerical examples verify the effectiveness of the model and the numerical method. Numerical simulations of microaneurysm formation in straight and bifurcated vessels reveal critical insights into the role of hemodynamics and AGEs accumulation.

It should be noted that the model framework proposed in this paper is not only applicable to vascular pathologies, but can also be easily adapted to corrosion processes and flows containing surfactants. In addition, this model  can be combined with vesicle models \cite{du2009energetic} to describe vesicle motion under vascular deformation. The development of efficient numerical methods for the model is critically important, as it can significantly enhance computational scalability for large-scale 3-D simulations.

\bibliographystyle{siam}
\bibliography{Ref}
\end{document}

%% file: ex_shared.tex
\usepackage{lipsum}
\usepackage{amsfonts}
\usepackage{graphicx}
\usepackage{epstopdf}
\usepackage{algorithmic}
\ifpdf
  \DeclareGraphicsExtensions{.eps,.pdf,.png,.jpg}
\else
  \DeclareGraphicsExtensions{.eps}
\fi

\newsiamremark{remark}{Remark}
\newsiamremark{hypothesis}{Hypothesis}
\newsiamremark{assumption}{Assumption}
\newsiamthm{claim}{Claim}

\headers{Phase-field model with interfacial reaction and deformation}{Zhaoyang Wang, Huaxiong Huang, Ping Lin and Shixin Xu}

\title{A thermodynamically consistent phase-field model for mass transport with interfacial reaction and deformation\thanks{Submitted to the editors \today.
\funding{Z. Wang is partially supported by the China Postdoctoral Science Foundation under grant 2024M760239. H. Huang is partially supported by  the National Natural Science Foundation of China 12231004. P. Lin is partially supported by  the National Natural Science Foundation of China 12371388, 11861131004 and 11771040. S. Xu is partially supported by the National Natural Science Foundation of China 12071190.}}}

\author{Zhaoyang Wang\thanks{School of Mathematical Sciences, Laboratory of Mathematics and Complex Systems, MOE, Beijing Normal University, Beijing 100875, China; Research Center for Mathematics, Advanced Institute of Natural Sciences, Beijing Normal University, Zhuhai, Guangdong 519087, China.
(\email{zhaoyang584520@163.com}).}
\and Huaxiong Huang\thanks{Corresponding author. Zu Chongzhi Center for Mathematics and Computational Sciences (CMCS), Duke Kunshan University, 8 Duke Ave., Kunshan, Jiangsu, 215316, China; Guandong Key Lab for Data Science and Technology, BNU-HKBU United International College, Zhuhai, Guangdong, China; Department of Mathematics and Statistics, York University, Toronto, Ontario, Canada M3J 1P3. (\email{hh321@duke.edu}).}
\and Ping Lin\thanks{Division of Mathematics, University of Dundee, Dundee DD1 4HN, United Kingdom (p.lin@dundee.ac.uk).}
\and Shixin Xu\thanks{Zu Chongzhi Center for Mathematics and Computational Sciences (CMCS), Duke Kunshan University, 8 Duke Ave., Kunshan, Jiangsu, 215316, China. (\email{shixin.xu@dukekunshan.edu.cn}).}}

\usepackage{amsopn}

%% file: Phasefield.bbl
\begin{thebibliography}{10}

\bibitem{bernardo2020recent}
{\sc G.~Bernardo, T.~Ara{\'u}jo, T.~da~Silva~Lopes, J.~Sousa, and A.~Mendes},
  {\em Recent advances in membrane technologies for hydrogen purification},
  International Journal of Hydrogen Energy, 45 (2020), pp.~7313--7338.

\bibitem{claesson2021permeability}
{\sc L.~Claesson-Welsh, E.~Dejana, and D.~M. McDonald}, {\em Permeability of
  the endothelial barrier: identifying and reconciling controversies}, Trends
  in molecular medicine, 27 (2021), pp.~314--331.

\bibitem{de2013non}
{\sc S.~R. De~Groot and P.~Mazur}, {\em Non-equilibrium thermodynamics},
  Courier Corporation, 2013.

\bibitem{du2009energetic}
{\sc Q.~Du, C.~Liu, R.~Ryham, and X.~Wang}, {\em Energetic variational
  approaches in modeling vesicle and fluid interactions}, Physica D: Nonlinear
  Phenomena, 238 (2009), pp.~923--930.

\bibitem{du2004phase}
{\sc Q.~Du, C.~Liu, and X.~Wang}, {\em A phase field approach in the numerical
  study of the elastic bending energy for vesicle membranes}, Journal of
  Computational Physics, 198 (2004), pp.~450--468.

\bibitem{engblom2013diffuse}
{\sc S.~Engblom, M.~Do-Quang, G.~Amberg, and A.-K. Tornberg}, {\em On diffuse
  interface modeling and simulation of surfactants in two-phase fluid flow},
  Communications in Computational Physics, 14 (2013), pp.~879--915.

\bibitem{girault2012finite}
{\sc V.~Girault and P.-A. Raviart}, {\em Finite element methods for
  Navier-Stokes equations: theory and algorithms}, vol.~5, Springer Science \&
  Business Media, 2012.

\bibitem{guo2015thermodynamically}
{\sc Z.~Guo and P.~Lin}, {\em A thermodynamically consistent phase-field model
  for two-phase flows with thermocapillary effects}, Journal of Fluid
  Mechanics, 766 (2015), pp.~226--271.

\bibitem{hawkins2012numerical}
{\sc A.~Hawkins-Daarud, K.~G. van~der Zee, and J.~Tinsley~Oden}, {\em Numerical
  simulation of a thermodynamically consistent four-species tumor growth
  model}, International journal for numerical methods in biomedical
  engineering, 28 (2012), pp.~3--24.

\bibitem{hirai2007retinopathy}
{\sc F.~E. Hirai, S.~E. Moss, M.~D. Knudtson, B.~E. Klein, and R.~Klein}, {\em
  Retinopathy and survival in a population without diabetes: The beaver dam eye
  study}, American journal of epidemiology, 166 (2007), pp.~724--730.

\bibitem{hudson2018targeting}
{\sc B.~I. Hudson and M.~E. Lippman}, {\em Targeting rage signaling in
  inflammatory disease}, Annual review of medicine, 69 (2018), pp.~349--364.

\bibitem{jin2022corrosion}
{\sc S.~Jin, H.~Tian, Q.~Wang, T.~Ren, P.~Liu, and Z.~Peng}, {\em Corrosion
  reaction kinetics and high-temperature corrosion testing of contact element
  strips in ultra-high voltage bushing based on the phase-field method}, IET
  generation, transmission \& distribution, 16 (2022), pp.~2947--2958.

\bibitem{kou2023thermodynamically}
{\sc J.~Kou, A.~Salama, and X.~Wang}, {\em Thermodynamically consistent
  phase-field modelling of activated solute transport in binary solvent
  fluids}, Journal of Fluid Mechanics, 955 (2023), p.~A41.

\bibitem{laradji1992effect}
{\sc M.~Laradji, H.~Guo, M.~Grant, and M.~J. Zuckermann}, {\em The effect of
  surfactants on the dynamics of phase separation}, Journal of Physics:
  Condensed Matter, 4 (1992), p.~6715.

\bibitem{meng2025convergence}
{\sc X.~Meng, Y.~Qin, and G.~Hu}, {\em The convergence analysis of a class of
  stabilized semi-implicit isogeometric methods for the cahn-hilliard
  equation}, Journal of Scientific Computing, 102 (2025), pp.~1--35.

\bibitem{oden2010general}
{\sc J.~T. Oden, A.~Hawkins, and S.~Prudhomme}, {\em General diffuse-interface
  theories and an approach to predictive tumor growth modeling}, Mathematical
  Models and Methods in Applied Sciences, 20 (2010), pp.~477--517.

\bibitem{osher2004level}
{\sc S.~Osher, R.~Fedkiw, and K.~Piechor}, {\em Level set methods and dynamic
  implicit surfaces}, Appl. Mech. Rev., 57 (2004), pp.~B15--B15.

\bibitem{peskin2002immersed}
{\sc C.~S. Peskin}, {\em The immersed boundary method}, Acta numerica, 11
  (2002), pp.~479--517.

\bibitem{qin2022phase}
{\sc Y.~Qin, H.~Huang, Y.~Zhu, C.~Liu, and S.~Xu}, {\em A phase field model for
  mass transport with semi-permeable interfaces}, Journal of Computational
  Physics, 464 (2022), p.~111334.

\bibitem{quarteroni2008numerical}
{\sc A.~Quarteroni and A.~Valli}, {\em Numerical approximation of partial
  differential equations}, vol.~23, Springer Science \& Business Media, 2008.

\bibitem{shen2018scalar}
{\sc J.~Shen, J.~Xu, and J.~Yang}, {\em The scalar auxiliary variable (sav)
  approach for gradient flows}, Journal of Computational Physics, 353 (2018),
  pp.~407--416.

\bibitem{shen2019new}
\leavevmode\vrule height 2pt depth -1.6pt width 23pt, {\em A new class of
  efficient and robust energy stable schemes for gradient flows}, SIAM Review,
  61 (2019), pp.~474--506.

\bibitem{shen2010numerical}
{\sc J.~Shen and X.~Yang}, {\em Numerical approximations of allen-cahn and
  cahn-hilliard equations}, Discrete Contin. Dyn. Syst, 28 (2010),
  pp.~1669--1691.

\bibitem{shen2015decoupled}
\leavevmode\vrule height 2pt depth -1.6pt width 23pt, {\em Decoupled, energy
  stable schemes for phase-field models of two-phase incompressible flows},
  SIAM Journal on Numerical Analysis, 53 (2015), pp.~279--296.

\bibitem{shen2022energy}
{\sc L.~Shen, Z.~Xu, P.~Lin, H.~Huang, and S.~Xu}, {\em An energy stable \(
  {C}^0 \) finite element scheme for a phase-field model of vesicle motion and
  deformation}, SIAM Journal on Scientific Computing, 44 (2022),
  pp.~B122--B145.

\bibitem{strathmann1981membrane}
{\sc H.~Strathmann}, {\em Membrane separation processes}, Journal of membrane
  science, 9 (1981), pp.~121--189.

\bibitem{temam2024navier}
{\sc R.~Temam}, {\em Navier--Stokes equations: theory and numerical analysis},
  vol.~343, American Mathematical Society, 2024.

\bibitem{tryggvason2001front}
{\sc G.~Tryggvason, B.~Bunner, A.~Esmaeeli, D.~Juric, N.~Al-Rawahi, W.~Tauber,
  J.~Han, S.~Nas, and Y.-J. Jan}, {\em A front-tracking method for the
  computations of multiphase flow}, Journal of computational physics, 169
  (2001), pp.~708--759.

\bibitem{unverdi1992front}
{\sc S.~O. Unverdi and G.~Tryggvason}, {\em A front-tracking method for
  viscous, incompressible, multi-fluid flows}, Journal of computational
  physics, 100 (1992), pp.~25--37.

\bibitem{wang2003level}
{\sc M.~Y. Wang, X.~Wang, and D.~Guo}, {\em A level set method for structural
  topology optimization}, Computer methods in applied mechanics and
  engineering, 192 (2003), pp.~227--246.

\bibitem{wang2025stability}
{\sc Z.~Wang, P.~Lin, and J.~Yang}, {\em Stability and error analysis of
  structure-preserving schemes for a diffuse-interface tumor growth model},
  SIAM Journal on Scientific Computing, 47 (2025), pp.~B59--B86.

\bibitem{wang2023fast}
{\sc Z.~Wang, P.~Lin, and L.~Zhang}, {\em A fast front-tracking approach and
  its analysis for a temporal multiscale flow problem with a fractional order
  boundary growth}, SIAM Journal on Scientific Computing, 45 (2023),
  pp.~B646--B672.

\bibitem{wautier2015advanced}
{\sc M.~Wautier and J.~Wautier}, {\em Advanced glycation end products and
  retinal vascular lesions in diabetes mellitus}, Austin J. Endocrinol.
  Diabetes, 2 (2015), p.~1034.

\bibitem{xu2018osmosis}
{\sc S.~Xu, B.~Eisenberg, Z.~Song, and H.~Huang}, {\em Osmosis through a
  semi-permeable membrane: a consistent approach to interactions}, arXiv
  preprint arXiv:1806.00646,  (2018).

\bibitem{yang2017linear}
{\sc X.~Yang and L.~Ju}, {\em Linear and unconditionally energy stable schemes
  for the binary fluid--surfactant phase field model}, Computer Methods in
  Applied Mechanics and Engineering, 318 (2017), pp.~1005--1029.

\bibitem{yang2016mathematical}
{\sc Y.~Yang, W.~J{\"a}ger, M.~Neuss-Radu, and T.~Richter}, {\em Mathematical
  modeling and simulation of the evolution of plaques in blood vessels},
  Journal of mathematical biology, 72 (2016), pp.~973--996.

\bibitem{zhu2019thermodynamically}
{\sc G.~Zhu, J.~Kou, B.~Yao, Y.-s. Wu, J.~Yao, and S.~Sun}, {\em
  Thermodynamically consistent modelling of two-phase flows with moving contact
  line and soluble surfactants}, Journal of Fluid Mechanics, 879 (2019),
  pp.~327--359.

\end{thebibliography}
